\documentclass[11pt]{article}

\usepackage{amsmath}
\usepackage{graphicx}
\usepackage{latexsym}
\usepackage{amssymb}
\usepackage{amsbsy}
\usepackage[]{authblk}
\usepackage{array}
\usepackage{color}
\usepackage{lineno}
\usepackage{ulem}
\usepackage{multirow}
\usepackage{rotating}
\DeclareSymbolFont{msbm}{U}{msb}{m}{n}
\DeclareMathSymbol{\C}{\mathalpha}{msbm}{'103}
\DeclareMathSymbol{\R}{\mathalpha}{msbm}{'122}
\DeclareMathSymbol{\Z}{\mathalpha}{msbm}{'132}
\DeclareMathSymbol{\N}{\mathalpha}{msbm}{'116}

\newtheorem{remark}{Remark}

\newtheorem{theorem}{Theorem}

\newcommand{\Sum}{\displaystyle\sum}

\def\P{{\cal P}}
\def\R{{\cal R}}
\def\T{{\cal T}}

\def\be{\begin{equation}}
\def\ee{\end{equation}}
\def\bea{\begin{eqnarray}}
\def\ba{\begin{array}{l}\displaystyle}
\def\eea{\end{eqnarray}}
\def\ea{\end{array}}
\def\E{{\cal E}}
\def\ca{\\[+0.3cm]\displaystyle}

\newfont{\numerikEleven}{ecrm1000}
\newfont{\numerikTen}{cmss10}
\newfont{\numerikNine}{cmss9}
\newfont{\numerikEight}{cmss8}
\newfont{\numerikSeven}{cmss7}
\begin{document}

% Put line numbers to help reviewing
\linenumbers

%==========================================================================================
% FIRST PAGE
%------------------------------------------
% TITLE
\title{ \Large
Towards an ultra efficient kinetic scheme\\ Part II: The high order case \thanks{This work was supported by
the ANR Blanc project BOOST}}
%------------------------------------------
% AUTHORS
%
\author[1]{Giacomo Dimarco\footnote
{Corresponding author address: Universit\'{e} de Toulouse; UPS, INSA, UT1, UTM ; CNRS, UMR 5219
Institut de Math\'{e}matiques de Toulouse ; F-31062 Toulouse,
France. \\
\emph{E-mail}: giacomo.dimarco@math.univ-toulouse.fr, raphael.loubere@math.univ-toulouse.fr}}
\author[1]{Rapha\"{e}l Loubere}
\affil[1]{Universit\'{e} de Toulouse; UPS, INSA, UT1, UTM ; CNRS, UMR 5219
Institut de Math\'{e}matiques de Toulouse ; F-31062 Toulouse,
France.}
%
% DATE
%
\date{\today}
\maketitle
%-------------------------------------------
% ABSTRACT
%
\begin{abstract}
In a recent paper we presented a new ultra efficient numerical method for
solving kinetic equations of the Boltzmann type \cite{dimarco1}.
The key idea, on which the method relies, is to solve the
collision part on a grid and then to solve exactly the transport
part by following the characteristics backward in time. On the contrary to classical semi-Lagrangian
methods one does not need to reconstruct the distribution
function at each time step. This allows to tremendously reduce the
computational cost and to perform efficient
numerical simulations of kinetic equations up to the six dimensional case
without parallelization.
However, the main drawback of the method developed was the loss of spatial accuracy close to the fluid limit.
In the present work, we modify the scheme in such a way that it is able to preserve the high order spatial
accuracy for extremely rarefied and fluid regimes. In particular, in the fluid limit, the method
automatically degenerates into a high order method for the compressible Euler equations.
Numerical examples are presented which validate the method, show the higher accuracy with respect to
the previous approach and measure its efficiency with respect to well known schemes (Direct Simulation Monte Carlo,
Finite Volume, MUSCL, WENO).
\end{abstract}
%------------------------------------------
% KEY WORDS
%

%%%%%%%%%%%%%%%%%%%%%%%%%%%%%%%%%%%%%%%%%%%%%%%%%%%%%%%%%%%%%%%%%%%%%%%%%%%%%%%%%%%%%%%%%%%%%%
%
% S E C T I O N
%
\section{Introduction}
The kinetic equations provide a description of non equilibrium gases and
more generally of particle systems \cite{bird, cercignani}. The distribution function,
which describes the evolution of the system, depends, in the most general case, on seven independent variables:
the time, the physical and the velocity space.
It turns out that the numerical simulation of these kind of equations with deterministic
techniques presents several drawbacks due to the large dimension of
the problem. On the other side of the spectrum of the numerical techniques used to approximate kinetic equations,
there are the probabilistic methods \cite{bird, Cf, CPima, Nanbu80}. These methods and in particular Monte Carlo
methods (DSMC) are extensively used due to their very low computational cost especially in the multidimensional cases,
compared to finite volume, finite difference or spectral methods \cite{Filbet2, Filbet, Mieussens, Mouhot, Pal, Pieraccini, Plat, Plat1}.
However, DSMC furnishes only poorly accurate and fluctuating solutions which cannot be easily ameliorated. This is especially
true in non stationary situations in which time averages techniques turn to be useless.

Many different works have been dedicated to reduce
some of the disadvantages of Monte Carlo methods. We quote \cite{Cf}
for an overview on efficient and low variance Monte Carlo methods.
Let us remind to the works of Homolle and Hadjiconstantinou
\cite{Hadji} and \cite{Hadji1} and of Dimarco, Pareschi and Degond
\cite{dimarco2, dimarco3, dimarco4, dimarcoN} for some applications of variance reduction techniques to kinetic
equations in transitional and general regimes. We recall also the works of Boyd and
Burt \cite{Boyd} and of Pullin \cite{Pullin78} who developed a low
diffusion particle methods for simulating compressible inviscid
flows.

In this work, we continue the development of a new deterministic
ultra fast method which permits to solve kinetic equations of the Boltzmann type \cite{dimarco1}. The scheme is
based on the classical discrete velocity models (DVM) approach
\cite{bobylev, Pal, Plat, Plat1}. The DVM models are obtained by discretizing the velocity
space into a set of fixed discrete velocities \cite{bobylev,
Mieussens, Pal, Pal1}. As a result of this discretization, the
original kinetic equation is then represented as a set of linear
transport equations plus a source term. The source term describes the collisions or the interactions between the particles and couples all the
equations of the resulting system. In order to solve the transport part of the DVM model, many different techniques can be employed like finite difference, finite volume or fast methods \cite{Besse, Filbet, Filbet2, Mouhot, Mieussens, Gu}. One of the most common strategies for solving this kind of problems is the semi Lagrangian approach \cite{Cheng, CrSon, Filbet, Shoucri} which will be also the basis of the method here developed. Unfortunately, for each of the method cited, we recall that, the computational effort needed for solving the full six dimensional equation, prevents still nowadays realistic simulations even with parallel machines.

To overcome the problem of the excessive computational cost, we recently proposed in \cite{dimarco1} to use a splitting method to separate the transport from the collision step \cite{Des, Strang}. Then, we used a Lagrangian technique
which exactly solves the transport step on the entire domain and
we projected the solution on a grid to compute the contribution
of the collision operator. The resulting scheme (Fast Kinetic Scheme, FKS) shares many
analogies with semi-Lagrangian methods \cite{Cheng, CrSon, CrSon1, Filbet, Shoucri}
and with Monte Carlo schemes \cite{sliu}, but on
the contrary to them, the method is as fast as a particle method
while the numerical solution remains fully deterministic, which
means that there is no source of statistical error. When used to solve the limiting fluid model, the FKS method shares also some analogies with the so-called Lattice Boltzmann methods \cite{Abe}, but on the contrary to them its application is not limited to dense flows, all the regimes from rarefied to dense can be studied with such approach. Thanks to this new scheme, we were able to compute the solution of the full six
dimensional kinetic equation on a laptop for
acceptable mesh sizes and in a reasonable amount of time (about ten
hours for $100^{3}$ space $\times$ $12^{3}$ velocity space mesh points for $110$ time steps
up to $t=0.1$ on the spherical Sod-like problem).

However, the method developed, exhibited some limitations in term of spatial accuracy. In fact, it was only poorly accurate when used to compute solutions close to the fluid limit. We showed in \cite{dimarco1} that, in these cases, the computed solutions laid between a first order and a second order MUSCL scheme. In the present work, we developed a strategy which permits to preserve a desired high-order of accuracy in space for all the different regimes from the extremely rarefied to the high dense cases. In particular, it permits to recover the solution of the compressible Euler equations, when the number of collisions tends to infinity, with an high order shock capturing scheme. The modification introduced consists in coupling the fast kinetic scheme (FKS) to a solver for the compressible Euler equations, then to match the moments obtained from the solution of the macroscopic equations with those obtained from the solution of the equilibrium part of the kinetic equation. Finally, the solution, in term of the moments, is recovered as a convex combination of the two contributions: the macroscopic and the microscopic parts. We will show that the introduction of a macroscopic solver will not increase dramatically the computational cost, instead this modification will represent only a fraction of the time employed for computing the solution. In this work, the interaction term between the particles is the BGK collision operator \cite{Gross}. However, the high order version of the Fast Kinetic Scheme can in principle be extended to other collisional operators as the Boltzmann one.  \\

The article is organized as follows. In section \ref{sec_Boltzmann},
we introduce the Boltzmann-BGK equation and its properties. In
section \ref{sec_DVM}, we present the discrete velocity model (DVM).
Then, in section \ref{sec_num_approx}, we present the fast kinetic scheme (FKS) and the High Order Fast Kinetic Scheme (HOFKS).
Several test problems which demonstrate the accuracy and the strong efficiency of the new
method are presented and discussed in section \ref{sec_tests}. Some
final considerations and future developments are finally drawn in
the last section.

%%%%%%%%%%%%%%%%%%%%%%%%%%%%%%%%%%%%%%%%%%%%%%%%%%%%%%%%%%%%%%%%%%%%%%%%%%%%%%%%%
%%%%%%%%%%%%%%%%%%%%%%%%%%%%%%%%%%%%%%%%%%%%%%%%%%%%%%%%%%%%%%%%%%%%%%%%%%%%%%%%%
%%%%%%%%%%%%%%%%%%%%%%%%%%%%%%%%%%%%%%%%%%%%%%%%%%%%%%%%%%%%%%%%%%%%%%%%%%%%%%%%%
%%%%%%%%%%%%%%%%%%%%%%%%%%%%%%%%%%%%%%%%%%%%%%%%%%%%%%%%%%%%%%%%%%%%%%%%%%%%%%%%%
%
% S E C T I O N
%
\section{Boltzmann-BGK Equation}
\label{sec_Boltzmann}

The equation to be solved is the following: \be
\partial_t f + v\cdot\nabla_{x}f = \frac1{\tau} (M_{f}-f),
\label{eq:B}
\ee
with the initial condition $f(x,v,t=0)=f_{0}(x,v)$. The non negative function $f=f(x,v,t)$ describes the time
evolution of the distribution of particles which move with velocity
$v \in \R^d $ in the space $x \in \Omega \subset \R^{d}$ at time
$ t > 0$. For simplicity, in the description of the method we will do the hypothesis that the dimension of the physical space is the same of the dimension of the velocity space $d$. However, the method is not restricted to this particular choice and it is possible to consider different
dimensions between the space and the velocity in order to obtain different simplified models. The type of interactions term which characterizes the kinetic equation in (\ref{eq:B}) is the so-called BGK relaxation operator. With this choice the collisions are modeled by a relaxation towards the
local thermodynamical equilibrium defined by the Maxwellian
distribution function $M_{f}$ \be
 M_f=M_{f}[\rho,u,T](v)=\frac{\rho}{(2\pi \theta)^{d/2}}\exp\left(\frac{-|u-v|^{2}}{2\theta}\right) ,
\label{eq:M} \ee where $\rho \in \R, \ \rho>0$ and $u \in \R^d$ are respectively the
density and mean velocity while $\theta=RT$ with $T$ the temperature
of the gas and $R$ the gas constant. The macroscopic values
$\rho$,$u$ and $T$ are related to $f$ by: \be \rho=\int_{\R^d} fdv,
\qquad u=\frac{1}{\rho}\int_{\R^d} vfdv, \qquad  \theta=\frac{1}{\rho d}
\int_{\R^d}|v-u|^{2}fdv, \label{eq:Mo} \ee while the energy $E$ is defined
by \be E=\frac{1}{2} \int_{\R^d}|v|^{2}fdv = \frac{1}{2} \rho |u|^2
+ \frac{d}{2} \rho \theta, \label{eq:E} \ee %while the kinetic
%entropy of $f$ by \be H(f)=\int_{\R^{d}}f\log f dv.
% \label{eq:H}
%\ee
The parameter $\tau> 0$ in (\ref{eq:B}) is the relaxation time.
We refer to section \ref{sec_tests} for the numerical values chosen.

Formally, when the number of collision goes to infinity, which means
$\tau \rightarrow 0$, the function $f$ converges towards the
Maxwellian distribution. In this limit, if we consider the BGK equation (\ref{eq:B}) and we multiply it by $1$, $v$,
$\frac{1}{2}|v^{2}|$, and then we integrate with respect to $v$, we get the so-called Euler system of compressible gas dynamics
equations \be \ba \frac{\partial \rho}{\partial t} + \nabla_{x}
\cdot(\rho u) = 0, \ca \frac{\partial \rho u}{\partial t} +
\nabla_{x} \cdot (\rho u \otimes u+pI) = 0, \ca \frac{\partial
E}{\partial t} +\nabla_{x} \cdot((E+p)u) = 0, \ca p=\rho \theta,
\quad E=\frac{d}{2}\rho \theta +\frac{1}{2} \rho |u|^{2}. \ea
\label{eq:sys1} \ee
In the following, we will combine the BGK equation (\ref{eq:B}) with the compressible Euler equation (\ref{eq:sys1}) to get our High Order Fast Kinetic scheme.

%%%%%%%%%%%%%%%%%%%%%%%%%%%%%%%%%%%%%%%%%%%%%%%%%%%%%%%%%%%%%%%%%%%%%%%%%%%%%%%%%
%%%%%%%%%%%%%%%%%%%%%%%%%%%%%%%%%%%%%%%%%%%%%%%%%%%%%%%%%%%%%%%%%%%%%%%%%%%%%%%%%
%%%%%%%%%%%%%%%%%%%%%%%%%%%%%%%%%%%%%%%%%%%%%%%%%%%%%%%%%%%%%%%%%%%%%%%%%%%%%%%%%
%%%%%%%%%%%%%%%%%%%%%%%%%%%%%%%%%%%%%%%%%%%%%%%%%%%%%%%%%%%%%%%%%%%%%%%%%%%%%%%%%
%
% S E C T I O N
%
\section{The Discrete Velocity Models (DVM)} \label{sec_DVM}
The principle of Discrete Velocity Model (DVM) is to
set a grid in the velocity space and thus to transform the kinetic
equation (\ref{eq:B}) in a set of linear hyperbolic equations with source terms.
We refer to the works of Platkowski \cite{Plat} and of Mieussens \cite{Mieussens} for the
description of this approach and we remind to them for the details. In the following, we will briefly describe the idea and we will introduce the notations which will be used.

Let $\mathcal{K}$ be a set of $M$ multi-indices of $\mathbb{N}^{d}$,
defined by $\mathcal{K}=\left\{k=(k^{(i)})_{i=1}^{d},\ k^{(i)}\leq
K^{(i)}\right\}$, where $\{K^{(i)}\}$ are some given bounds. We
introduce a Cartesian grid $\mathcal{V}$ of $\mathbb{R}^{d}$ by \be
\mathcal{V}=\left\{ v_{k}=k\Delta v+a, \ k \in
\mathcal{K}\right\},\ee where $a$ is an arbitrary vector of
$\mathbb{R}^{d}$ and $\Delta v$ is a scalar which represents the
grid step in the velocity space. We denote the discrete collision
invariants on $\mathcal{V}$ by
$m_{k}=(1,v_{k},\frac{1}{2}|v_{k}|^{2})$.

Now, in this setting, the continuous distribution function $f$ is
replaced by a $N-$vector $f_{\mathcal{K}}(x,t)$, where each
component is assumed to be an approximation of the distribution
function $f$ at location $v_{k}$: \be
f_{\mathcal{K}}(x,t)=(f_{k}(x,t))_{k},\qquad f_{k}(x,t) \approx
f(x,v_{k},t). \ee The fluid quantities are then obtained from
$f_{k}$ thanks to discrete summations on $\mathcal{V}$: \be
U(x,t)=\sum_{k}m_{k}f_k(x,t)\, \Delta v=\langle m_k f_k(x,t)\rangle_{\mathcal{K}}. \label{eq:DM} \ee The
discrete velocity BGK model consists of a set of $N$ evolution
equations for $f_k$ of the form \be
\partial_t f_{k} + v_{k} \cdot\nabla_{x}f_{k} = \frac1{\tau} (\E_{k}[U]-f_{k}), \ k=1,..,N
\label{eq:DM1} \ee where $\E_{k}[U]$ is a suitable approximation of
$M_{f}$ defined next. The DVM approach deserves some remarks.
\begin{remark}{~}
\begin{itemize}
  \item When dealing with discrete velocity methods, one needs to truncate the velocity space and to fix
some bounds. This gives the number $N$ of evolution equations
(\ref{eq:DM1}). Of course, the number $N$ is chosen as a compromise
between the desired precision in the discretization of the velocity
space and the computational cost, while the bounds are chosen to
give a correct representation of the flow.
%It is well known \cite{Mieussens}, in fact, that the
%macroscopic velocity and temperature are bounded above by velocity
%bounds. Moreover, as a consequence of the velocity
%discretization, we have that the temperature is bounded from below.
This implies that the discrete velocity set must be large
enough to take into account large variations of the macroscopic
quantities which may appear during the evolution of the problem. On the other hand, the number of mesh points should be sufficiently large to guarantee that the small variations of the macroscopic quantities are well described.
%We summarize the above remarks by the following statement. Let $f$
%be a non negative distribution function, then the macroscopic
%velocity and temperature associated to $f$ in $\mathcal{V}$ by \be
%u=\frac{1}{\rho}\langle vf\rangle_{\mathcal{K}}, \qquad
%T=\frac{1}{dR\rho}\langle|v-u|^{2}f\rangle_{\mathcal{K}}, \ee where
%$\langle .\rangle_{\mathcal{K}}$ denotes the summation over the set
%of multi-indices $\mathcal{K}$, satisfy the bounds \cite{Mieussens}
%\bea \min_{\mathcal{K}} v_{k}^{(i)}  \leq
%& u^{(i)} & \leq \max_{\mathcal{K}}v_{k}^{(i)}, \; \; \forall i=1,\ldots,d \qquad  \\
%\frac{1}{dR}\min_{\mathcal{K}}|v-u|^{2}\leq
%&  T     & \leq \frac{1}{dR}\max_{\mathcal{K}}|v-u|^{2} .
%\eea
\item The exact conservation of macroscopic quantities is impossible, because
in general the support of the distribution function is non compact. This is the case for instance of the Maxwellian equilibrium distribution.
Thus, in order to conserve macroscopic variables, different
strategies can be adopted, two possibilities are described in
\cite{gamba, Mieussens}. Moreover, the approximation of the
equilibrium distribution $M_f$ by $\E_{k}[U]$ must be carefully
chosen in order to satisfy conservations of physical quantities. In the following
section we will discuss our choices in details.
%\item Starting from the discrete velocity model derived above, the natural choices which permits to solve the entire problem is then
%to discretize the $N$ evolution equations with the preferred finite volume or finite difference
%method \cite{Mieussens, Pal, Pal1, Pieraccini}. Alternatively, one
%can reconstruct the distribution function in space and then follows
%the characteristics backward in time to obtain the solution of the
%linear transport equation \cite{CrSon, CrSon1, Filbet2, Filbet}. However, both choices are extremely time consuming when dealing with multi dimensional problems. On the contrary, the fast kinetic scheme
%enables to drastically decrease the computational cost, because it is based on the exact solution of
%the linear transport equation but it avoids the reconstruction of the
%distribution function.
\end{itemize}
\end{remark}

%%%%%%%%%%%%%%%%%%%%%%%%%%%%%%%%%%%%%%%%%%%%%%%%%%%%%%%%%%%%%%%%%%%%%%%%%%%%%%%%%
%%%%%%%%%%%%%%%%%%%%%%%%%%%%%%%%%%%%%%%%%%%%%%%%%%%%%%%%%%%%%%%%%%%%%%%%%%%%%%%%%
%%%%%%%%%%%%%%%%%%%%%%%%%%%%%%%%%%%%%%%%%%%%%%%%%%%%%%%%%%%%%%%%%%%%%%%%%%%%%%%%%
%%%%%%%%%%%%%%%%%%%%%%%%%%%%%%%%%%%%%%%%%%%%%%%%%%%%%%%%%%%%%%%%%%%%%%%%%%%%%%%%%
%
%
% S E C T I O N
%

\section{Fast Kinetic Schemes (FKS) and High Order Fast Kinetic Schemes (HOFKS)}
\label{sec_num_approx} In this section we recall the FKS method and then we will introduce a new class of schemes which enables to get high order spatial accuracy (HOFKS). Before, we will discuss and propose a solution for the problem of lack of conservation of the macroscopic quantities which characterizes the class of discrete velocity models we are dealing with. We first of all introduce a Cartesian uniform grid in the physical space, FKS schemes are in fact based on uniform meshes. The extension of this class of methods to general meshes is not trivial but nevertheless under study.
The mesh is defined by the set $\mathcal{J}$ of $N$
multi-indices of $\mathbb{N}^{d}$, which is
$\mathcal{J}=\{j=(j^{(i)})_{i=1}^{d},\ j^{(i)}\leq J^{(i)}\}$, where
$\{J^{(i)}\}$ are some given bounds which represent the boundary
points in the physical space. The grid $\mathcal{X}$ of
$\mathbb{R}^{d}$ is then given by
\be
\mathcal{X}=\{ x_{j}=j\Delta x+b, \ j \in \mathcal{J}\},
\ee
where $d$ represents at the same time the
dimension of the physical space and the dimension of the velocity
space. The form of the physical space is determined by the vector $b$ of $\mathbb{R}^{d}$ and $\Delta x$ is a scalar which represents the grid step in the physical space. We consider a third discretization which is the time
discretization $t^{n}=n\Delta t$ with $n\in\mathbb{N}$.
We will at the end of the section discuss the time step limitations and the CFL condition.

\subsection{Conservative Discrete Velocity Models (DVM)}\label{sec:conserv}
Suppose a continuous in phase space distribution function is given, \textit{i.e.} $f(x_j,v,t^{n})$, with moments $U(x_j,t^{n})$ for every $j\in\mathcal{J}$ and $n\geq 0$. We proceed into two steps. First we define \be \widetilde{f}^{n}_{j,k}=f(x_{j},v_{k},t_{n}),\ee which is the pointwise distribution value in phase space. Observe that, due to the truncation of the velocity space and to the finite number of points with which $f$ is discretized, the moments of $\widetilde{f}^{n}_{j,k}$ differ from the original moments $U(x_j,t^{n})$. In fact the discrete moments of this distribution are
\be
\widetilde{U}_{j}^{n}=
\langle m_{k}\widetilde{f}_{j,k}^{n}\, \Delta v \rangle_{\mathcal{K}} \neq U_{j}^{n},
\qquad
\widetilde{\rho}^{n}_{j}\leq \rho^{n}_{j}, \qquad \widetilde{\theta}^{n}_{j}\leq \theta^{n}_{j}.
\ee
Different strategies can be adopted to restore the correct moments. Our choice, which is the second step of the conservative DVM model, consists in defining the approximated distribution function $f^{n}_{j,k}$ as the distribution closer in the discrete $L_2$ norm to $\widetilde{f}^{n}_{j,k}=f(x_{j},v_{k},t_{n})$ for which the moments are exactly the macroscopic quantities we want to preserve , \textit{i.e}
\be U_{j}^{n}=\langle m_{k}f_{j,k}^{n}\, \Delta v\rangle_{\mathcal{K}}. \ee

In order to find this distribution we make use of a simple
constrained Lagrange multiplier method \cite{gamba}, where the
constraints are mass, momentum and energy of the solution.
Let us recall the technique from \cite{gamba}.
For each spatial cell, let define the
pointwise distribution vector \be \widetilde{f}^{n}_{j} = \left(\widetilde{f}^{n}_{j,1},
\widetilde{f}^{n}_{j,2},\ldots,\widetilde{f}^{n}_{j,N} \right)^{T}, \ee let also define the vector containing the corrected distribution which fulfills the conservation of moments we are searching for
\be
f^{n}_{j} = \left(f^{n}_{j,1}, f^{n}_{j,2},\ldots, f^{n}_{j,N}\right)^{T},
\ee
and the matrix which contains the discretization parameters $ C \in \mathbb{R}^{(d+2)\times N} $.
%\be
%C_{(d+2)\times N}=\left(
%\begin{array}{ll}
%& (\Delta v)^{d}\\
%& v_k(\Delta v)^{d} \\
%& |v_{k}|^{2}(\Delta v)^{d}\\
%\end{array}
%\right).
%\ee
At this point, conservation can be imposed in each cell and at any time index $n$ solving the following
constrained optimization problem:
\begin{eqnarray}
% \nonumber to remove numbering (before each equation)
  &\nonumber  \mbox{ Given } \widetilde{f}^{n}_{j}\in \mathbb{R}^{N}, \ C \in \mathbb{R}^{(d+2)\times
  N}, \mbox{ and } U^{n}_{j} \in\mathbb{R}^{(d+2)\times 1},\label{eq:minim}\\
  & \mbox{ find } f^{n}_{j} \in \mathbb{R}^{N} \mbox{ such that } \\
  &\nonumber \|\widetilde{f}^{n}_{j} -f^{n}_{j}\|^{2}_{2} \mbox{ is minimized subject to the constraint } Cf^{n}_{j} = U^{n}_{j}.
\end{eqnarray}

Thus, let $\lambda\in
\mathbb{R}^{d+2}$ be the Lagrange multiplier vector, the
objective function to be minimized, in each cell, is given by
\be
L(f^{n}_{j}, \lambda) = \sum_{k=1}^{N} |\widetilde{f}^{n}_{j,k}-f^{n}_{j,k} |^{2} + \lambda^{T} (Cf^{n}_{j}- U^{n}_{j}).
\ee
The above equation can be solved explicitly. %In fact, taking the derivative of $L(f^{n}_{j}, \lambda)$ with
%respect to $f^{n}_{k,j}$, for all $k = 1, ...,N$ and $\lambda_i$, for all
%$i = 1, ..., d + 2$, we obtain that the value of $\lambda$ is uniquely determined by
%\be
%\lambda= 2(CC^{T})^{-1}(U^{n}_{j}-C \widetilde{f}^{n}_{j}),
%\ee
The searched distribution function is then $f^{n}_{j}$ is
\be
f^{n}_{j} = \widetilde{f}^{n}_{j} + C^{T} (CC^{T})^{-1}(U^{n}_{j}-C\widetilde{f}^{n}_{j}),
\quad \quad \forall  j\in\mathcal{J}.
\label{eq:minim1}
\ee

We end this part defining the approximated equilibrium distribution $\E_k[U_{j}^{n}]$, or
equivalently $\E^{n}_{j,k}[U]$. The discretization of the Maxwellian distribution $M_f(x,v,t)$,
should satisfy the same properties of conservation of the
distribution $f$, \textit{i.e.}
\be
U_{j}^{n}=\langle m_{k}f_{j,k}^{n}\, \Delta v\rangle_{\mathcal{K}}=\langle m_{k}\E_{k}[U_{j}^{n}]\, \Delta v\rangle_{\mathcal{K}}.
\ee
To this
aim, observe that the natural approximation $
\E_k[U^{n}_{j}]=M_{f}(x_{j},v_{k},t_{n})=M_f[U^n_j]$ cannot satisfy these requirements. Thus, the
calculation carried out above for the definition of the approximated
distribution $f$, should also be performed for the equilibrium
distribution $M_f$. This should be done each time we invoke the
equilibrium distribution during the computation as explained in the next subsection. The function
$\E[U]$ is therefore given by the solution of the same minimization
problem defined in (\ref{eq:minim}), and its explicit value is given
mimicking (\ref{eq:minim1}) by
\be
\E[U^{n}_{j}] = M_f[U^{n}_{j}] + C^{T} (CC^{T})^{-1}(U^{n}_{j}-CM_f[U^{n}_{j}]), \quad \quad \forall  j\in\mathcal{J}.
\label{eq:minimMax}
\ee
Notice that the computation of
the new distributions $f$ and $\E$ only involves a matrix-vector
multiplication. In fact, matrix $C$ only depends on the parameter of
the discretization and thus it is constant in time. In other words
matrices $C$ and $C^{T}(CC^{T})^{-1}$ can be precomputed and stored
in memory while initializing the problem.

%Some remarks on the conservative DVM method follow
\begin{remark}{~}
\begin{itemize}
\item For FKS schemes, we need to solve the above
minimization problem for the initial data $f(x_{j},v_{k},t^0=0)$ and for the distribution $\E[U^{n}_{j}]$ at each time
index $n$. In fact,
once the conservation is guaranteed for $f$ for $t=0$, this is
also guaranteed for the entire computation because the exact
solution is used for solving the transport step.
\item The only possible source of loss of conservation for this type of schemes is due to the solution of collision term and to the way in which the equilibrium distribution is discretized. This aspects will be detailed in the next subsection.
\item The conservation technique described in (\ref{eq:minim}) does not assure the positivity of the distribution function $f_{j,k}^{n}$. It may happen in fact that during the constrained minimization procedure $f_{j}^{n}$ becomes negative for some values of $k$. In practice, we did not observe this phenomenon to create instability in the solution. However, in the cases in which positivity is strictly demanded, as for instance for the full Boltzmann operator discretized with spectral methods, alternative techniques should be designed.
\end{itemize}
\end{remark}

\subsection{FKS schemes}
The main features of the method developed in \cite{dimarco1} can be summarized as follows:
\begin{itemize}
\item The BGK equation is discretized in velocity space by using the DVM model. The distribution $f$ as well as the Maxwellian $M_f$ are initialized by the conservative DVM method detailed in Section \ref{sec:conserv}.
\item A time splitting procedure is employed between the transport
and the relaxation operators for each of the resulting $N$ evolution
equations (\ref{eq:DM1}). First order time
splitting is considered. In principle, others more sophisticated splitting can be employed \cite{Strang}.
\item The transport part is solved exactly and continuously in space, this means that
no spatial mesh is involved. The initial data of this step is given by the
solution of the relaxation operator.
\item The relaxation part is solved on the spacial grid. The initial data
for this step is given by the value of the distribution function in
the center of the cells after the previous transport step. Each time the equilibrium distribution is invoked, conservation is retrieved through equation (\ref{eq:minimMax}). We need to impose conservation of the macroscopic quantities for the equilibrium distribution only.
\end{itemize}

Let us give the details of the method. We recall that the FKS methods are constructed on uniform grids.
%the extension of this approach to general grid is actually under study.
Let $f^{0}_{j,k}$ be the pointwise initial data, solution of
equation (\ref{eq:minim1}) with
$\widetilde{f}^{0}_{j,k}=f(x_{j},v_{k},t=0)$. Let also $\E^{0}_{j,k}[U]$ be the initial
equilibrium distribution solution of equation (\ref{eq:minimMax})
with $M^{0}_{j,k}=M_f(x_{j},v_{k},t=0)$ defined at points $x_j$ at $t=0$ as the distribution $f$. We start describing the first step of the method $[t^0;t^1]$ starting at $t^0=0$. The scheme is then generalizable to the generic time step $[t^n;t^{n+1}]$ starting from $t^n$.

\paragraph{First time step $[t^0;t^1]$.} Let us describe the transport
and relaxation stages.
\begin{description}
\item \textit{Transport stage.} We solve $N$ linear transport
equations of the form: \be
\partial_t f_{k} + v_{k} \cdot\nabla_{x}f_{k}=0, \quad k=1,\ldots,N
\label{eq:DVM}. \ee The idea of FKS schemes is to solve the transport part continuously in space instead of solving it only on the mesh points. To this aim, we define for each of the $N$ equations a piecewise constant function
in space as
\be
\overline{f}_{k}(x,t^0=0)=f^{0}_{j,k} \quad \forall x \in [x_{j-1/2},x_{j+1/2}], \quad  k=1,\ldots,N.
\ee
%The choice of a piecewise constant function in space is not mandatory
%for the method, on the contrary other choices can in principle improve the spatial accuracy of the method. However, in this work, we stick ourselves to the piecewise constant case and we remind to a future work for an analysis of the different initial shapes in space for the distribution $f$.
Now, the exact solution of the $N$ equations at time $t^1=t^0+\Delta t =
\Delta t$ is given by
\be
\overline{f}^{*}_{k}(x)=\overline{f}(x-v_{k}\Delta t), \quad k=1,\ldots,N.
\ee
Observe that, with this choice, we do not need to reconstruct the distribution function $f$ as for instance in the semi-Lagrangian schemes \cite{Filbet2, Filbet}; the shape of the function in space is in fact known
 at the beginning of the computation and it remains so through the duration of the computation. To be more precise the distribution function is transported in time with constant velocity so the discontinuities remain at the same relative locations. It remains also a piecewise constant function.
The relaxation step, as finite difference methods, is solved only on the
grid points. This means that only the value of the distribution function $f$ and the macroscopic quantities in the
centers of the cells are needed for this step. From the exact solution of the
function $f_{k}$ we can immediately recover these values at the cost
of one simple vector multiplication.
\item \textit{Relaxation stage.} This step is local to
the grid, this means that we solve the following
ordinary differential equation:
\be
\partial_t f_{j,k} =\frac{1}{\tau}(\E_{j,k}[U]-f_{j,k}), \ \ \ k=1,\ldots,N, \ \ \ j=1,\ldots,M
\label{eq:relax},
\ee
where the initial datum is the result of the
transport step at points $x_j$ at time $t^{1}=t^{0}+\Delta t$
\be \label{eq:f_star}
\overline{f}^{*}_{k}(x_{j}) = \bar{f}( x_j-v_k \Delta t), \ \ \ k=1,\ldots,N , \ \ j=1,\ldots,M.
\ee
To solve equation (\ref{eq:f_star}) we need the value of the equilibrium distribution $\E$ at
the center of the cell after the transport stage. In order to
compute the Maxwellian, the macroscopic quantities in the center of
the cells, \textit{i.e.} the density, the mean velocity and the
temperature, are given by summing the local value of the discrete
distribution $f$ over the velocity set: $\langle
m_{k}f^{*}_{j,k}\Delta v\rangle_{\mathcal{K}}=U^{*}_{j}$, for all
$j=1,\ldots,M$, where $f^{*}_{j,k}=\overline{f}^{*}_{k}(x_{j})$. The
discrete equilibrium distribution at time $t^1=t^0+\Delta t$, $\E^{*}_{j,k}=\E^{1}_{j,k}$, is the
solution of equation (\ref{eq:minimMax}) with moments $U^{*}_{j}=U^{1}_{j}$, for all
$j=1,\ldots,M$. Observe in fact that, the Maxwellian distribution does not change during the
relaxation step. In other words during this step the macroscopic
quantities remain constants. We can now compute the solution of the relaxation
stage as
\be
f^{1}_{j,k}= \exp(-\Delta t/\tau)f^{*}_{j,k}+(1-\exp(-\Delta t/\tau))\E^{1}_{j,k}[U].
\ee
The above equation
furnishes only the new value of the distribution $f$ at time
$t^1=t^0+\Delta t=\Delta t$ in the center of each spatial cell for
each velocity $v_k$. However in order to continue the
computation, we need the value of the distribution $f$ in all points of
the space. Let us assume that the equilibrium
distribution $M_f$ has the same shape than the distribution $f$ in
space. Thus, starting from the pointwise value of $\E$ we define a piecewise constant function in
space $\overline{\E}_k$ for each velocity $v_k$. The values of this piecewise constant
function are the values computed in the center of the spatial cells. In other words,
one defines
\be \overline{\E}^{*}_{k}(x)=\overline{\E}_{k}(x,t^1)=\E^{1}_{j,k},
\; \; \forall x \; \; \mbox{ such that } \; \;
\overline{f}^{*}_{k}(x)=\overline{f}^{*}_{k}(x_j), \ j=1,\ldots,M.
\ee
We can now rewrite the relaxation term directly in term of
spatial continuous function as
\be
\overline{f}^{1}_{k}(x)=
\overline{f}_{k}(x,\Delta t) =
  \exp(-\Delta t/\tau)\overline{f}^{*}_{k}(x)
  +(1-\exp(-\Delta t/\tau))\overline{\E}^{*}_{k}(x)[U].
\ee
For each velocity $v_k$ the original shape in space for
the distribution $f_k$ is preserved throughout the computation,
and, as a consequence it drastically reduces the computational cost because no reconstruction is needed.
\end{description}
The time marching procedure can be now be described.
\paragraph{Generic time step $[t^n;t^{n+1}]$.}
Given the value of the distribution
function $\overline{f}^{n}_{k}(x)$, for all $k=1,\ldots,N$, and all
$x \in \mathbb{R}^{d}$ at time $t^n$, the value of the distribution
at time $t^{n+1}$, $\overline{f}_{k}^{n+1}(x)$, is given by \be
\overline{f}_{k}^{*}(x)=\overline{f}^{n}_{k}(x-v_k\Delta t), \quad k=1,\ldots,N
\ee \be \overline{f}^{n+1}_{k}(x)=\exp(-\Delta
t/\tau)\overline{f}^{*}_{k}(x)+(1-\exp(-\Delta
t/\tau))\overline{\E}^{n+1}_{k}(x)[U], \quad k=1,\ldots,N,
\label{relaxdvm} \ee where $\overline{\E}^{n+1}_{k}(x)[U]$ is a
piecewise constant function with the discontinuities located in the same positions as the distribution $f^{*}_k$. It is computed considering the solution of the minimization problem (\ref{eq:minimMax}) relative to the moments
value in the center of each spatial cell after the transport stage:
$U^{n+1}_{j}, \ j=1,\ldots,M$. These moments are given by computing
$\langle m_{k}f^{*}_{j,k}\Delta v\rangle_{\mathcal{K}}$ where
$f^{*}_{j,k}$ is the value that the distribution function takes
after the transport stage in the center of each spatial cell.
\begin{remark}{~}
\begin{itemize}
\item Due to the fact that the
relaxation stage preserves the macroscopic quantities, the scheme is
globally conservative by construction. In fact, at each time step, the change of
density, momentum and energy is only due to the transport step. This
latter being exact, does preserve the macroscopic quantities as
well as the distribution function.
\item For the same reason, except for the constrained optimization
procedure\footnote{Observe that the positivity of the constrained optimization step can be forced introducing an inequality constraint
 of the type $f_{j,k}^{n}\geq 0$ or $\E_{j,k}^{n}\geq 0$. However, the introduction of such a step will cause the minimization
step to be solved numerically instead of analytically. This will means that the computational cost of the method will increase.
In the present work, we did not attack this problem and we remind to the future for the development of strictly positivity preserving fast kinetic schemes.},
the scheme is
also unconditionally positive. More precisely, if
$f^{n}_{k}(x)\geq 0$, and $k=1,\ldots,M$ and the optimization procedure preserves positivity, then $f^{n+1}_{k}(x)\geq 0$ if the
initial datum is positive $f^{0}_{k}(x)\geq 0$ for all
$k=1,\ldots,M$. In fact, the transport maintains the shape of $f$
unchanged in space while the relaxation towards the Maxwellian
distribution is a convex combination of $M_f$ and $f(x-v_k\Delta t)$
both being positive.
\item The time step
$\Delta t$ is constrained by the CFL condition
\be
\Delta t \max_{k} \left( \frac{|v_{k}|}{\Delta x} \right) \leq 1 = \text{CFL}.
\label{eq:Time}
\ee
Observe that this choice is not mandatory, in
fact the scheme is stable for every choice of the time step. However,
being based on a time splitting technique the error is of the
order of $\Delta t$. This suggests to take the
usual CFL condition in order to maintain the time splitting error small enough.
\item Some experiments have been done on the influence of the CFL condition on the schemes. The results showed that, for the cases tested, up to $CFL=5$ in (\ref{eq:Time}), the FKS scheme provides a solution very close to the case $CFL=1$ for all the values of the Knudsen number.
Moreover, when the Knudsen number is large, \textit{i.e} the BGK equations are very close to a free transport equation,
using larger values of the CFL number does not cause any more degradation to the global accuracy because \textit{de facto} the FKS scheme solves the transport term exactly.
\end{itemize}
\end{remark}

\subsection{HOFKS schemes}
As observed in \cite{dimarco1} the FKS scheme performs very well in collisionless or
almost collisionless regimes. In these cases, in fact, the relaxation
stage is neglectable and only the exact transport does play a role.
However, when moving from rarefied to dense regimes the projection over the
equilibrium distribution becomes more important. Thus, the accuracy
of the scheme was expected to diminish in fluid regimes, because the
projection method is only first order accurate. These behaviors were, in fact, observed in the numerical simulations performed \cite{dimarco1}.
In the present paper, we developed a method which preserves the high spatial accuracy observed with the FKS schemes for rarefied regimes and which becomes a high order shock capturing scheme applied to the kinetic equation (\ref{eq:B}) in the fluid limit. This means that throughout all possible regimes, from fluid to extremely rarefied flows, the new scheme maintains high accuracy in space. Moreover, the new method does not cause the computational efficiency to drop down. As shown in the numerical test section the high order fast kinetic schemes (HOFKS) still works with computational costs close to the original FKS method and, for unsteady problems, in which time averaging are unusable, it is still much faster than DSMC methods. We recall that we only focus on the spatial accuracy in this paper, we postpone to the future the development of high order schemes both in time and space.

\subsubsection{The general methodology}
The idea, onto which the method is based, is to solve
the equilibrium part of the distribution function with a macroscopic scheme instead of a kinetic scheme.
In fact, observe that at each time step, the relaxation stage consists in computing the distribution function $f^n$ as a convex combination of the transported distribution $f^*$ and a Maxwellian distribution $\E^n$. The Maxwellian distribution $\E^n$ is computed through the moments of $f^*$. Then, in order to complete one step in time, the scheme solves the transport part which leads to the new intermediate distribution $f^*$ at time $n+1$. So now, in HOFKS scheme, we replace the moments at time index $n+1$, obtained from the solution of the transport stage at time $n$ with another set of moments. This new set of moments are computed through the same convex combination of the relaxation stage (\ref{relaxdvm}). However now, to the contrary of (\ref{relaxdvm}), the convex combination is performed between the moments which come from the solution of the transport part of the kinetic equation and the moments which are solutions of the compressible Euler equations. At this point, if the compressible Euler equations are solved with a high order shock capturing scheme in the limit $\tau\rightarrow 0$ the HOFKS corresponds to the same method for the macroscopic equations.
We detail this new scheme in the sequel.

In order to keep notations simple and compact we introduce three operators:
the projection operator $\P$, the relaxation
operator $\R_{\Delta t}$ and the transport operator $\T_{\Delta t}$ which act on a time step $\Delta t$.
\begin{description}
\item From the kinetic variable
$f$ (or $M_{f}$) the projection operator computes the macroscopic averages
$U(x_j,t^{n})=U^{n}_j$, thus
\be
\P^{n}_{j}(f)=\langle m_{k}f_{j,k}^{n}\, \Delta v\rangle_{\mathcal{K}}=\langle
m_{k}\E_{k}[U_{j}^{n}]\, \Delta v\rangle_{\mathcal{K}}=\P^{n}_{j}(\E)=U(x_{j},t^{n}),
\ee
since the local Maxwellian $M_f$ has the same moments of the distribution function $f$.
\item
The relaxation and transport operators
solve the relaxation and transport steps for piecewise constant functions $\overline{f}_{k}$ and $\overline{\E}_{k}$
or equivalently for pointwise functions $f_{j,k}$ and $\E_{j,k}[U]$ for all velocities $v_k, \ k=1,..,N$.
The relaxation operator has the form
\be
\R_{\Delta t}\left(\overline{f}\right)=\lambda \overline{f}+(1-\lambda) \overline{\E},
\ee
where $\lambda=\exp(-\Delta t/\tau)$, whereas the transport operator reads
\be
\T_{\Delta t}\left(\overline{f}\right)=\overline{f}(x-v\Delta t).
\ee
\end{description}
In order to describe the HOFKS scheme, let us start, to the contrary of FKS scheme, from the relaxation step.
Recall that starting either from the relaxation or from the transport step gives consistent splitting discretizations.
Suppose that the distribution function $\overline{f}^{n}_{k}(x)$ is known as well as the function $\overline{\E}^{n}_{k}(x)$ for all $k=1,\ldots,N$,
then, as first step of the splitting we have
\be
\overline{f}^{*}=\R_{\Delta t}\left(\overline{f}^{n} \right)=\lambda \overline{f}^{n}+ (1-\lambda) \overline{\E}^{n}.
\ee
The distribution function $\overline{f}^{*}$ is given by a convex combination of the transported distribution and the Maxwellian distribution.
Then, the transport step, applied to the solution of the relaxation step, produces the so called kinetic solution ($K$) at time index $n+1$
\be
\nonumber
\overline{f}_{K}^{n+1}=\T_{\Delta t}\left(\overline{f}^{*}\right)=\T_{\Delta t}\left(\lambda\overline{f}^{n}\right)+\T_{\Delta t}\left((1-\lambda) \overline{\E}^{n}\right).
\ee
On the other hand, the
kinetic solution in terms of the macroscopic moments furnishes the following values
\bea
\label{eq:UK}
U_{K}(x_j,t^{n+1})&=&\P^{n+1}_{j}\left(\T_{\Delta t}(\overline{f}^{*})\right)  \\
                &=& \P^{n+1}_{j}(\T_{\Delta t}(\lambda\overline{f}^{n}))+\P^{n+1}_{j}(\T_{\Delta t}((1-\lambda) \overline{\E}^{n}))\\
                &=& U_{*}(x_j,t^{n+1})+U_{M}(x_j,t^{n+1}).
\eea
In order to construct the HOFKS scheme we replace the moments $U_M(x_j,t^{n+1})$ by the moments obtained solving the compressible Euler equations
that we call $U_E(x_j,t^{n+1})$.
(The details of the numerical scheme used will be given next.)  \\
Thus the final moments used in the solution are given by
\be
U_H(x_j,t^{n+1})= U_*(x_j,t^{n+1})+U_E(x_j,t^{n+1})
\label{eq:mom}\ee
where $U_H$ stands for hybrid. Before describing the last step which ensures consistency
between the kinetic solver for the Maxwellian distribution and the macroscopic solver
for the compressible Euler equations, we state the following result
(see \cite{dimarco2, perthame2} for a proof):
\begin{theorem}
If we denote by $U_E(x,t+\Delta t)$ the solution of the Euler
equations (\ref{eq:sys1}) and with $U_M(x,t+\Delta t)$ the solution of the kinetic equation in which we consider initial thermodynamical equilibrium, \textit{i.e.} $f=\E[U]$. If in addition we consider as initial data
$U_E(x,t)=U_M(x,t)$ then
\be
U_E(x,t+\Delta t)=U_M(x,t+\Delta t)+O(\Delta t^2).  \label{eq:est}
\ee
\label{th:1}
\end{theorem}
By virtue of the above result, we can replace the moments after the transport $U_M(x_j,t^{n})$ with $U_{E}(x_j,t^{n})$ at each time step without
affecting the overall first order accuracy in time of the splitting method. %and get $U_{H}(x_j,t^{n})$.
However, to have consistency between the macroscopic solution and the kinetic discretized solution, it is necessary that the advected equilibrium satisfies
\be
\P^{n+1}_{j} \left(\T_{\Delta t}((1-\lambda) \overline{\E}^{n}) \right) = U_E(x_j,t^{n+1}),
\label{eq:cons}
\ee
namely the kinetic solution to the fluid equations in
one time step should match the direct solution to the limiting
fluid equations. This is not true in general. To solve this problem, we apply again the minimization method (\ref{eq:minimMax}) to find the new distribution $\T_{\Delta t}\left((1-\lambda)\E^{n}_{j}\right)$ which shares the same moments than the macroscopic solution $U_E(x_j,t^{n+1})$. Thus, we search for a distribution
$\T'_{\Delta t} \left((1-\lambda)\E^{n}_{j} \right)$ which satisfies the following minimization problem:
\begin{eqnarray}
% \nonumber to remove numbering (before each equation)
  &\nonumber  \mbox{ Given } \T_{\Delta t}((1-\lambda)\E^{n}_{j})\in \mathbb{R}^{N}, \ C \in \mathbb{R}^{(d+2)\times
  N}, \mbox{ and } U_E(x_j,t^{n+1}) \in\mathbb{R}^{(d+2)\times 1},\label{eq:minim_full}\\
  & \mbox{ find } \T'_{\Delta t}((1-\lambda)\E^{n}_{j}) \in \mathbb{R}^{N} \mbox{ such that } \\
  &\nonumber \|\T_{\Delta t}((1-\lambda)\E^{n}_{j}) -\T'_{\Delta t}((1-\lambda)\E^{n}_{j})\|^{2}_{2} \mbox{ is minimized subject to the constraint} \\ & \nonumber C\left(\T'_{\Delta t}((1-\lambda)\E^{n}_{j})\right) = U_E(x_j,t^{n+1}).
\end{eqnarray}
Then again, starting from the pointwise solution, we define as in the FKS schemes a new piecewise constant equilibrium function
$\T'_{\Delta t}((1-\lambda)\overline{\E})(x)$ sharing its shape with $\bar{f}$ as
\be
\T'_{\Delta t}\left((1-\lambda)\overline{\E}^{n}_{k} \right)(x)=\T'_{\Delta t} \left((1-\lambda)\E^{n}_{j,k} \right),
\; \; \forall x \; \mbox{s.t.} \; \;
\overline{f}^{n}_{k}(x)=\overline{f}^{n}_{k}(x_j), \ j=1,\ldots,M.
\ee
Finally, the new distribution $f$, at time index $n+1$, is defined as
\be
\overline{f}^{n+1}=\T_{\Delta t}(\lambda\overline{f}^{n})+\T'_{\Delta t}((1-\lambda)\overline{\E}^{n}_{k}),
\ee
while the new moments are given by (\ref{eq:mom}).
This somehow ends the methodology to design HOFKS schemes.
What remains to be detailed is how the
solution of the compressible Euler equations $U_E(x_j,t^{n+1})=U_{E,j}^{n+1}$ is computed.
Remark that at this point any solver for the compressible Euler equations
can be used. One example is proposed in the following.

\subsubsection{One example: MUSCL Finite Volume (FV) scheme}
\label{sssec:MUSCL}

As a first example we propose the MUSCL Finite Volume (FV) scheme. This is also the scheme used in the numerical test section. It reads, starting from $U_{E,j}^{n}$,
\be
\frac{U_{E,j}^{n+1}-U_{E,j}^{n}}{\Delta t} +\frac{\psi_{j+1/2}(U_E^n)-\psi_{j-1/2}(U_E^n)}{\Delta x} = 0,
\label{eq:discmom}
\ee
where discrete fluxes are defined as in \cite{leveque:numerical-methods} by
\be
\psi_{j+1/2}(U_E^n) =
 \frac{1}{2}(F(U^n_{E,j})+F(U^n_{E,j+1})) -
 \frac{1}{2}\alpha(U^n_{E,j+1}-U^n_{E,j}) +
 \frac{1}{4}(\sigma^{n,+}_j-\sigma^{n,-}_{j+1})
\ee
with $F(U)$ the flux of the compressible Euler equations and
\be
\sigma^{n,\pm}_j=\left(F(U^n_{E,j+1})\pm \alpha U^n_{E,j+1}-F(U^n_{E,j})\mp \alpha U^n_{E,j}\right)\varphi(\chi^{n,\pm}_j)
\ee
with $\varphi$ being the slope limiter, as instance we use the Van Leer slope limiter
\be
\varphi(\chi)=\frac{|\chi|+\chi}{1+\chi},
\ee
where the variable $\chi^{\pm}$ is defined as following
\be
\chi^{n,\pm}_j=\frac{F(U^n_{E,j})\pm \alpha
  U^n_{E,j}-F(U^n_{E,j-1})\mp \alpha U^n_{E,j-1}}{F(U^n_{E,j+1})\pm \alpha
  U^n_{E,j+1}-F(U^n_{E,j})\mp \alpha U^n_{E,j}}.
\ee
The above ratio of vectors is defined componentwise and
$\alpha$ represent the eigenvalues of the Euler system.

\subsubsection{Extensions and stability constraints}

We have proposed second order schemes to solve the compressible Euler equations.
However, any other solver can in principle be used, which could additionally increase
the spatial accuracy of the HOFKS method, for instance WENO methods \cite{Shu98}
or a genuine MOOD scheme \cite{MOOD1, MOOD2}. \\
Finally, the time step $\Delta t$ is chosen such that it satisfies the stability condition of the Euler solver,
in fact we recall that the FKS is stable for all choices of the time step.
This means that the time step is driven for the MUSCL scheme by
\be
\Delta t =\frac{1}{2}  \left(\frac{\Delta x}{\alpha_{\max}} \right).
\label{eq:Time2}
\ee
with $\alpha_{\max}$ the largest eigenvalue of the Euler system. If a more CFL restrictive Euler solver has to be used, as instance with a WENO scheme, then the time step will be accordingly reduced. \\
\section{Numerical tests}
\label{sec_tests}

%%%%%%%%%%%%%%%%%%%%%%%%%%%%%%%%%%%%%%%%%%
%\subsection{General setting}
%\label{subsec_gen_set}

In this section, we present several numerical tests to illustrate
the main features of the method and the improvements with respect to the FKS scheme.
The following methodology is adopted
\begin{description}
\item
  First, we test the HOFKS method on the one dimensional Sod shock tube.
  In this case, we compare two kinetic schemes (the FKS, a first order upwind scheme and the HOFKS
  a second order scheme) versus a finite volume (FV) upwind scheme, a second order finite volume MUSCL method (MUSCL)
  and a third order WENO finite difference scheme (WENO).
  For any fluid regime the goal is threefold:
  (i) the two kinetic schemes produce valide results,
  (ii) the HOFKS is more accurate than FKS and,
  (iii) the accuracy of FKS lays in between FV and MUSCL
  and the accuracy of HOFKS lays in between MUSCL and WENO.
\item
  In a second test case we use the exact smooth solution of the advected isentropic vortex of
  the 2D Euler equations to assess the effective numerical accuracy and rate of convergence of FKS, HOFKS and also
  the unlimited version of HOFKS. The solution being smooth an unlimited scheme can be used to measure the maximal accuracy
  that can be obtained with our choice of Euler solver.
\item
  Then, in a third series of tests we solve a two dimensional-two dimensional BGK equation and we compare
  our method with a Monte Carlo scheme (DSMC) for $\tau=10^{-3}$, and, in the fluid limit, $\tau=10^{-4}$, with FV and MUSCL schemes.
  The goal is to show that a genuine multi-dimensional solutions with shocks and interaction of waves can be accurately captured with HOFKS
  in any fluid regime.
  We also report the computational times for the two dimensional simulations for the HOFKS, the FKS and the DSMC methods.
  All simulations are performed on a mono-processor laptop machine.
\end{description}

%=== B E G I N   F I G ================================================
\begin{figure}[h!]
\begin{center}
\includegraphics[scale=0.43]{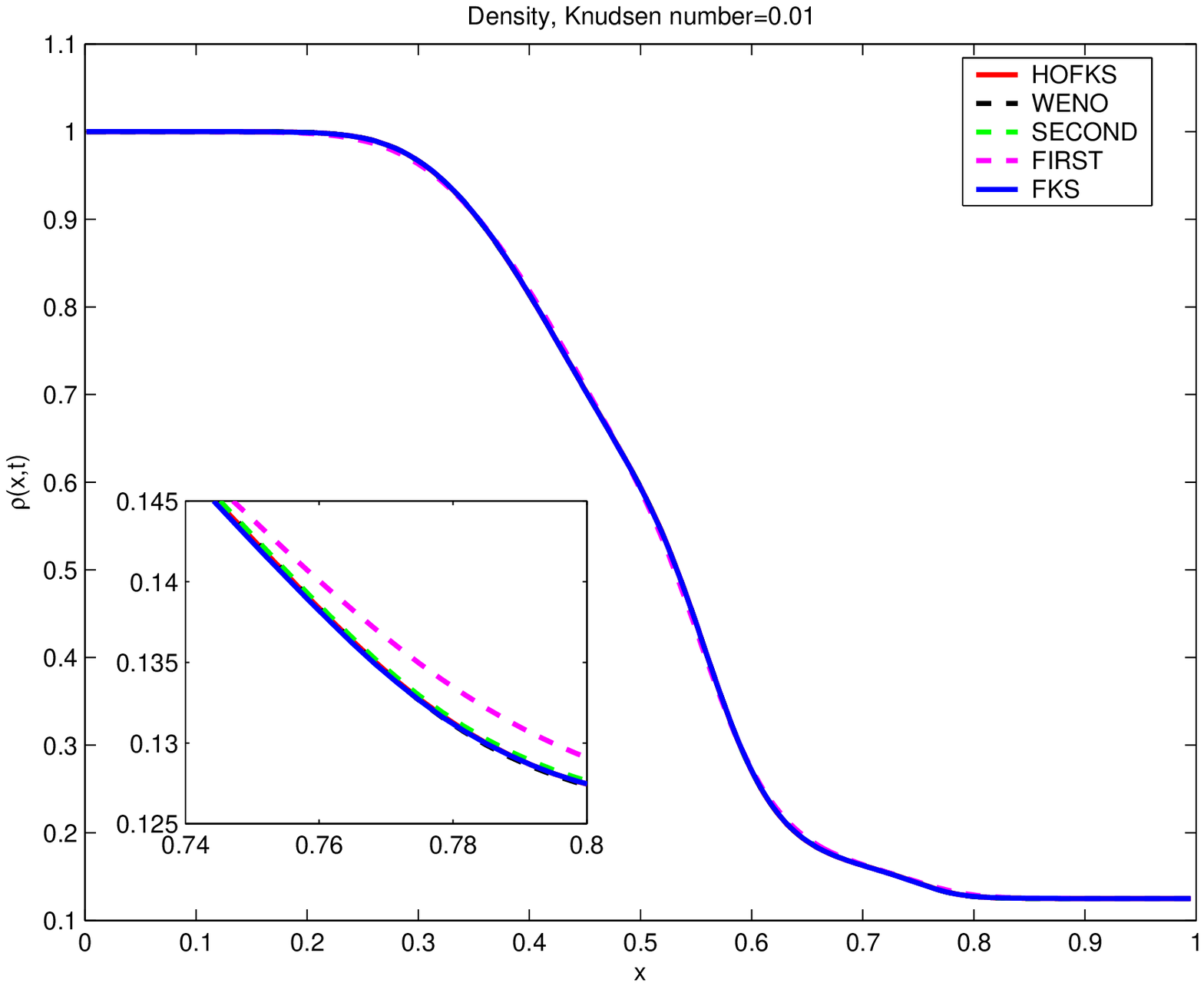}
\includegraphics[scale=0.43]{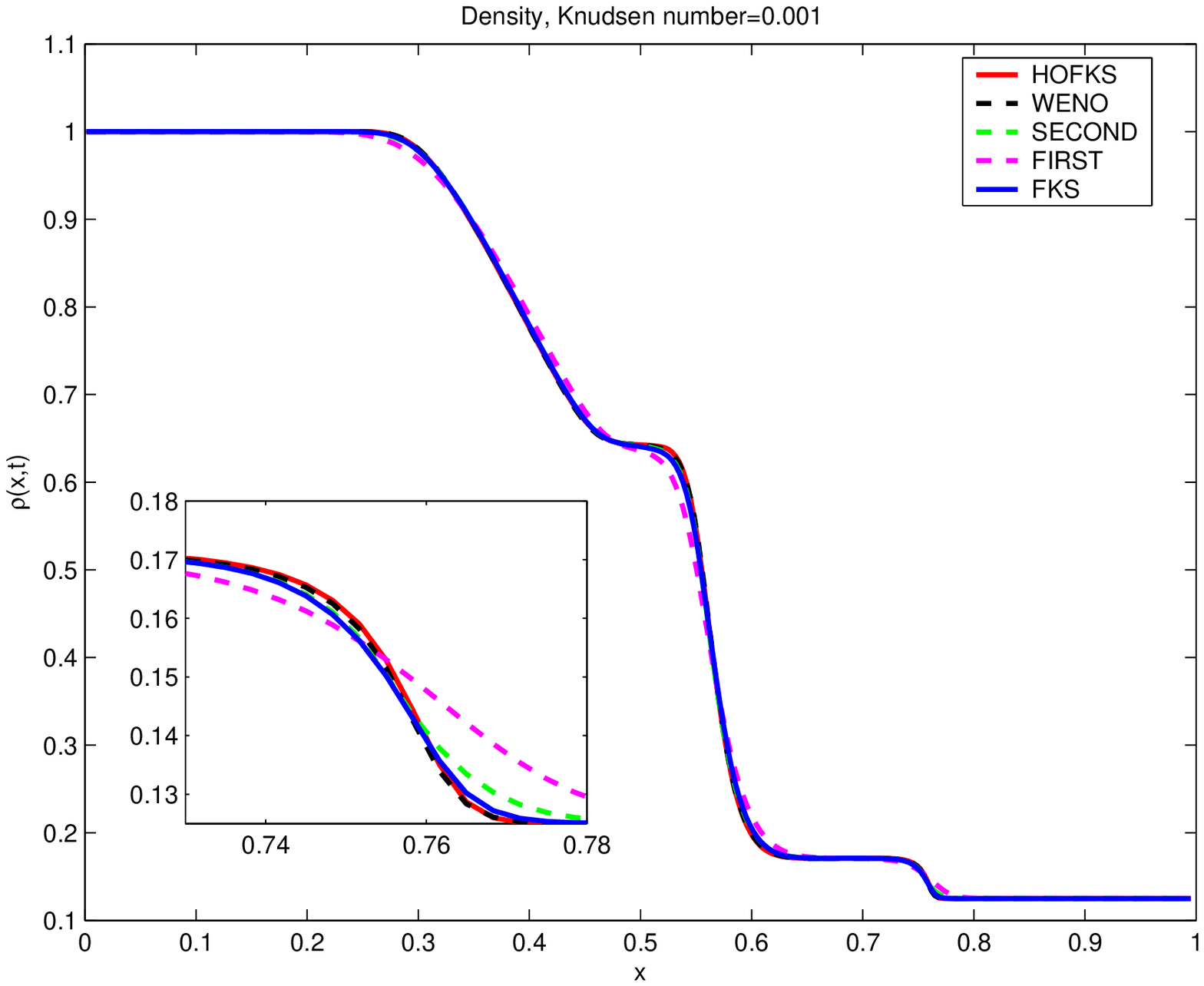}\\
\includegraphics[scale=0.43]{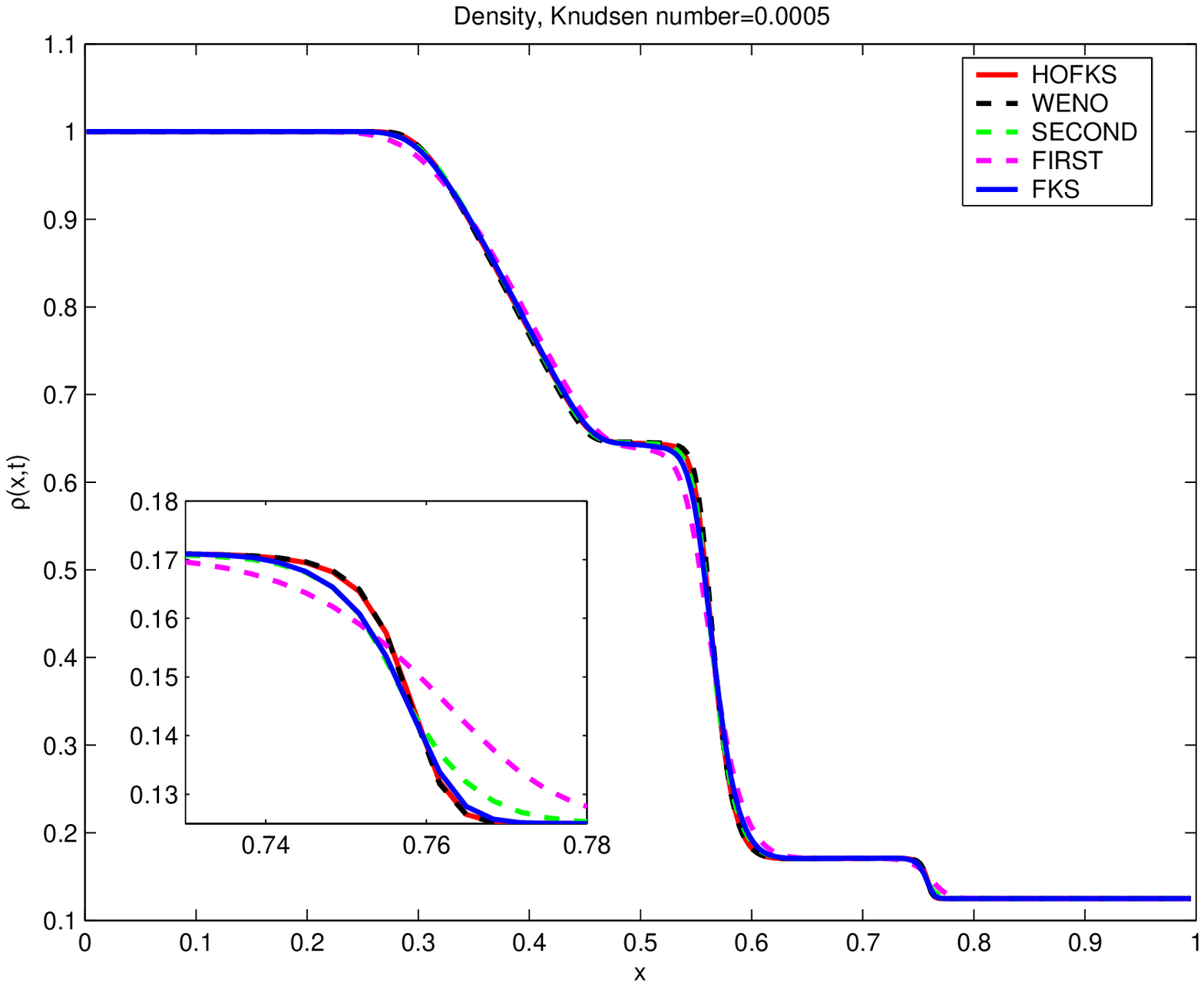}
\includegraphics[scale=0.43]{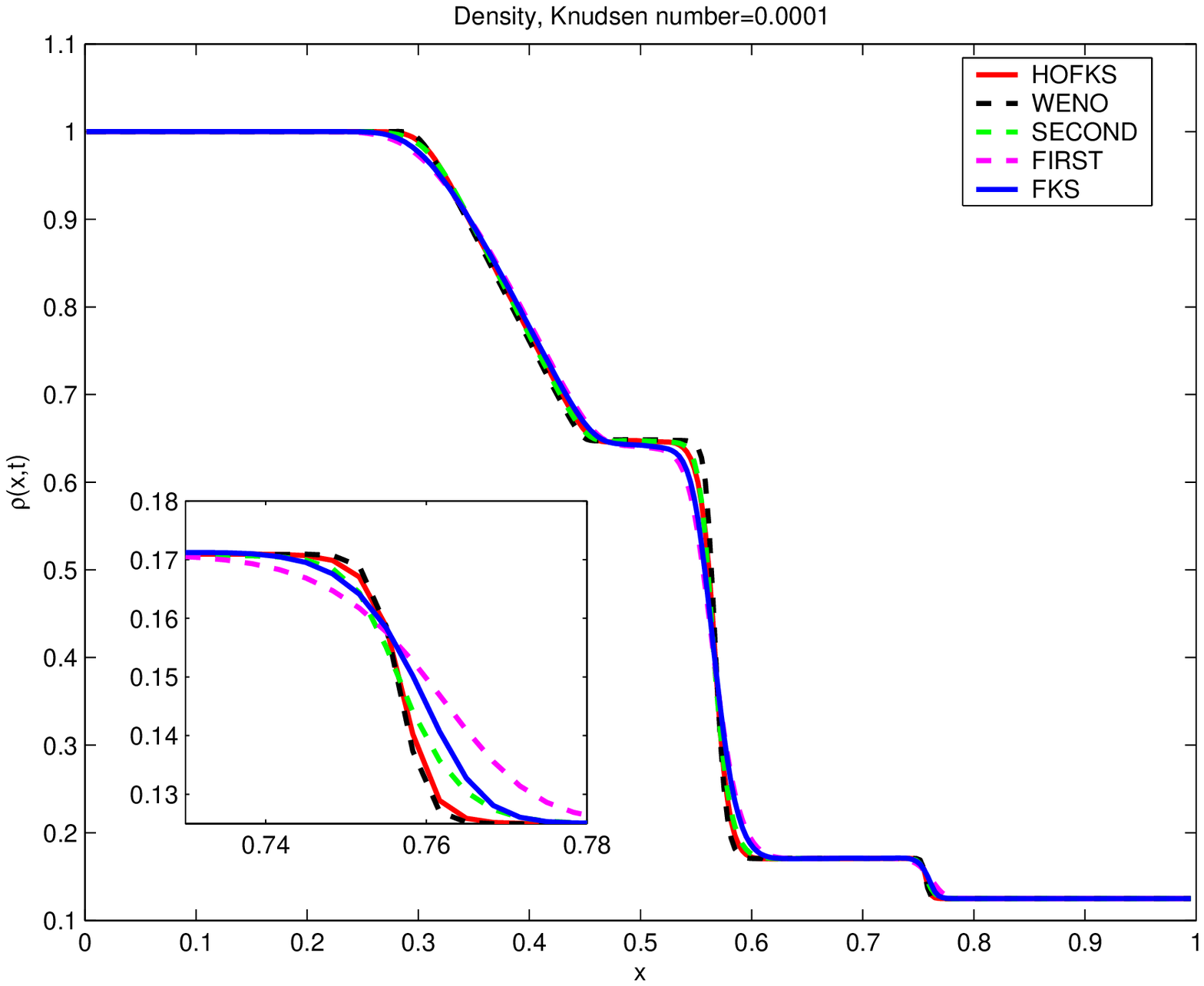}\\
\caption{Sod test: solution at $t_{\text{final}}=0.05$ for the
density, with $\tau=10^{-2}$ (top left), $\tau=10^{-3}$ (top right),
$\tau=5 \times 10^{-4}$ (bottom left) and $\tau=10^{-4}$ (bottom right).}
\label{sod1}
\end{center}
\end{figure}
%=== E N D   F I G ================================================

%=== B E G I N   F I G ================================================
\begin{figure}[h!]
\begin{center}
\includegraphics[scale=0.425]{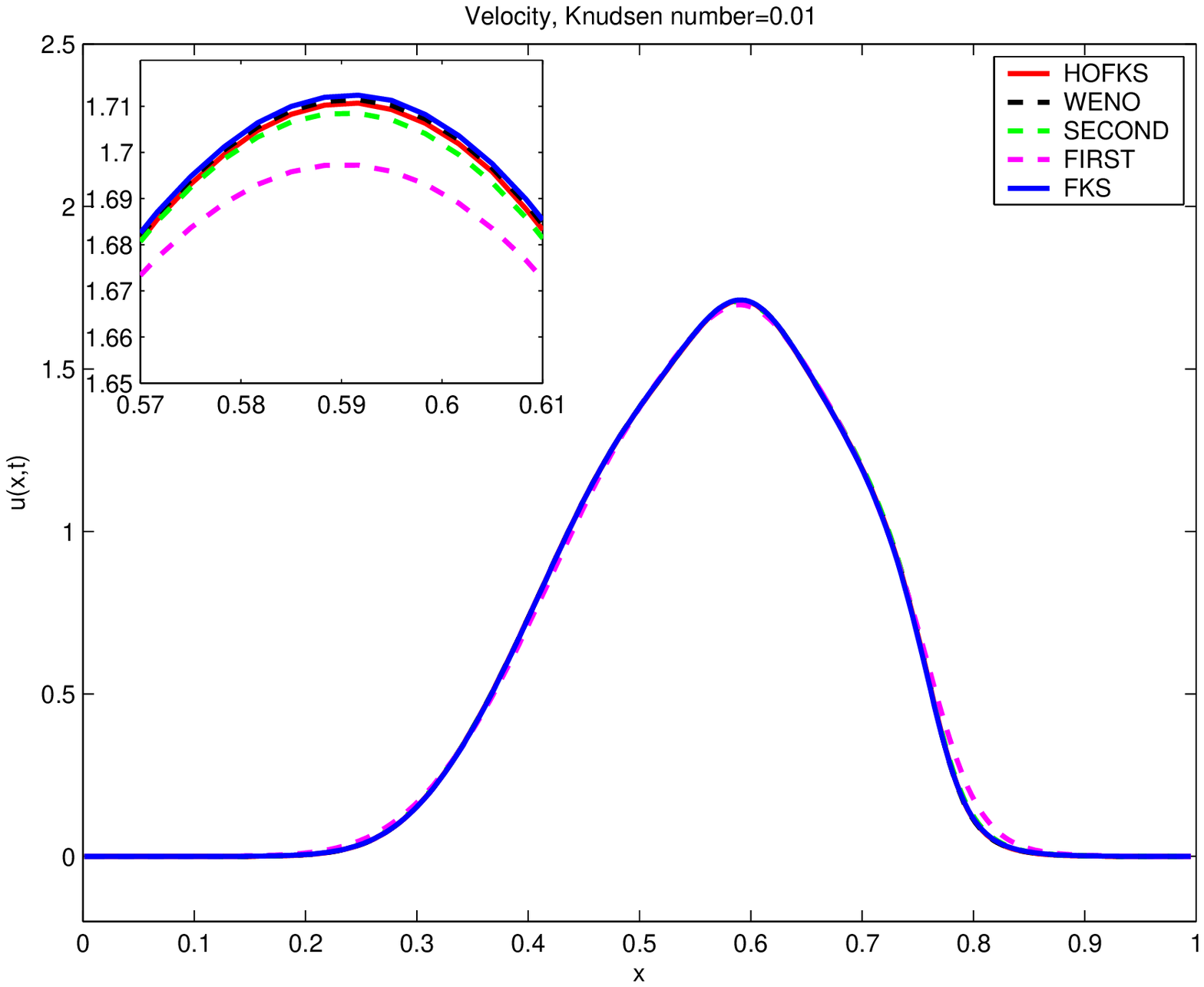}
\includegraphics[scale=0.425]{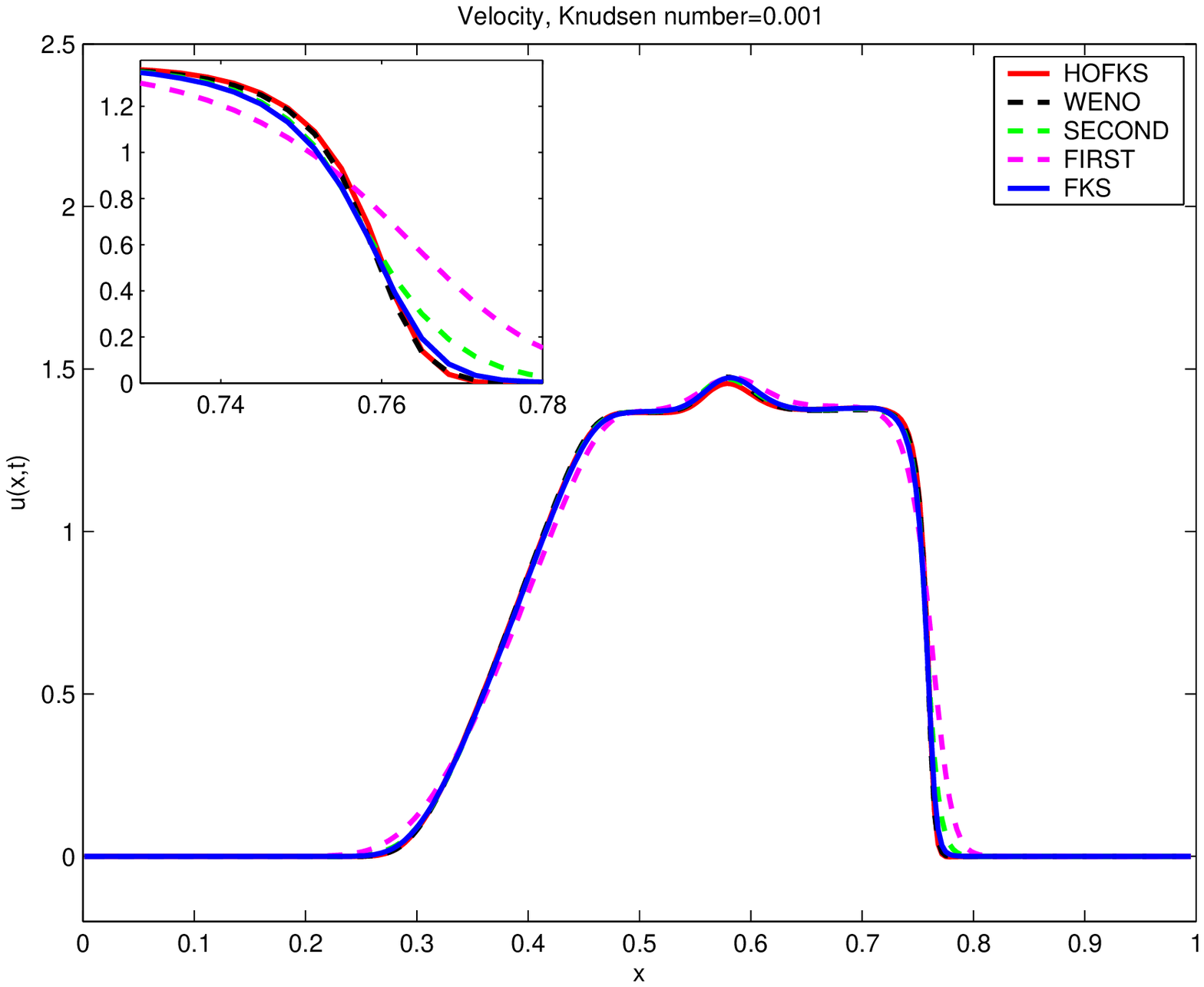}\\
\includegraphics[scale=0.425]{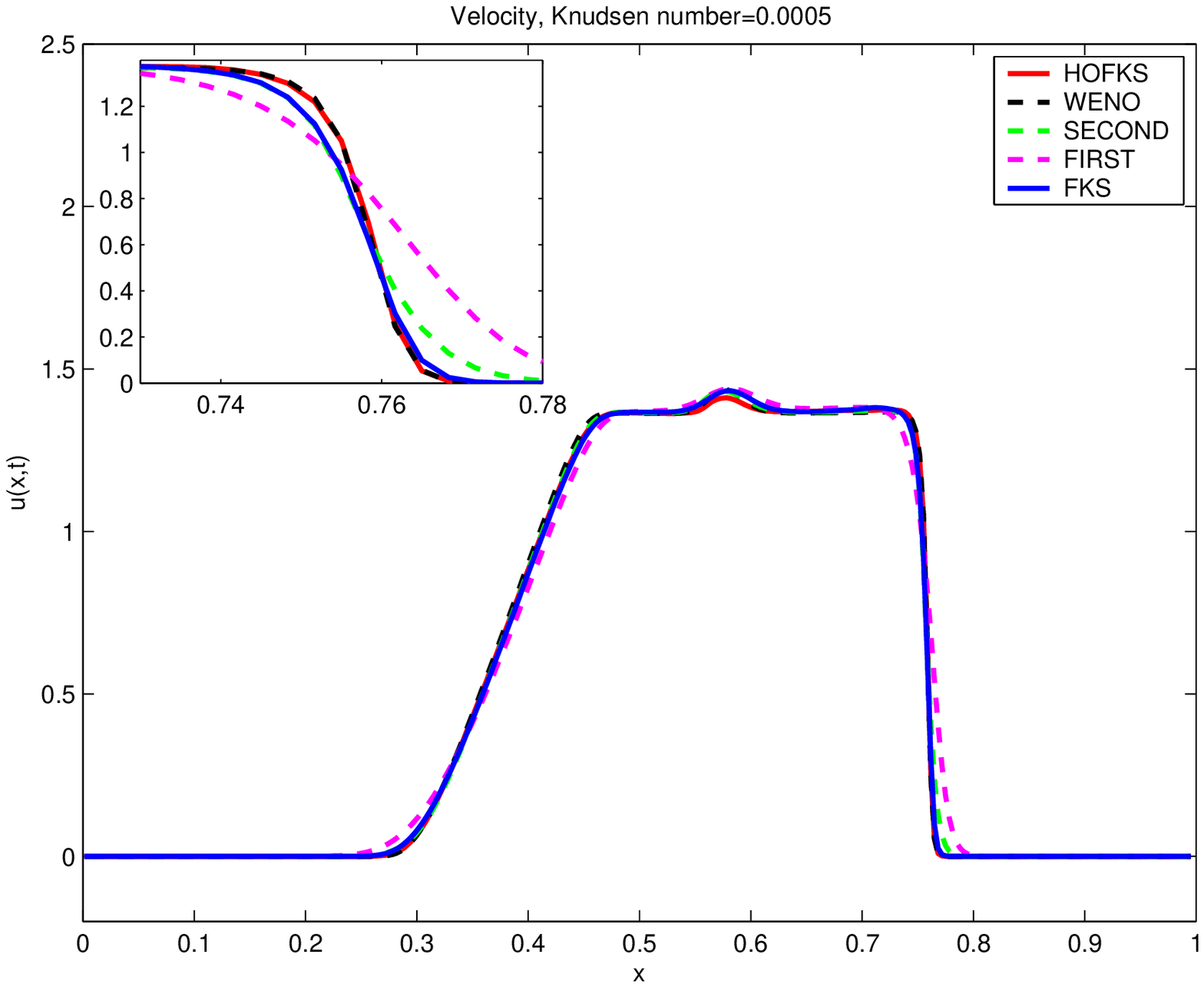}
\includegraphics[scale=0.425]{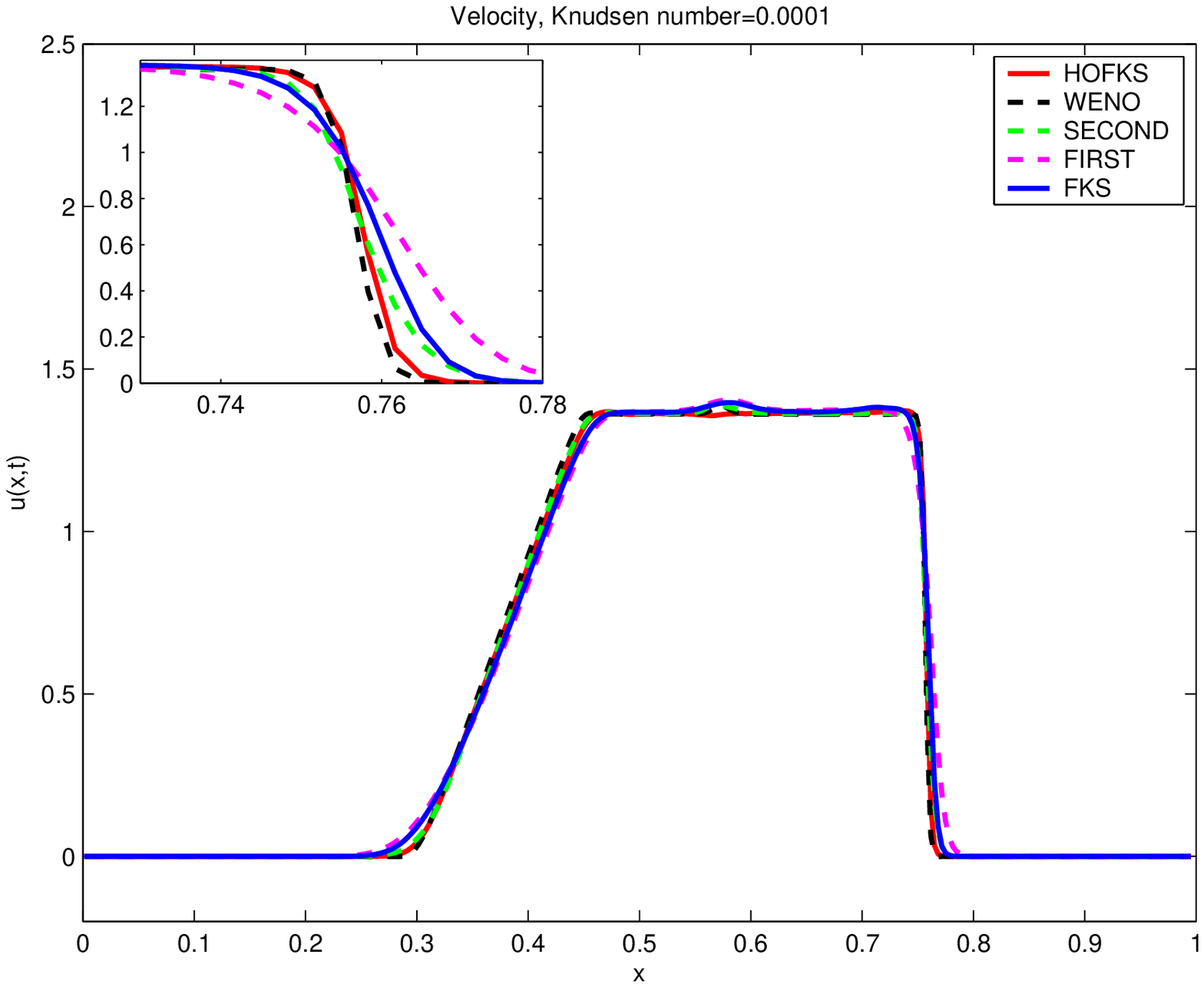}\\
\caption{1D Sod test: solution at $t_{\text{final}}=0.05$ for the
mean velocity, with $\tau=10^{-2}$ (top left), $\tau=10^{-3}$ (top
right), $\tau=5 \times 10^{-4}$ (bottom left) and $\tau=10^{-4}$ (bottom
right).} \label{sod2}
\end{center}
\end{figure}
%=== E N D   F I G ================================================
%=== B E G I N   F I G ================================================
\begin{figure}[h!]
\begin{center}
\includegraphics[scale=0.43]{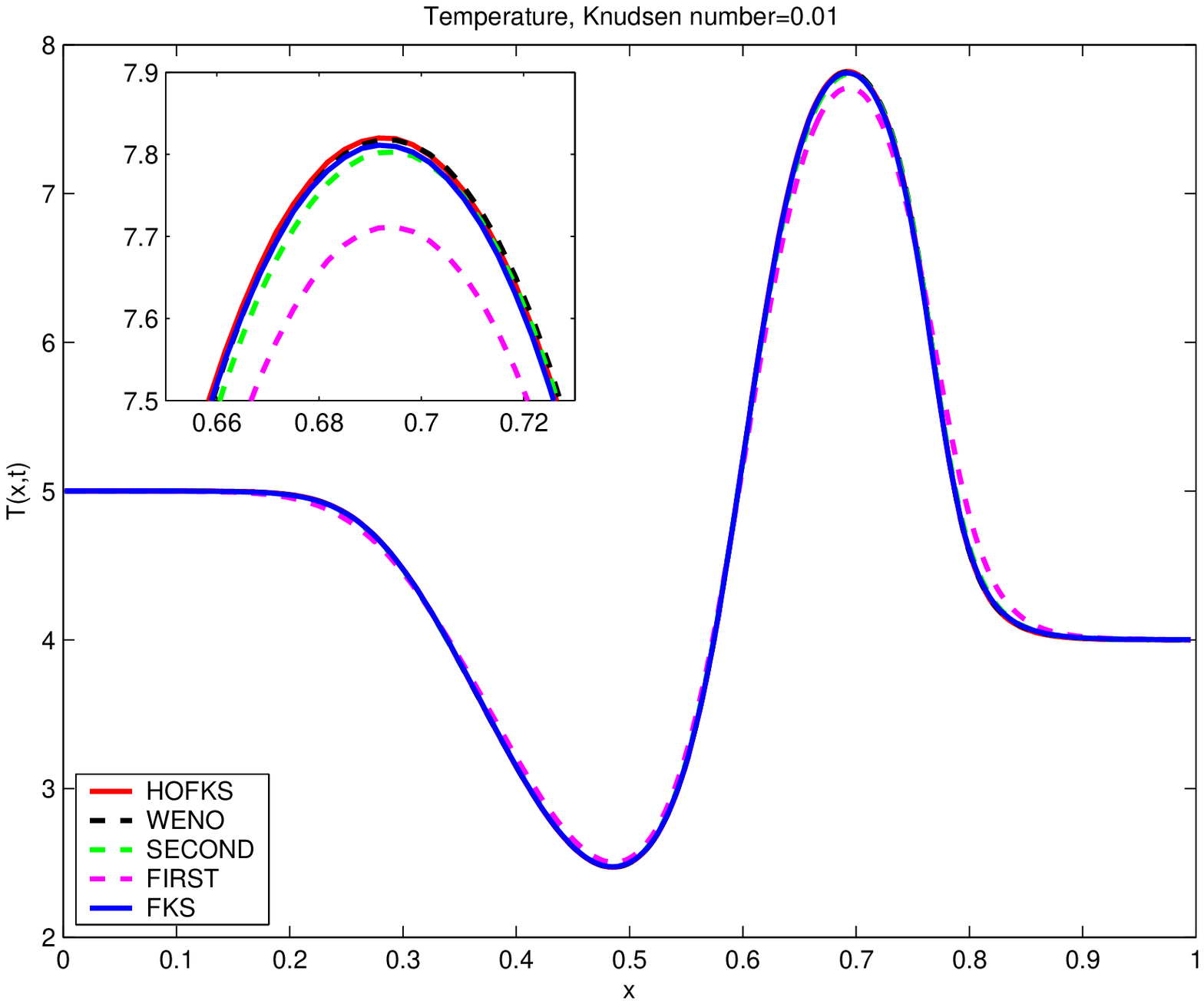}
\includegraphics[scale=0.43]{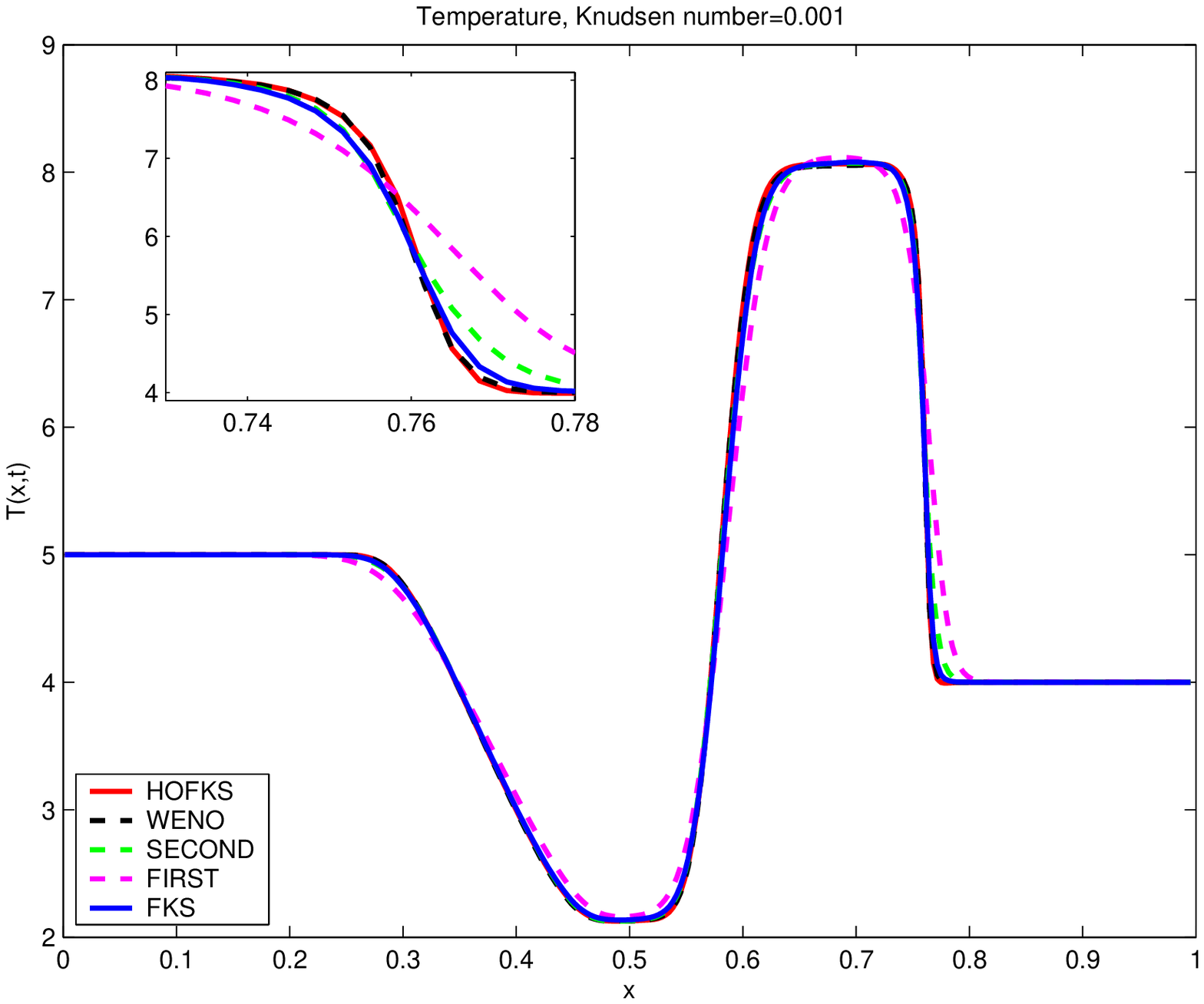}\\
\includegraphics[scale=0.43]{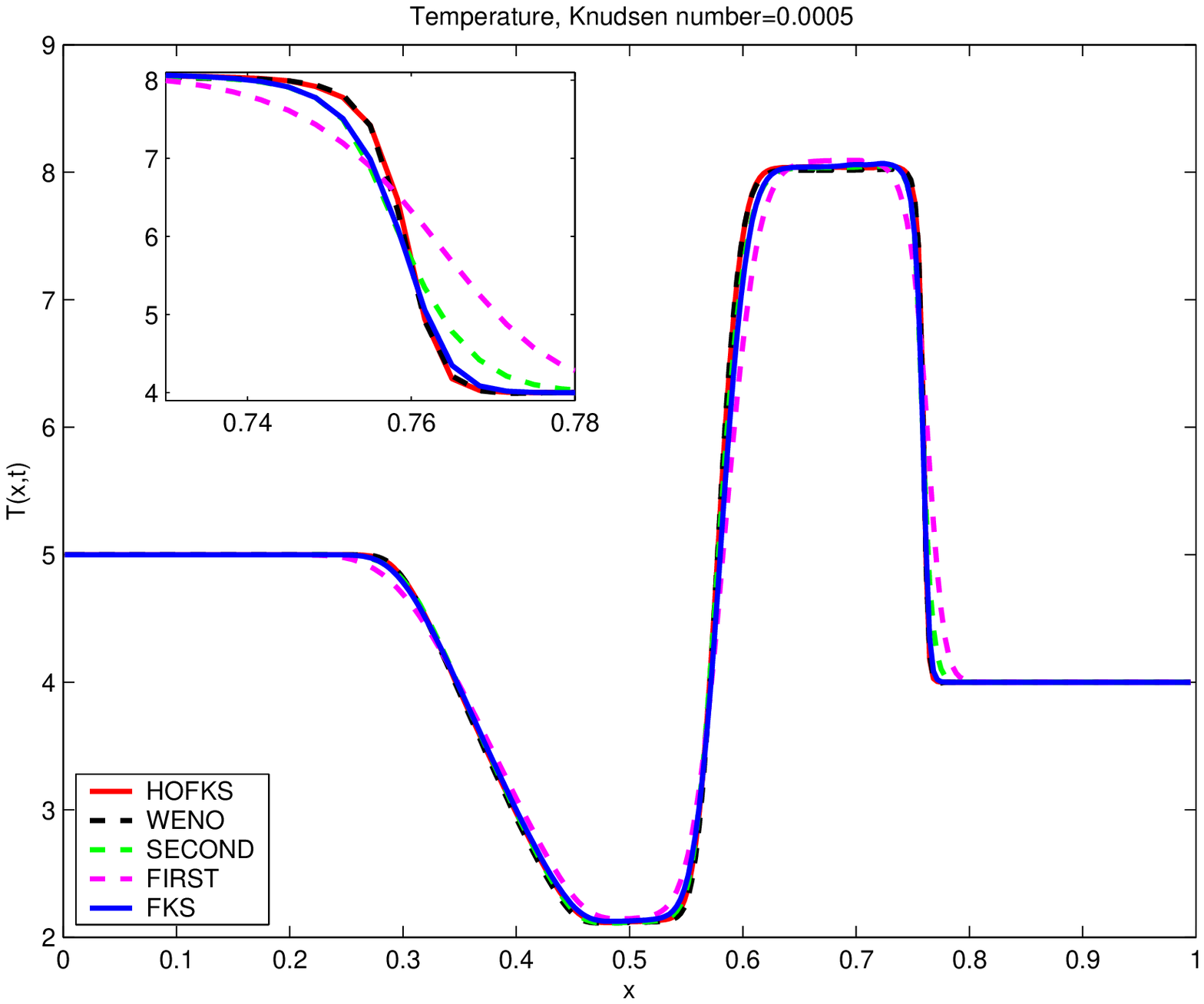}
\includegraphics[scale=0.43]{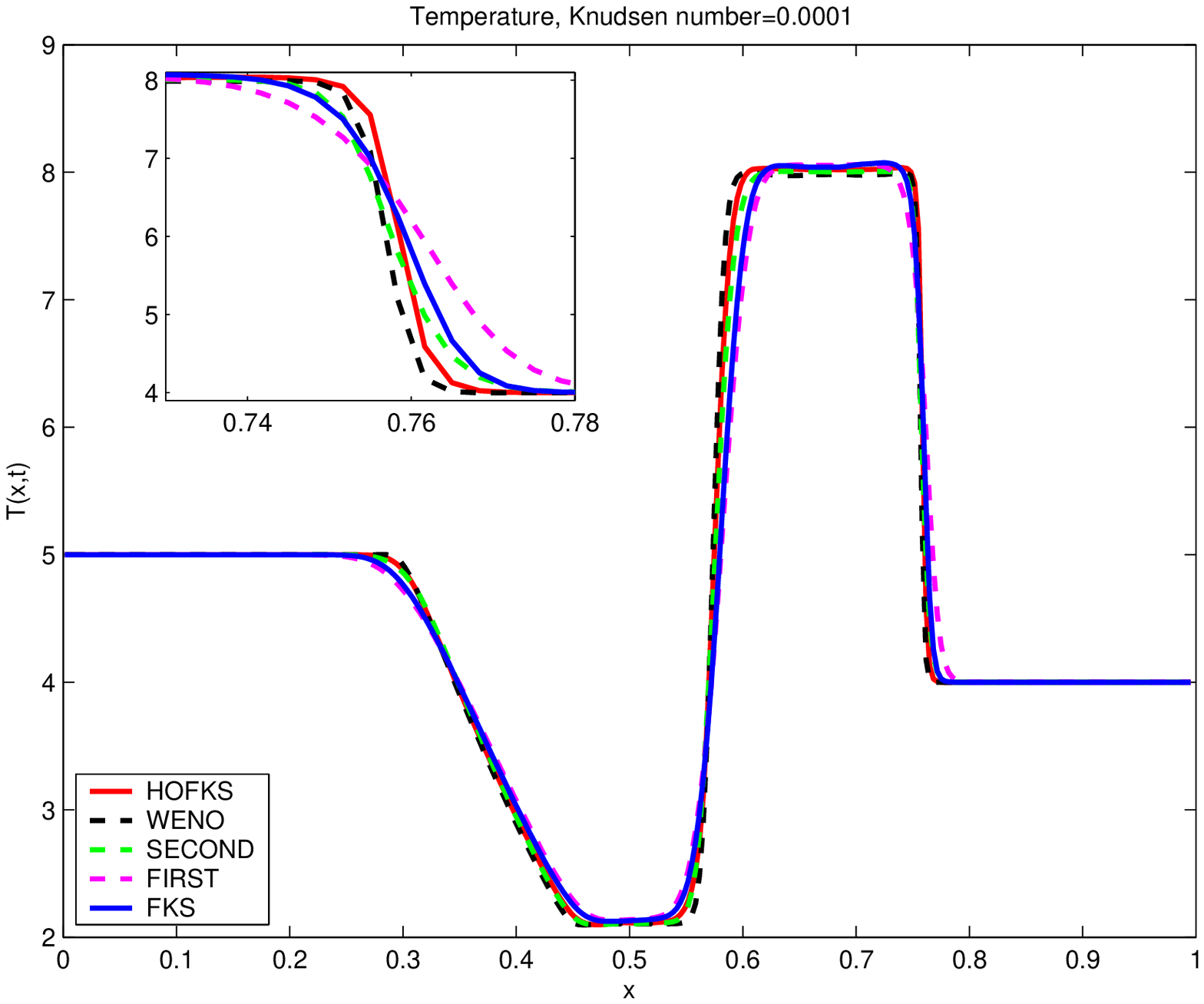}\\
\caption{1D Sod test: solution at $t_{\text{final}}=0.05$ for the
temperature, with $\tau=10^{-2}$ (top left), $\tau=10^{-3}$ (top
right), $\tau=5 \ 10^{-4}$ (bottom left) and $\tau=10^{-4}$ (bottom
right).} \label{sod3}
\end{center}
\end{figure}
%=== E N D   F I G ================================================

\subsection{1D Sod shock tube problem}
\label{subsec_sod}

We consider the 1D/1D Sod test with $300$ mesh points in physical
and $100$ points in velocity spaces. The boundaries in velocity
space are set to $-15$ and $15$. The left and right states are given
by a density $\rho_{L}=1$, mean velocity $u_{L}=0$ and temperature
$T_{L}=5$ if $0 \leq x \leq 0.5$, while $\rho_{R}=0.125$, $u_{R}=0$,
$T_{R}=4$ if $0.5 \leq x \leq 1$. The gas is, at the initial state, in thermodynamical
equilibrium. We repeat the same test with four different values of
the Knudsen number, \textit{i.e.} $\tau=10^{-2}$, $\tau=10^{-3}$, $\tau=5 \times 10^{-4}$ and $\tau=10^{-4}$.
We plot the results for the final time $t_{\text{final}}=0.05$ for
the density (Figure \ref{sod1}), the mean velocity (Figure
\ref{sod2}) and the temperature (figure \ref{sod3}). In each figure
we compare the HOFKS method with the FKS method. We reported also the solutions computed with a third order WENO method, a
second-order MUSCL method and a first-order upwind method
\cite{leveque:numerical-methods}. These numerical methods, used as
reference, employ the same discretization parameters, except for the
time steps which for each scheme is chosen in order to satisfy the stability conditions.

From Figures (\ref{sod1}) to (\ref{sod3}) we observe that the HOFKS, the FKS and the WENO methods give identical or almost identical results for $\tau=10^{-2}$, this result was expected. We build up the method in such a way that for larger $\tau$ it behaves like the original fast kinetic scheme, because we knew that in these regimes the FKS already gave very good results. For larger values of $\tau$ we found the same behaviors as for the case $\tau=10^{-2}$, thus we did not report simulations results. Starting from $\tau=10^{-3}$, some small differences arise between the HOFKS and the FKS methods, however both schemes are still very close to the WENO solution. For $\tau=5 \times 10^{-4}$, we clearly see differences between the high order fast kinetic scheme and the fast kinetic scheme. In particular, we see that the HOFKS remains stick to the third order WENO scheme while the FKS not. This aspect is made very clear for $\tau=10^{-4}$. In this latter case, the HOFKS gives very good results while the FKS scheme lays between a first and a second order space accurate method. The key point of the schemes developed is their very low CPU time consumption in comparison to other existing methods. This gain as expected is not so relevant for the one dimensional case, while it becomes very important for the two and the three dimensional cases. Thus, in the next subsection we report some two dimensional simulations together with their computational costs.

%\clearpage
\subsection{2D isentropic vortex}
\label{subsec_isentropicvortex}
This vortex test case is a classical test to assess the accuracy of numerical schemes because this
problem produces a genuine 2D smooth solution of Euler equations.
The overall error produced by our kinetic schemes associates spacial discretization error, time discretization error
and velocity discretization error.
In our case the first order time discretization dooms the overall scheme to remain first order accurate only.
Nevertheless using a time step $\Delta t$
of the order $\Delta x^2$ or even smaller allows to reduce the time discretization error to a negligible
quantity by respect to space discretization error. Doing so we can effectively measure the spacial accuracy of
the overall scheme when the velocity error is kept small enough. \\
To do so an isentropic vortex is introduced to a uniform mean flow, by small perturbations of velocity, density and temperature
variables and is detailed in \cite{Shu98}, \cite{YVD} as instance.
The simulation domain $\Omega$ is the square
$[0,10]\times [0,10]$ and we consider an initial gas flow given by the following background condition
$\rho_\infty=1.0$, $u_\infty=1.0$, $v_\infty=1.0$, $p_\infty=1.0$, with a normalized ambient temperature
$T^*_\infty = 1.0$ computed with the perfect gas equation of state and $\gamma=5/3$.
%=== B E G I N   F I G ================================================
\begin{figure}[h!]
\begin{center}
   %\begin{minipage}[c]{.63\textwidth}
   %   \includegraphics[width=1.05\textwidth]{figure/density_init}
   %\end{minipage} \hfill
   %\begin{minipage}[c]{.36\textwidth}
   %   \includegraphics[width=1.05\textwidth]{figure/velvec_init} \\
   %   \includegraphics[width=1.05\textwidth]{figure/vortex_init}
   %\end{minipage}
  \includegraphics[width=0.4\textwidth, height=0.4\textwidth]{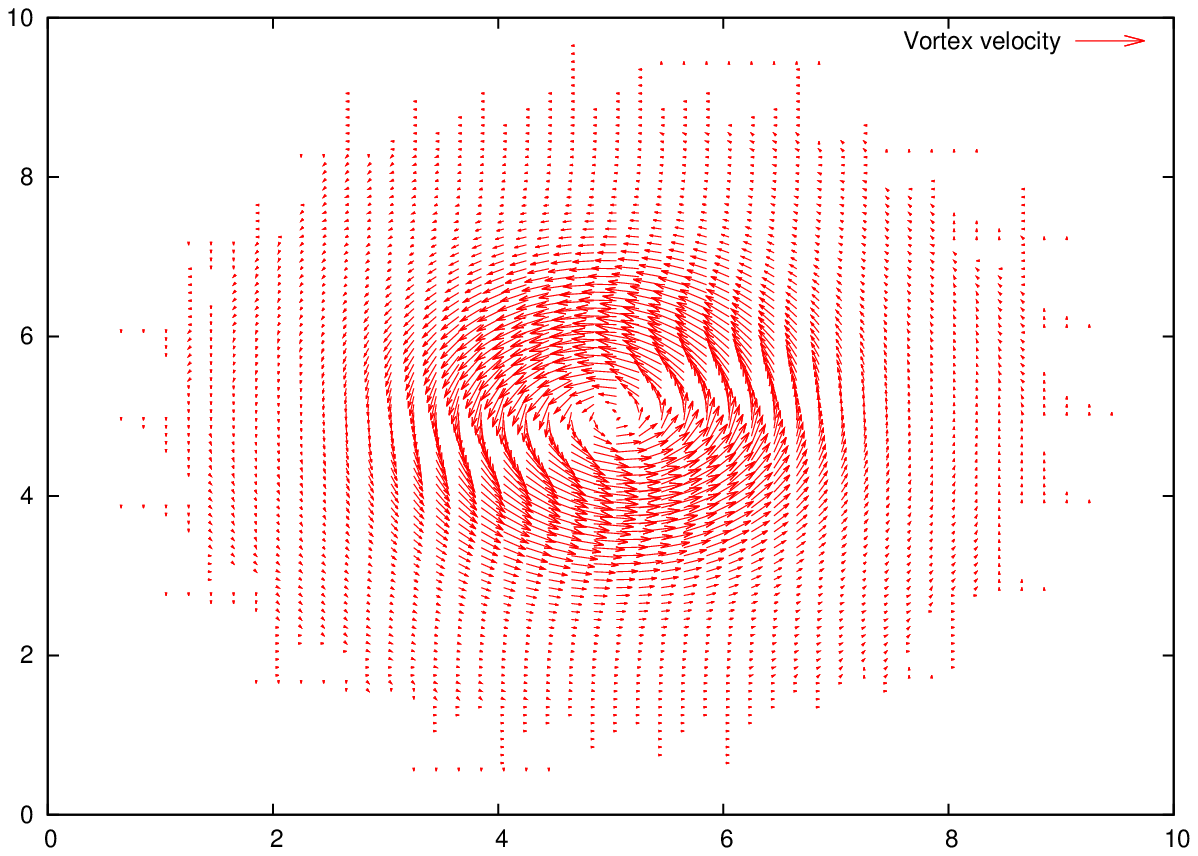}
  \includegraphics[width=0.4\textwidth, height=0.35\textwidth]{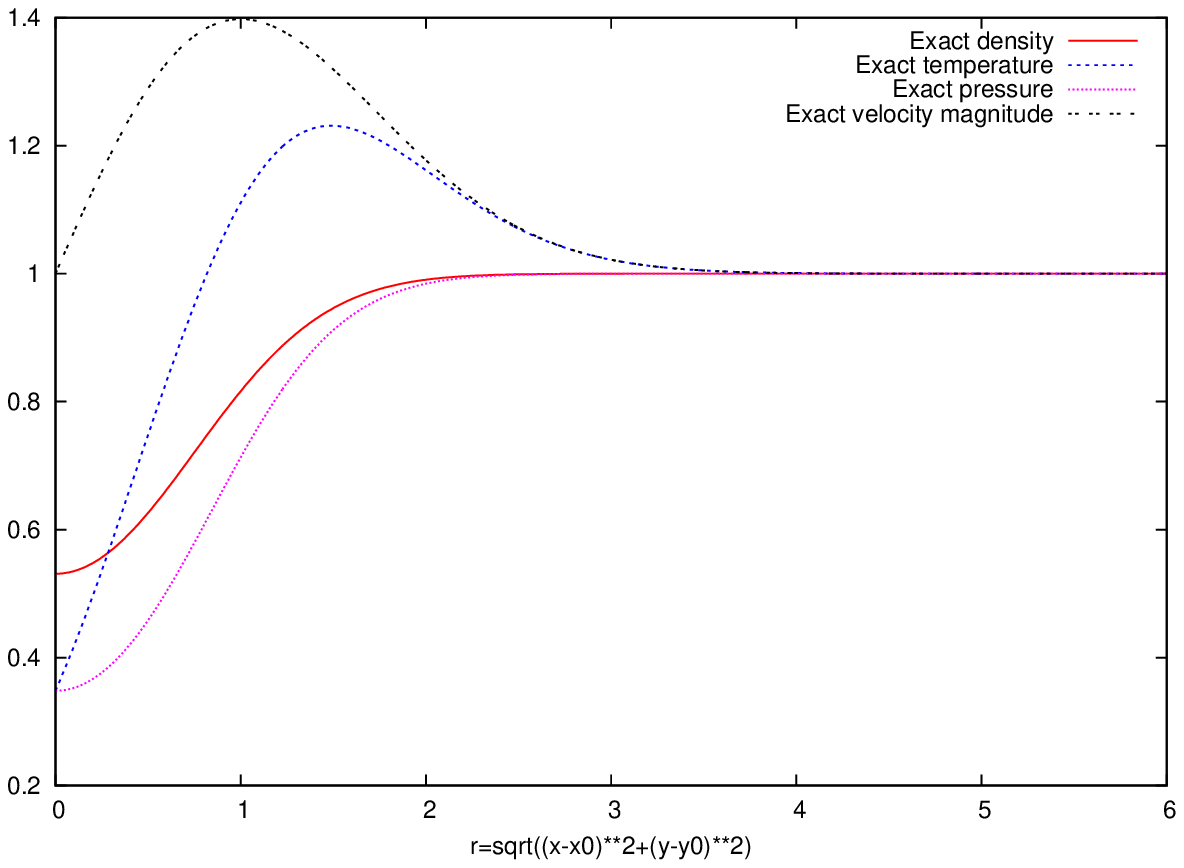}
   \caption{
     2D isentropic vortex. Initial data.
     %Left: density. Right-top:
     Vortex velocity vectors (background velocity is substracted) and
     exact density, temperature, pressure and velocity magnitude as a function of $r= \sqrt{{x'}^2 + {y'}^2}$, $ (x' = x - x_{0}, y' = y - y_{0})$
     with $(x_0,y_0)$ the center of the vortex.
   } \label{fig:vortex}
\end{center}
\end{figure}
%=== E N D   F I G ================================================

\begin{figure}[h!]
  \begin{center}
    \begin{tabular}{cc}
      \includegraphics[width=0.45\textwidth]{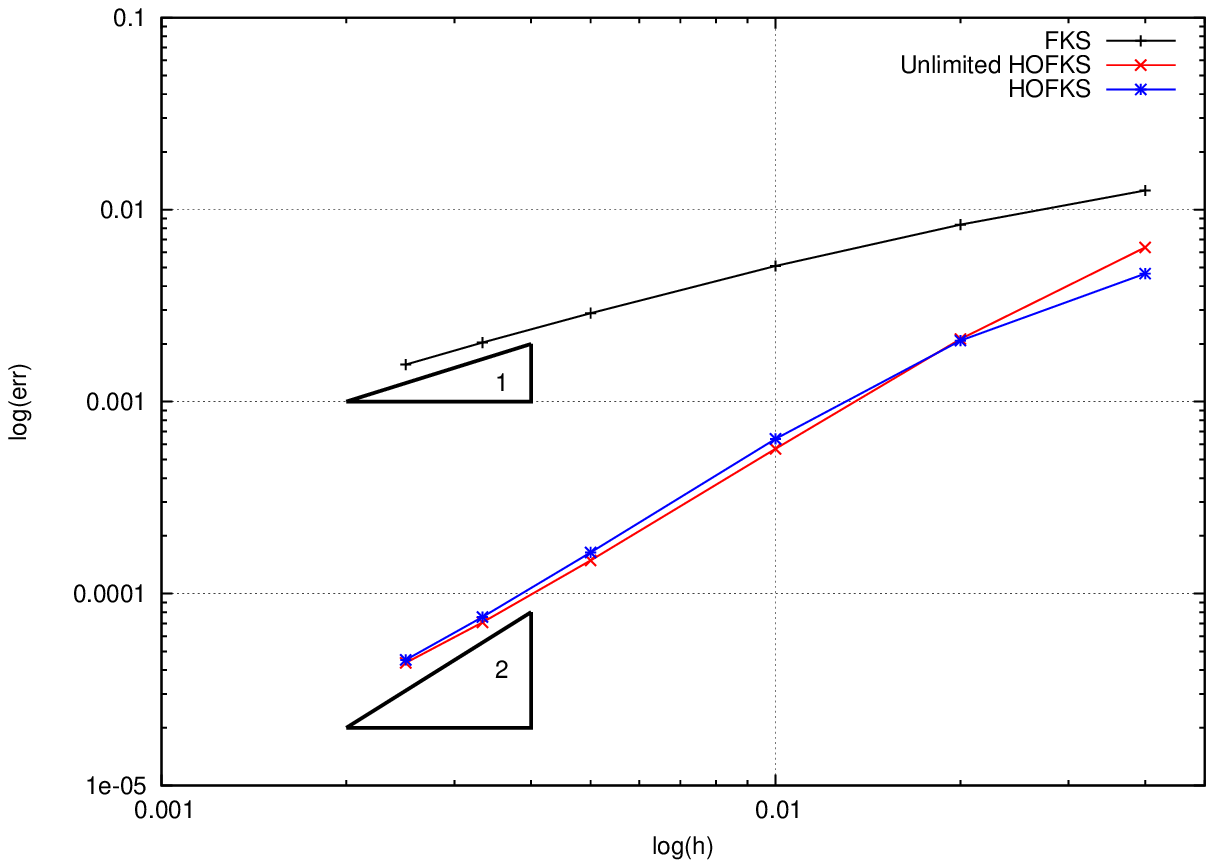} &
      \includegraphics[width=0.45\textwidth]{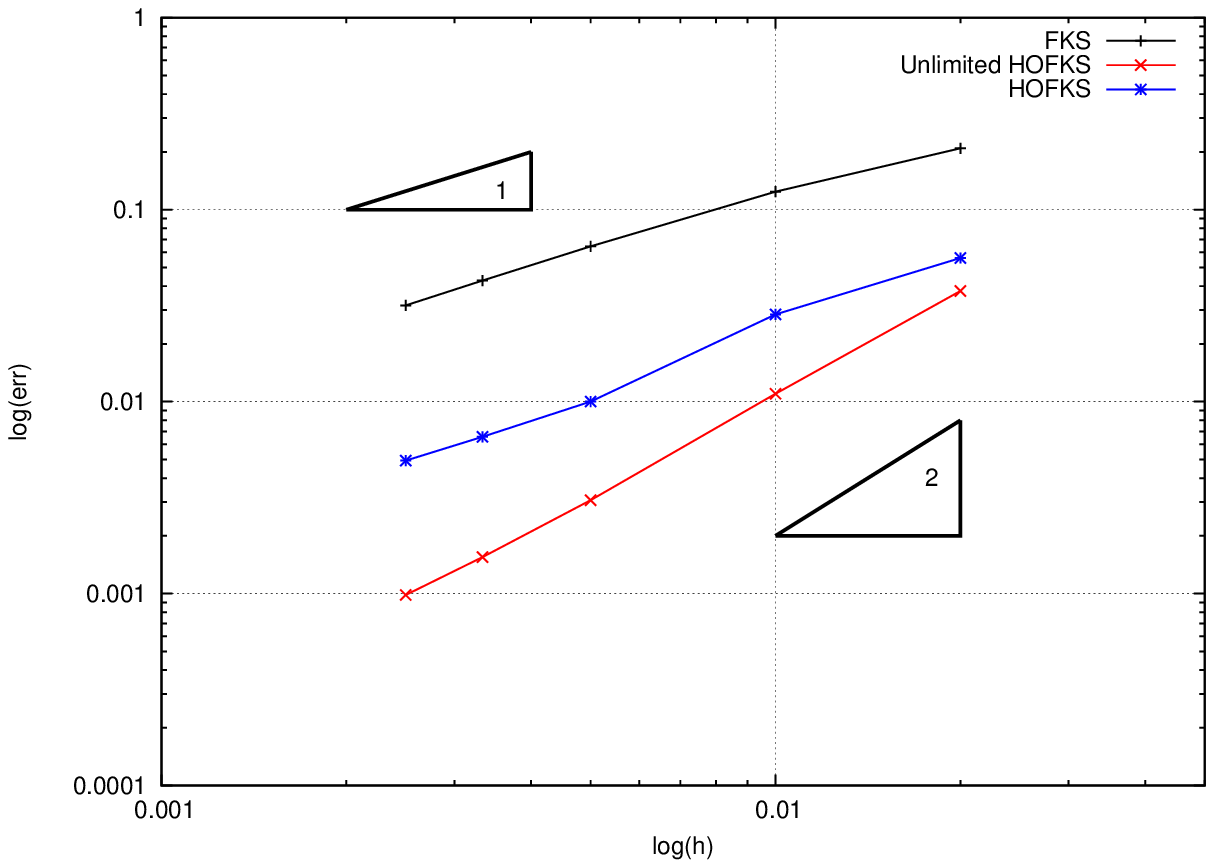}
    \end{tabular}
    \caption{ \label{fig:convergence_isentropic} $L^1$ and $L^\infty$ errors in logscale for the vortex problem. FKS,
    unlimited HOFKS and HOFKS are presented.}
  \end{center}
\end{figure}

A vortex centered at $X_0 =(x_0, y_0)=(0,0)$
is added to the ambient gas at the initial time $t=0$ with the following conditions
$u = u_\infty + \delta u$, $ v = v_\infty + \delta v$, and $ T^* = T^*_\infty + \delta T^*$
\begin{eqnarray}
  \nonumber \delta u   = -y' {\frac {\beta} {2 \pi}} \exp \left( {\frac {1-r^2} {2}} \right),  \quad
  \delta v   = x' {\frac {\beta} {2 \pi}} \exp \left( {\frac {1-r^2} {2}} \right), \quad
%& \\
%  \nonumber
  \delta T^* = - { \frac {(\gamma - 1 ) \beta} {8 \gamma \pi^2}} \exp \left( {1-r^2} \right).&
\end{eqnarray}
with $ r = \sqrt{{x'}^2 + {y'}^2}$, $ (x' = x - x_{0}, y' = y - y_{0})$
and vortex strength is given by $ \beta=5.0$.
Consequently, the initial density is given by
\begin{eqnarray}
  \rho
  =   \rho_\infty \left( {\frac{T^*}{T_\infty^*} } \right)^{\frac{1}{\gamma-1} }
  = \left(1 - { \frac {(\gamma - 1 ) \beta} {8 \gamma \pi^2}} \exp \left( {1-r^2} \right)
  \right)^{\frac{1}{\gamma-1} },
\end{eqnarray}
and the pressure ig given by $p=\rho^\gamma$.
We assume periodic conditions and the exact solution at any time $T>0$ is the same vortex but translated
in the direction $V_\infty=(u_\infty,v_\infty)$.
Note that $V_\infty=(0,0)$ generates a static vortex which is usually simpler to solve and can also be misleading.
The exact density function for any point at time $T$ is denoted by $\rho^{ex}(x,T)$, moreover in figure~\ref{fig:vortex} are plotted the exact
solution for density, temperature, pressure and velocity as a function of $r$.
A series of refined meshes (from $25\times 25$ up to $400\times 400$ cells) are successively used to compute the numerical solution where $\tau=10^{-4}$ is used.
The meshes are made of $20^{2}$ points in velocity space with bounds $[-15,15]^{2}$.
The errors at time $T=t^n$ for $M$ spatial cells in one direction are given by
\begin{equation}
\varepsilon_M^1     = \frac{\Sum_{j\in\mathcal{J}}  | \rho^{ex}(x_j,t^n) - \rho_j^n |}{\Sum_{j\in\mathcal{J}}  | \rho^{ex}(x_j,t^n) |}, \ \ \ \
\varepsilon_M^\infty = \frac{\max_{j\in\mathcal{J}}  | \rho^{ex}(x_j,t^n) - \rho_j^n |}{\max_{j\in\mathcal{J}}  | \rho^{ex}(x_j,t^n) |},
\end{equation}
The rates of convergence are computed as $\log(\varepsilon_{M'}/\varepsilon_M)/\log(M'/M)$ for two meshes with $M'$ and $M$ cells.
In Table \ref{tab:isentropic_vortex_Error} are gathered the $L^1$ and $L^\infty$ errors on the density variables and rates of convergence for
FKS, the unlimited HOFKS and HOFKS at final time $T=1$. As expected the FKS produces only first order accurate results in $L^1$ and $L^\infty$ norms.
Contrarily the unlimited HOFKS can reach a genuine higher order of accuracy in both norms; the high accuracy in $L^\infty$ norm
is due to the fact that the solution is smooth and no limiter is applied therefore extrema are only little diffused compared to limited schemes.
Finally the (limited) HOFKS behaves, as expected, like a high order accurate scheme in $L^1$ norm and like a first order scheme in $L^\infty$ norm.
%-------------------------------
 \begin{table}[h]
 \centering
 \numerikSeven
 \begin{tabular}{|cc|cc|cc|cc|cc|cc|cc|}
 \hline
 \multicolumn{14}{|c|}{\textbf{Advected isentropic vortex}}\\
 \cline{3-14}
 & & \multicolumn{4}{c|}{\textbf{FKS}} & \multicolumn{4}{c|}{\textbf{Unlim. HOFKS}} & \multicolumn{4}{c|}{\textbf{HOFKS}} \\
 \hline
 $\Delta x$ & $M$ &  $L^1$ err & &$L^\infty$ err&  & $L^1$ err & &$L^\infty$ err & & $L^1$ err &  &$L^\infty$ err & \\
 \hline \hline
 1/25   &  25$^2$ & 1.26E-02 &  ---   & 2.78E-01  &  ---  & 6.36E-03  & ---  & 1.16E-01  &  --- & 4.64E-03 &  --- & 9.21E-02 & \\
 1/50   &  50$^2$ & 8.36E-03 &  0.59  & 2.07E-01  &  0.43 & 2.12E-03  & 1.59 & 3.77E-02  & 1.62 & 2.08E-03 & 1.16 & 5.60E-02 & 0.72\\
 1/100  & 100$^2$ & 5.09E-03 &  0.72  & 1.22E-01  &  0.76 & 5.68E-04  & 1.90 & 1.10E-02  & 1.78 & 6.40E-04 & 1.70 & 2.85E-02 & 0.97\\
 1/200  & 200$^2$ & 2.86E-03 &  0.82  & 6.34E-02  &  0.95 & 1.49E-04  & 1.93 & 3.07E-03  & 1.84 & 1.64E-04 & 1.97 & 1.00E-02 & 1.51\\
 1/300  & 300$^2$ & 2.03E-03 &  0.88  & 4.27E-02  &  1.01 & 7.08E-05  & 1.83 & 1.55E-03  & 1.69 & 7.55E-05 & 1.91 & 6.56E-03 & 1.04\\
 1/400  & 400$^2$ & 1.56E-04 &  0.91  & 3.17E-02  &  1.03 & 4.36E-05  & 1.69 & 9.85E-04  & 1.57 & 4.52E-05 & 1.79 & 4.93E-03 & 0.99\\
 \cline{1-14}
 \multicolumn{2}{|l|}{Expected order} &  & 1 &  & 1 &  & 2 &  & 2 & & 2 & & 1 \\
 \hline  													
%MOVING VORTEX --- phi=0 FV SCHEME – dt/2									
%25	1.26E-002			1.27E-003			2.80E-001		
%50	8.36E-003	1.50	0.59	6.03E-004	2.10	1.07	2.09E-001	1.34	0.42
%100	5.09E-003	1.64	0.72	2.22E-004	2.72	1.44	1.24E-001	1.68	0.75
%200	2.89E-003	1.76	0.81	6.85E-005	3.23	1.69	6.44E-002	1.93	0.95
%300	2.03E-003	1.43	0.88	3.29E-005	2.08	1.81	4.27E-002	1.51	1.01
%400	1.56E-003	1.85	0.91	1.92E-005	3.56	1.86	3.17E-002	2.03	1.03	
%MOVING VORTEX --- phi=1 UNLIMITED SCHEME HOFKS – dt/2									
%25	6.36E-003			2.75E-004			1.16E-001		
%50	2.12E-003	3.01	1.59	2.84E-005	9.67	3.27	3.77E-002	3.07	1.62
%100	5.68E-004	3.73	1.90	2.05E-006	13.86	3.79	1.10E-002	3.43	1.78
%200	1.49E-004	3.81	1.93	1.45E-007	14.19	3.83	3.07E-003	3.57	1.84
%300	7.08E-005	2.10	1.83	3.35E-008	4.31	3.60	1.55E-003	1.98	1.69
%400	4.36E-005	3.42	1.69	1.28E-008	11.26	3.34	9.85E-004	3.12	1.57
%	
%MOVING VORTEX --- LIMITED SCHEME HOFKS – dt/2									
%25	4.64E-003			1.42E-004			9.21E-002		
%50	2.08E-003	2.24	1.16	3.44E-005	4.14	2.05	5.60E-002	1.65	0.72
%100	6.40E-004	3.24	1.70	4.44E-006	7.75	2.95	2.85E-002	1.96	0.97
%200	1.64E-004	3.92	1.97	3.56E-007	12.45	3.64	1.00E-002	2.85	1.51
%300	7.55E-005	2.17	1.91	8.27E-008	4.31	3.60	6.56E-003	1.53	1.04
%400	4.52E-005	3.62	1.79	2.82E-008	12.65	3.74	4.93E-003	2.03	0.99									
 \end{tabular}
 \caption{ \label{tab:isentropic_vortex_Error}
   $L^1$ and $L^\infty$ errors and convergence rates for the
   isentropic vortex problem with FKS, unlimited HOFKS and HOFKS schemes. }
 \end{table}
%-------------------------------
We also display in figure \ref{fig:convergence_isentropic} the convergence curves corresponding to the errors of Table
\ref{tab:isentropic_vortex_Error}.

%-------------------------------
%---------------------------------

\begin{figure}[h!]
  \begin{center}
    \includegraphics[width=0.75\textwidth]{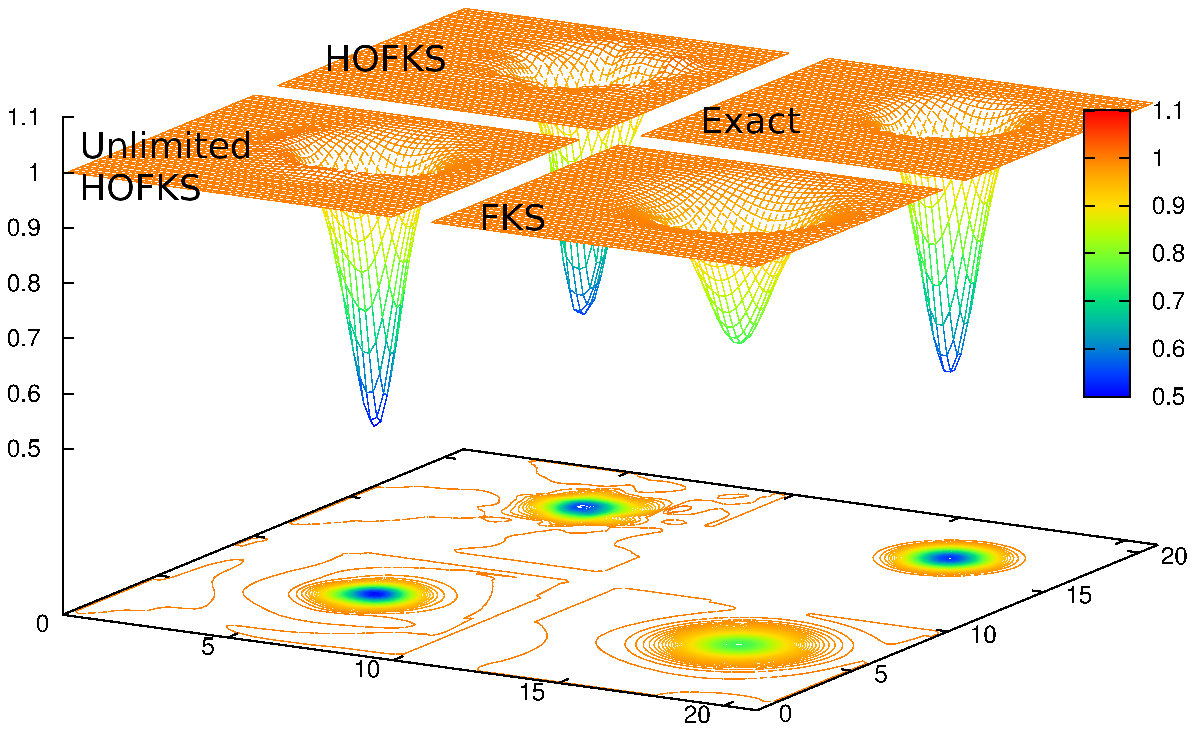} \\
    \begin{tabular}{cc}
      \includegraphics[width=0.5\textwidth]{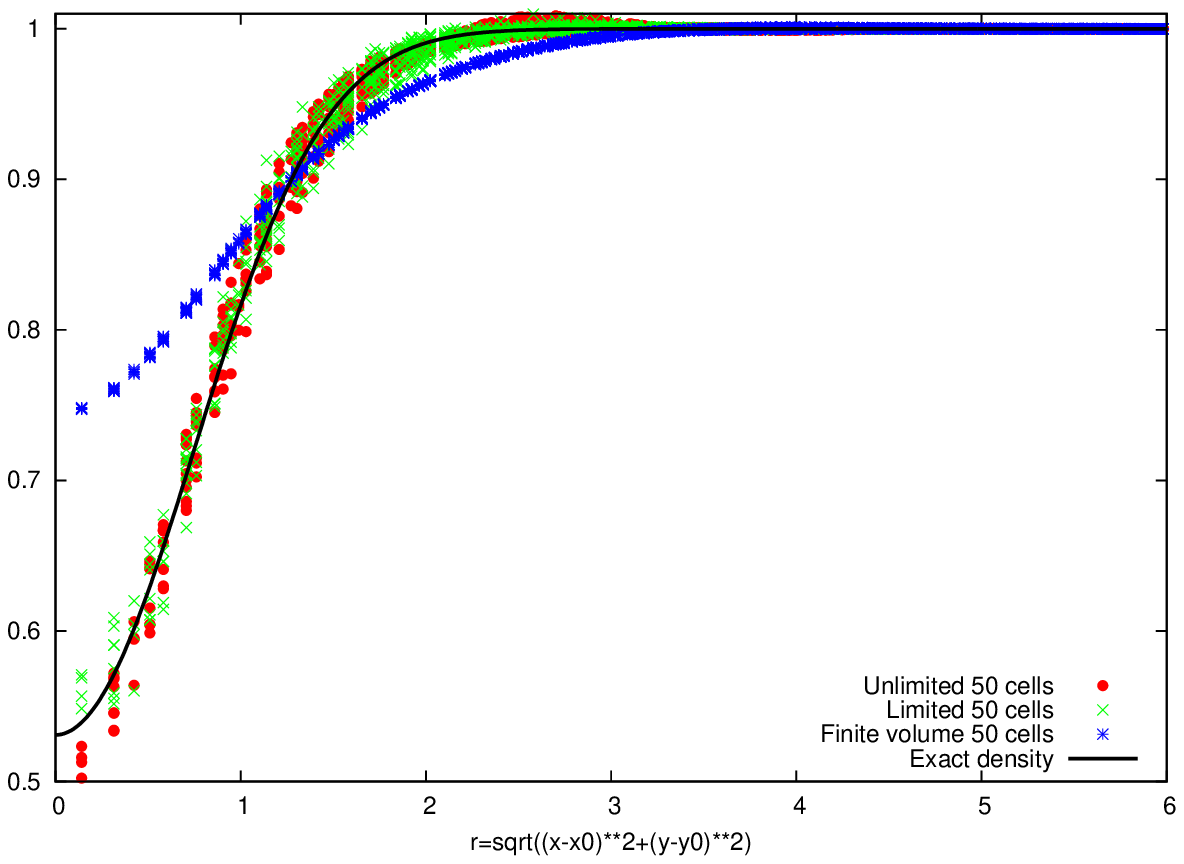} &
      \includegraphics[width=0.5\textwidth]{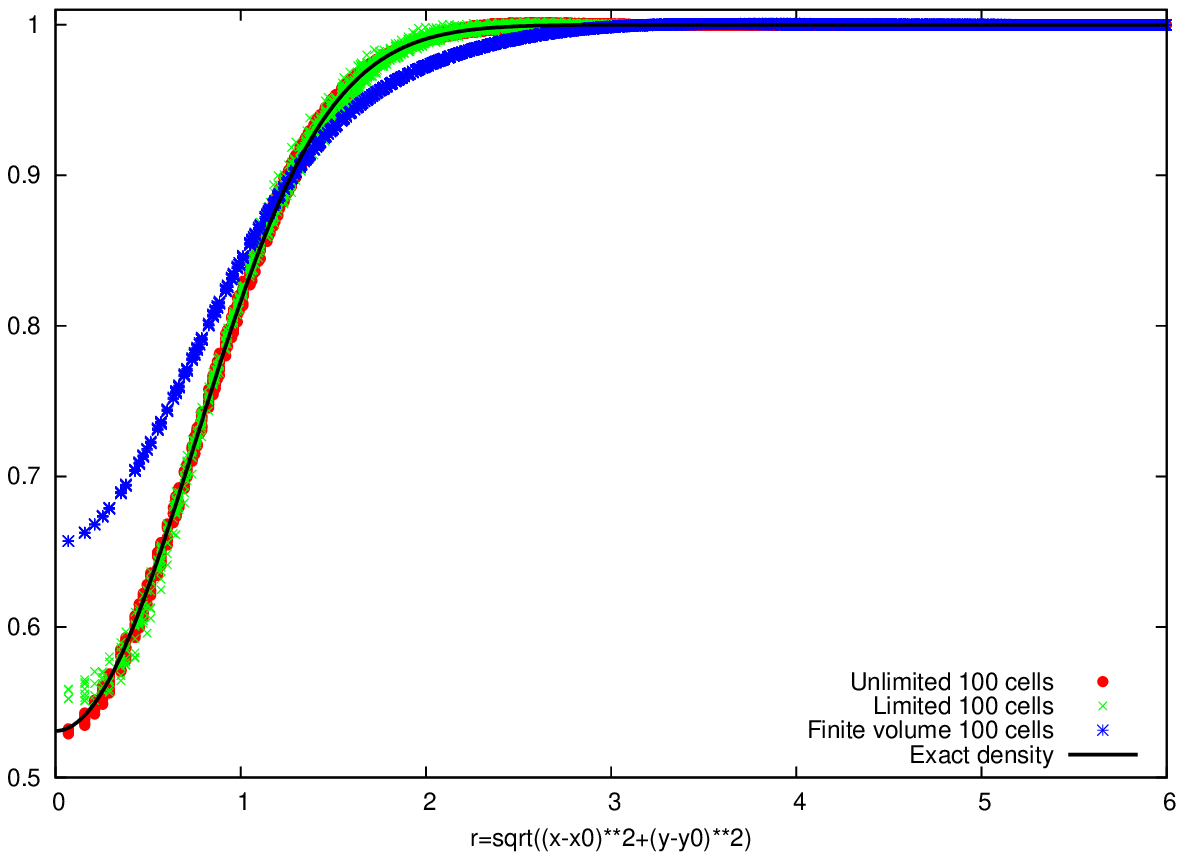} \\
    \end{tabular}
    \caption{ \label{fig:solution_vortex}
      2D Vortex problem at $t=1$. Density results obtained by the three schemes: FKS, unlimited HOFKS and HOFKS.
      Top: 3D view of the density on the $50\times 50$ mesh versus the exact solution.
      Bottom:
      $50\times 50$ (left panel) and $100 \times 100$ (right panel) cells meshes, density as a function of
      $r= \sqrt{{x'}^2 + {y'}^2}$, $(x' = x - x_{0}(T), y' = y - y_{0}(T))$,
      where $x_0(T)=x_0+ u\ T$ and $y_0(T)=y_0 + v \ T$ are the exact coordinates of the vortex center, is plotted for all cells.
    }
  \end{center}
\end{figure}

Finally in figure \ref{fig:solution_vortex} we present the density results obtained by the three schemes: FKS, unlimited HOFKS and HOFKS.
 The top panel presents the density on the $50\times 50$ mesh for all schemes versus the exact solution. The bottom panels present
 the $50\times 50$ and $100 \times 100$ cell mesh results: the density as a function of
$r= \sqrt{{x'}^2 + {y'}^2}$, $(x' = x - x_{0}(T), y' = y - y_{0}(T))$,
where $x_0(T)=x_0+ u\ T$ and $y_0(T)=y_0 + v \ T$ are the exact coordinates of the vortex center at final time $T$, is plotted for all cells
in the domain. Doing so we can measure the ``convergence'' of the results when a finer mesh is used, the excessive diffusion of the first
order scheme and the tendency of undershooting of the unlimited scheme.

%\clearpage
\subsection{2D Sod shock tube problem}
\label{subsec_sod2D}
%-------------------------------
%---------------------------------

We consider now the 2D/2D Sod test on a square $[0,2]\times[0,2]$.
The velocity space is also a square with bounds $-15$ and $15$,
\textit{i.e.} $[-15,15]^{2}$, discretized with $N_v=20$ points in
each direction which gives $20^{2}$ points. The domain is divided into two parts,
a disk centered at point $(1,1)$ of radius $R_d=0.2$ is filled with a
gas with density $\rho_{L}=1$, mean velocity $u_{L}=0$ and
temperature $T_{L}=5$, whereas the gas in the rest of the domain
is initiated with $\rho_{R}=0.125$, $u_{R}=0$, $T_{R}=4$.
The final time is $t_{\text{final}}=0.07$.

We report results for two different values of the Knudsen number $\tau=10^{-3}$ and $\tau=10^{-4}$. For these regimes, we can appreciate the differences between the high order fast kinetic scheme and the fast kinetic scheme. For larger $\tau$, as expected, the solutions furnished by the two methods are very close and thus we do not show the figures. In the case of $\tau=10^{-4}$
we compare the results between the HOFKS and the FKS method with a first order and a second order MUSCL scheme for the compressible Euler equations.
In the case of $\tau=10^{-3}$ we compare the results between the HOFKS and the FKS method with a DSMC method for the Boltzmann-BGK equation.
For this latter, we employed on average $100$ particles per cell and we averaged the solution over $100$ realizations.

\begin{figure}[h!]
\begin{center}
\includegraphics[scale=0.4]{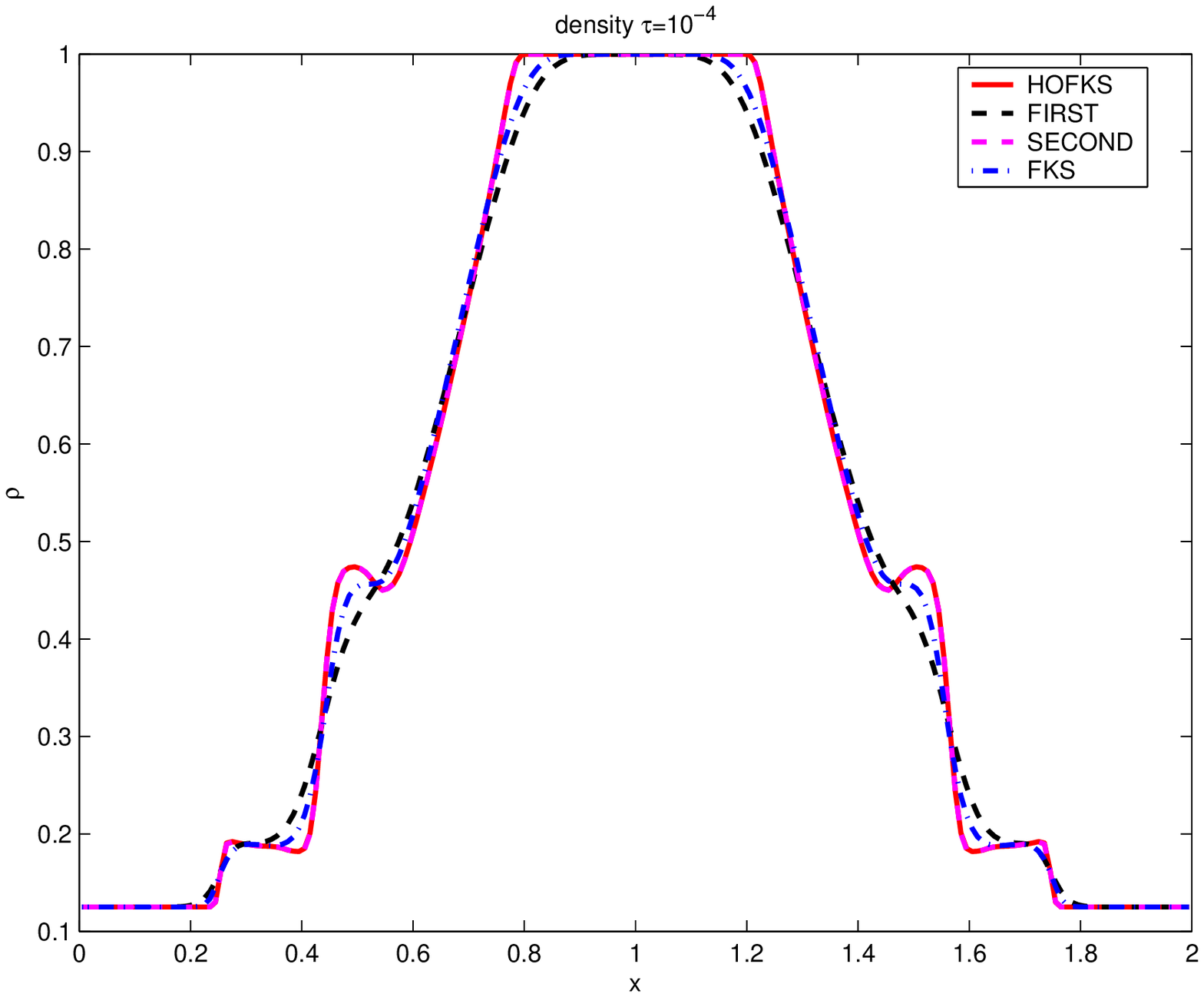}
\includegraphics[scale=0.4]{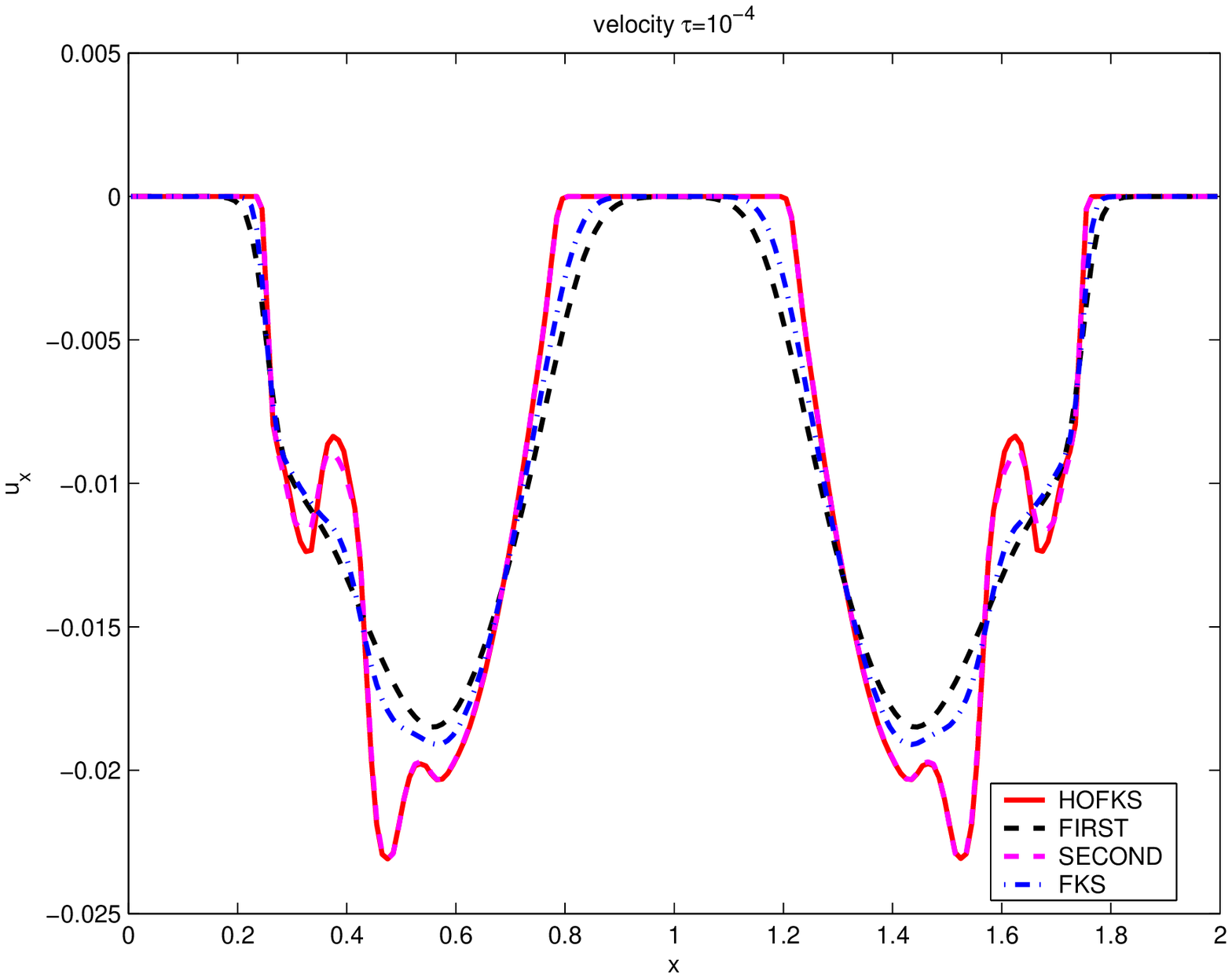}\\
\includegraphics[scale=0.4]{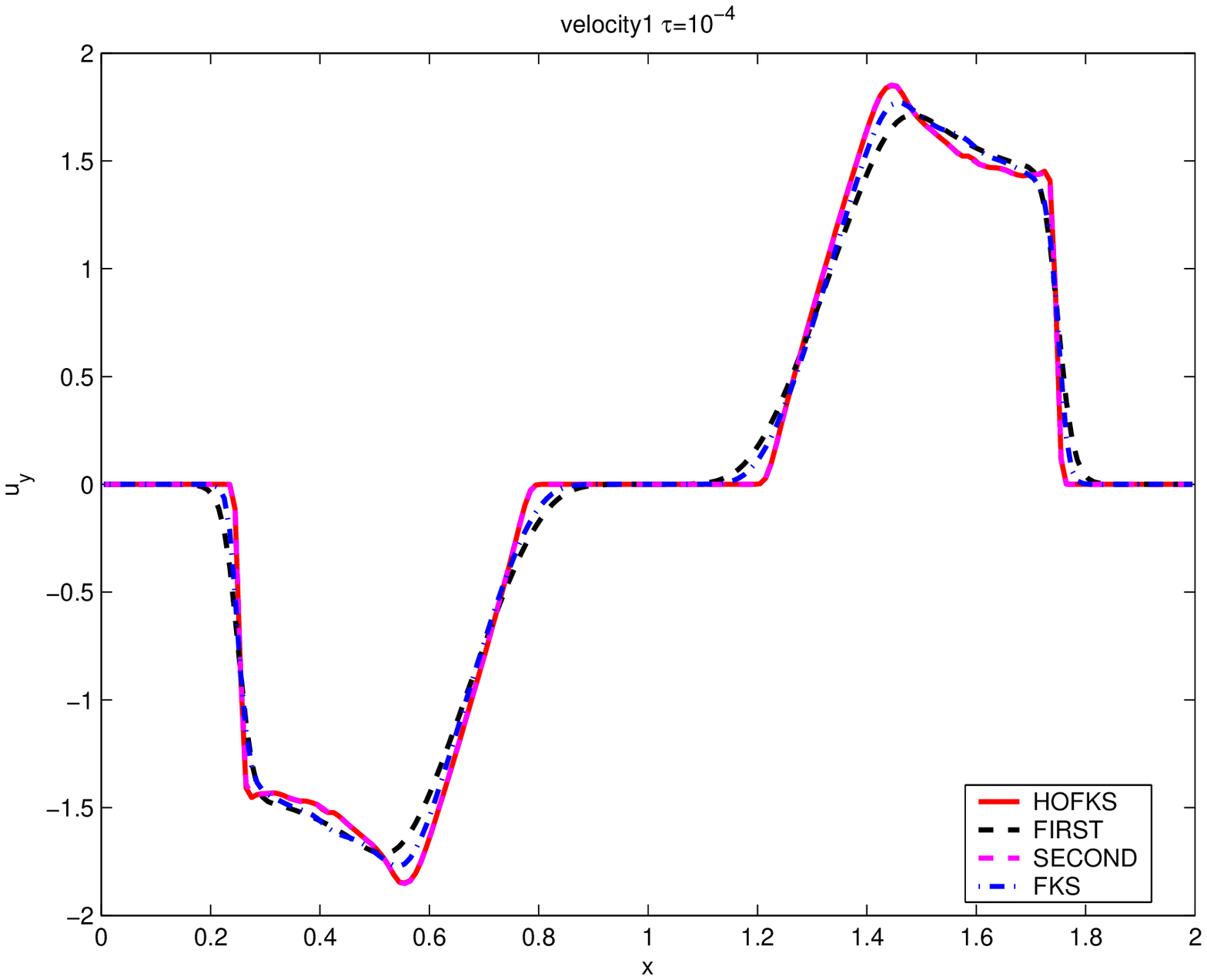}
\includegraphics[scale=0.4]{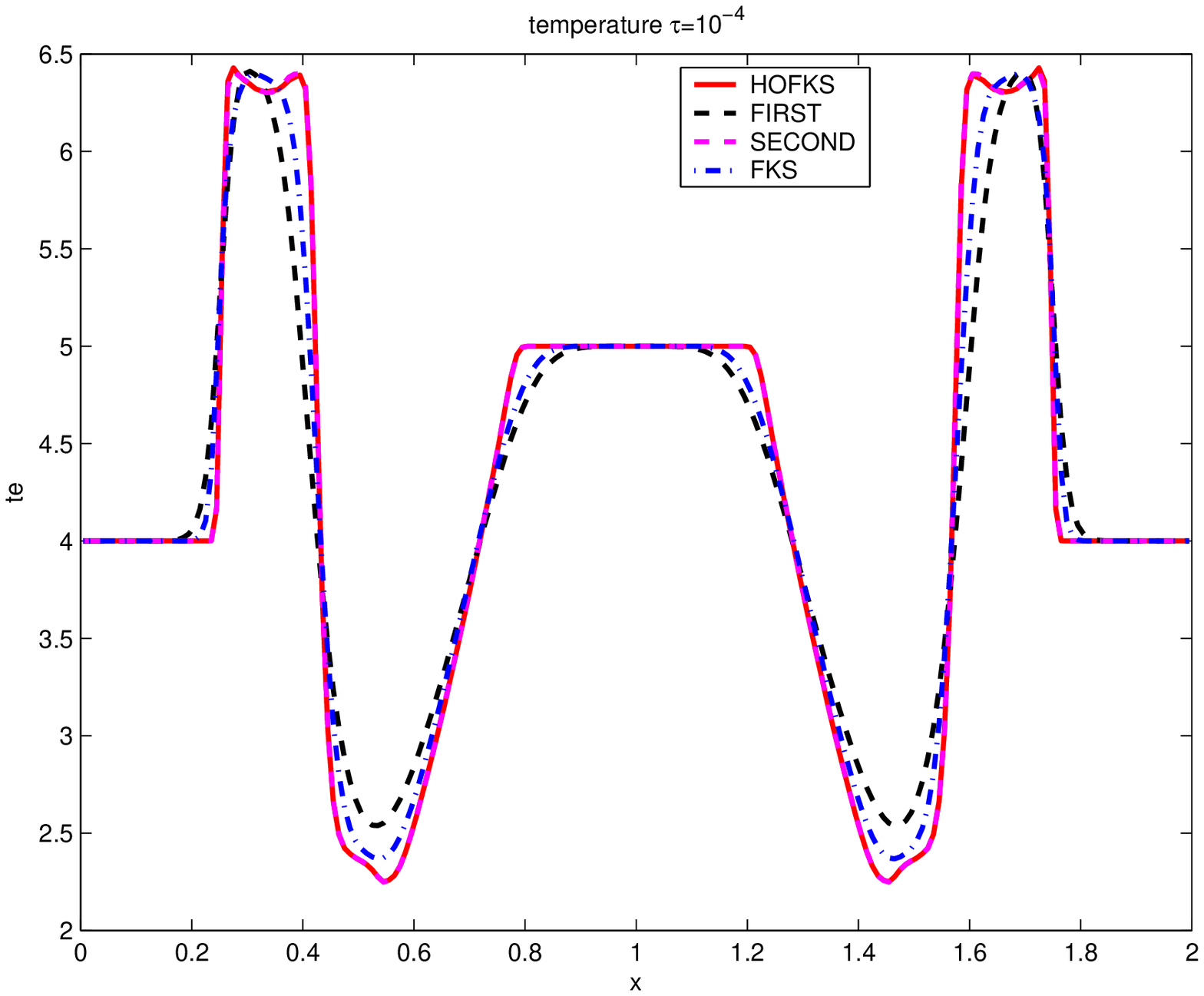}\\
\includegraphics[scale=0.45]{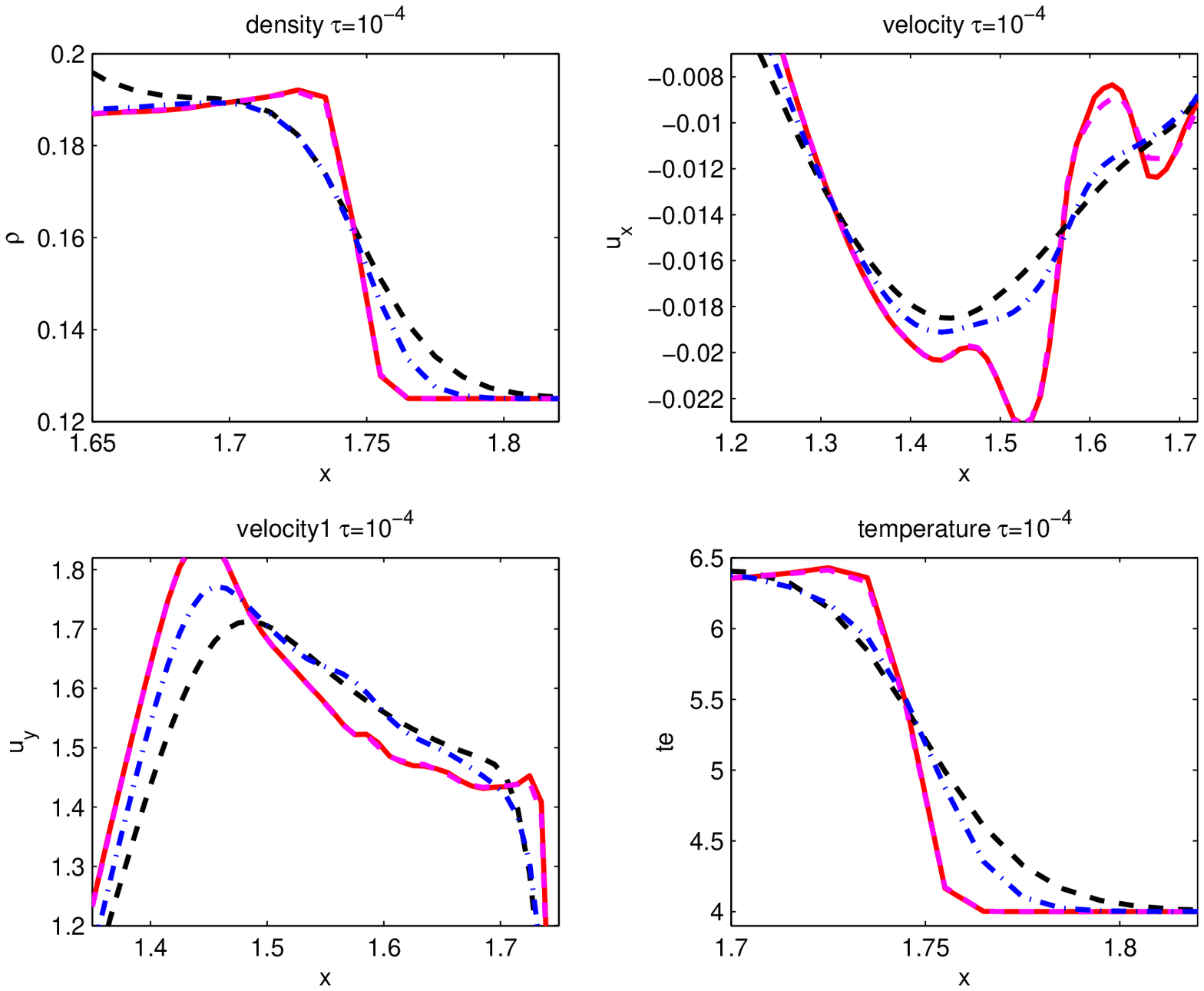}
\caption{2D Sod test: solution at
$t_{\text{final}}=0.07$ and $x=1$ for $\tau=10^{-4}$.
Top-Middle: Density (top left), velocity in the $x$-direction (top right), velocity in the
$y$-direction (middle left) and temperature (middle right).
Bottom: magnification of the solution.
HOFKS method continuous line, FKS method dash dotted line, first order and second order MUSCL methods dotted
lines.} \label{sod2D1}
\end{center}
\end{figure}

\begin{figure}[h!]
\begin{center}
\includegraphics[scale=0.4]{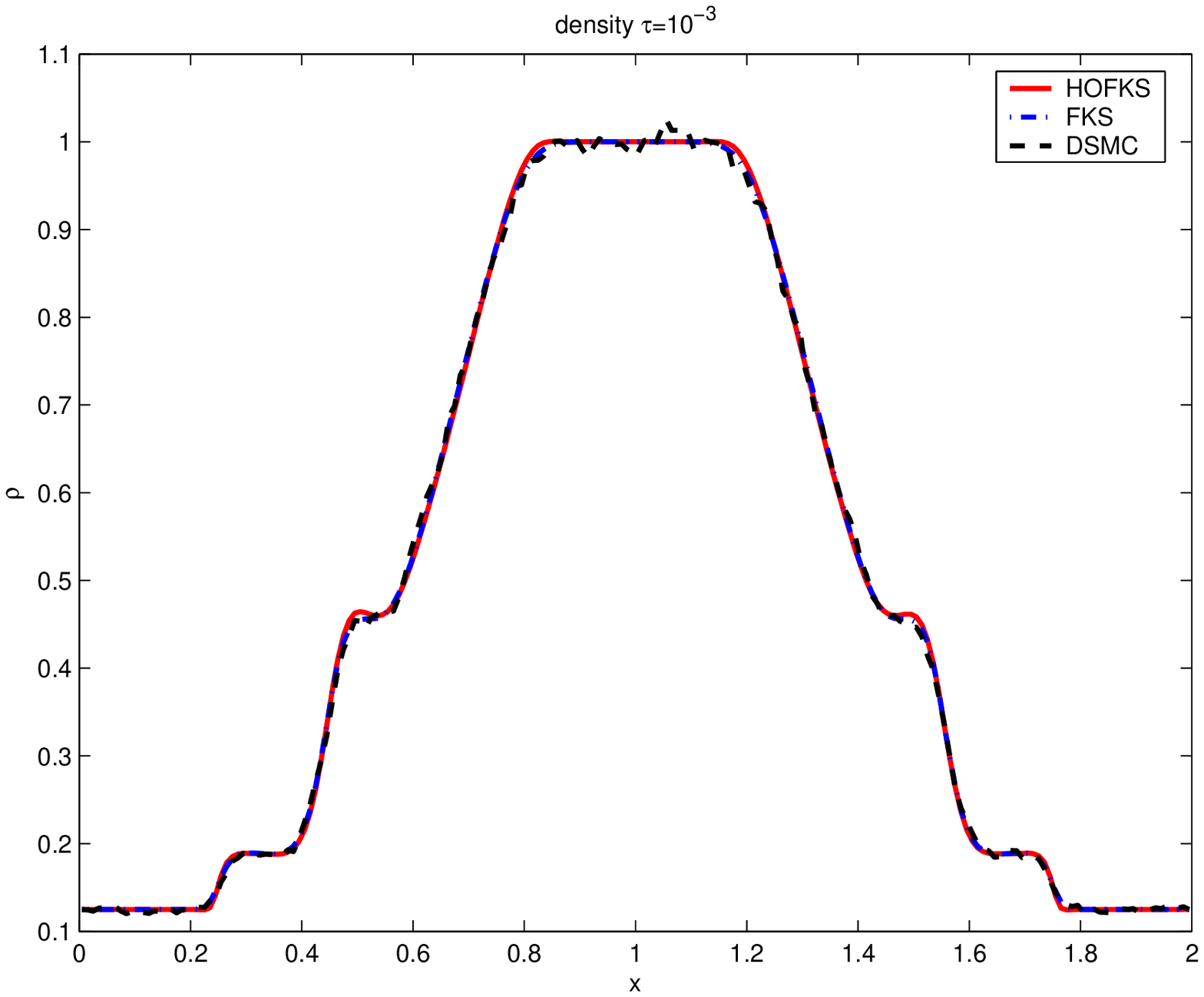}
\includegraphics[scale=0.4]{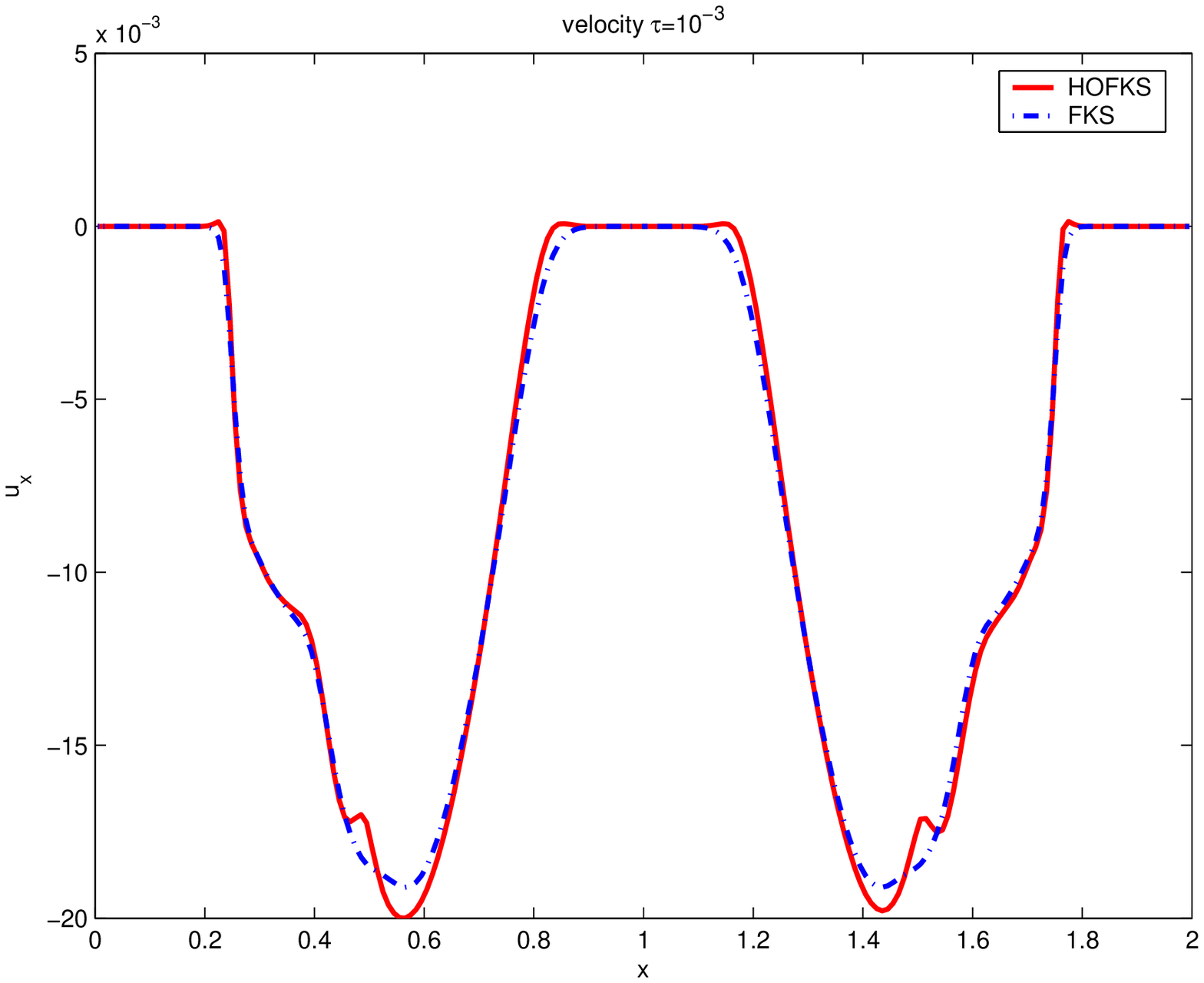}\\
\includegraphics[scale=0.4]{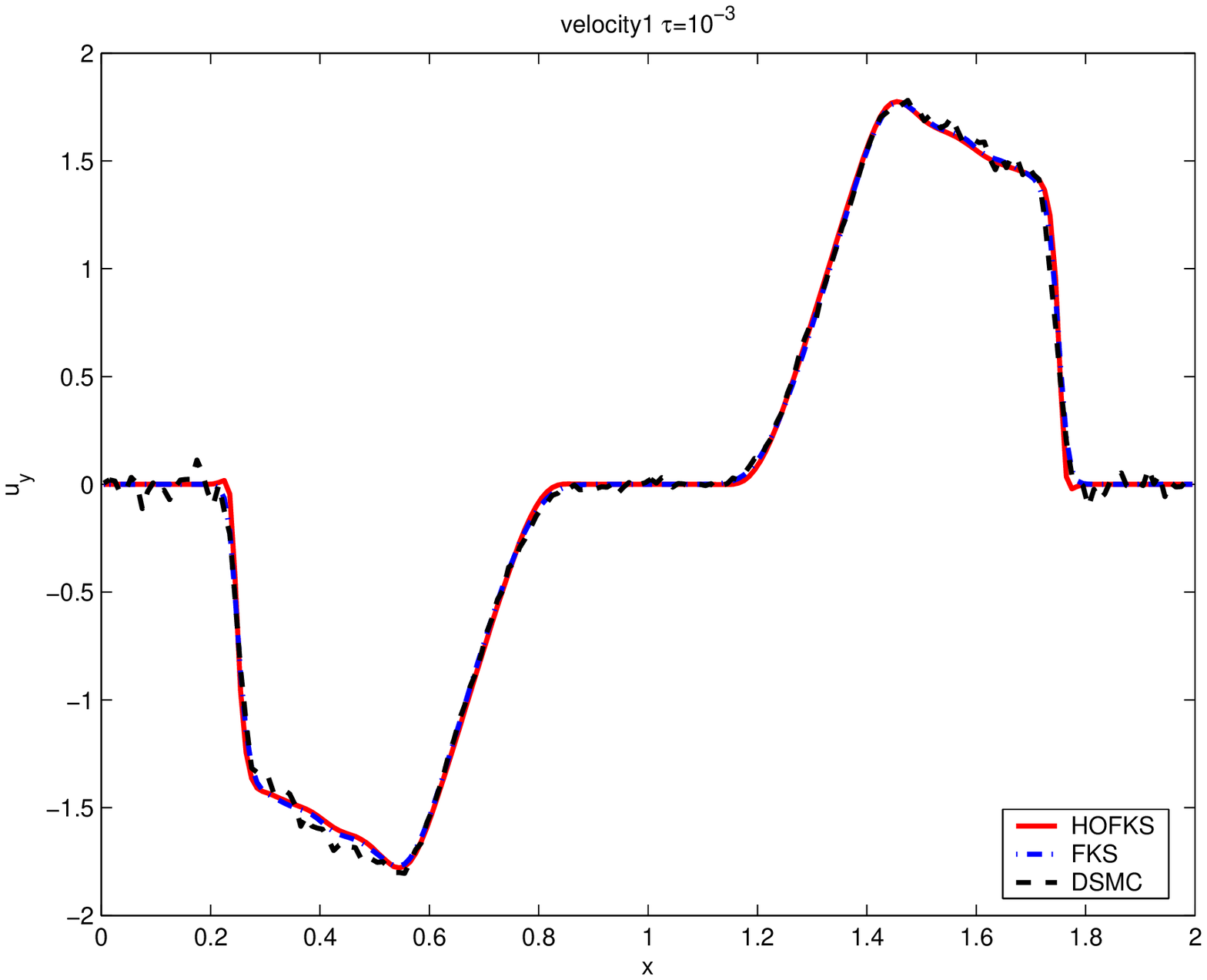}
\includegraphics[scale=0.4]{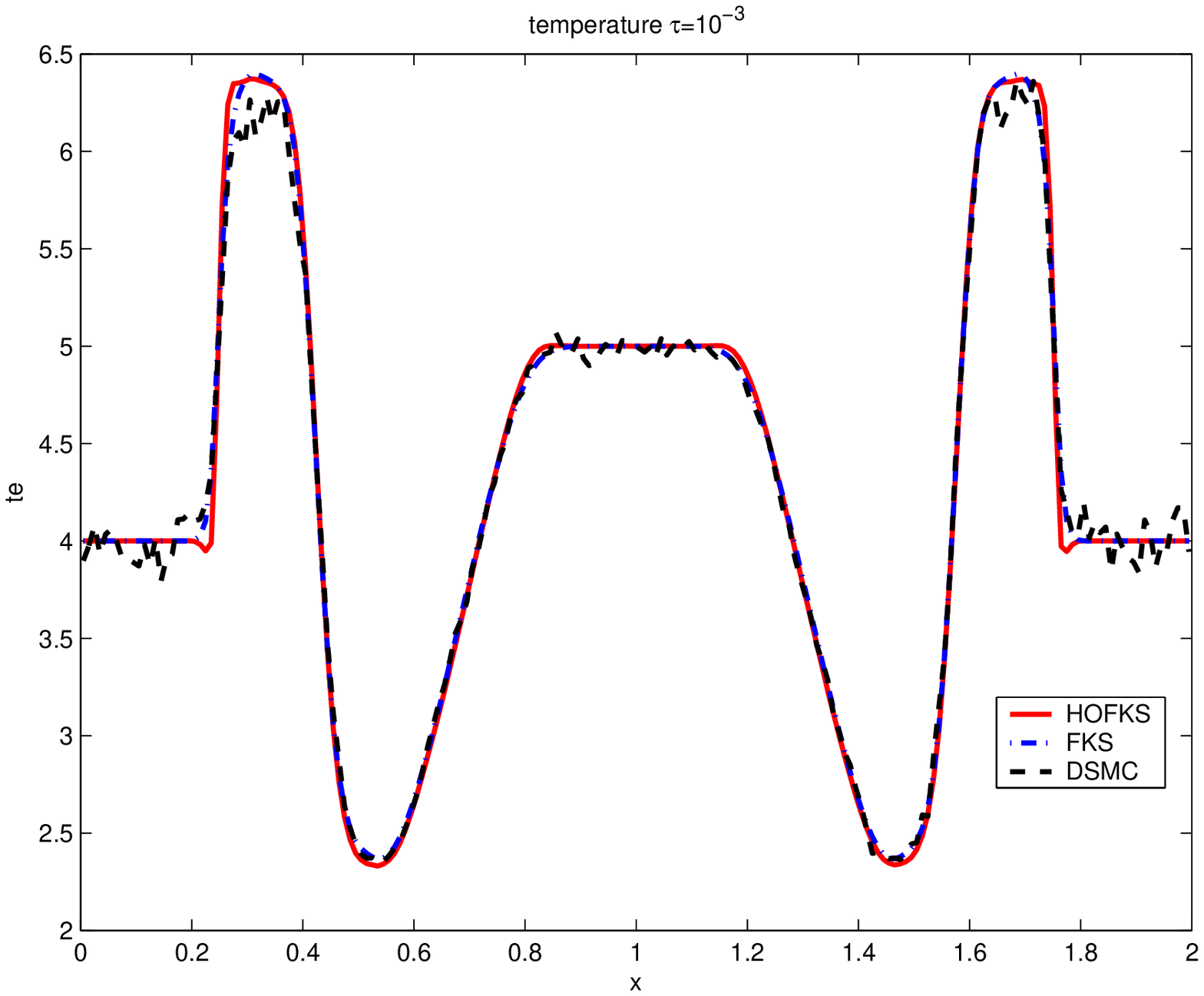}\\
\includegraphics[scale=0.45]{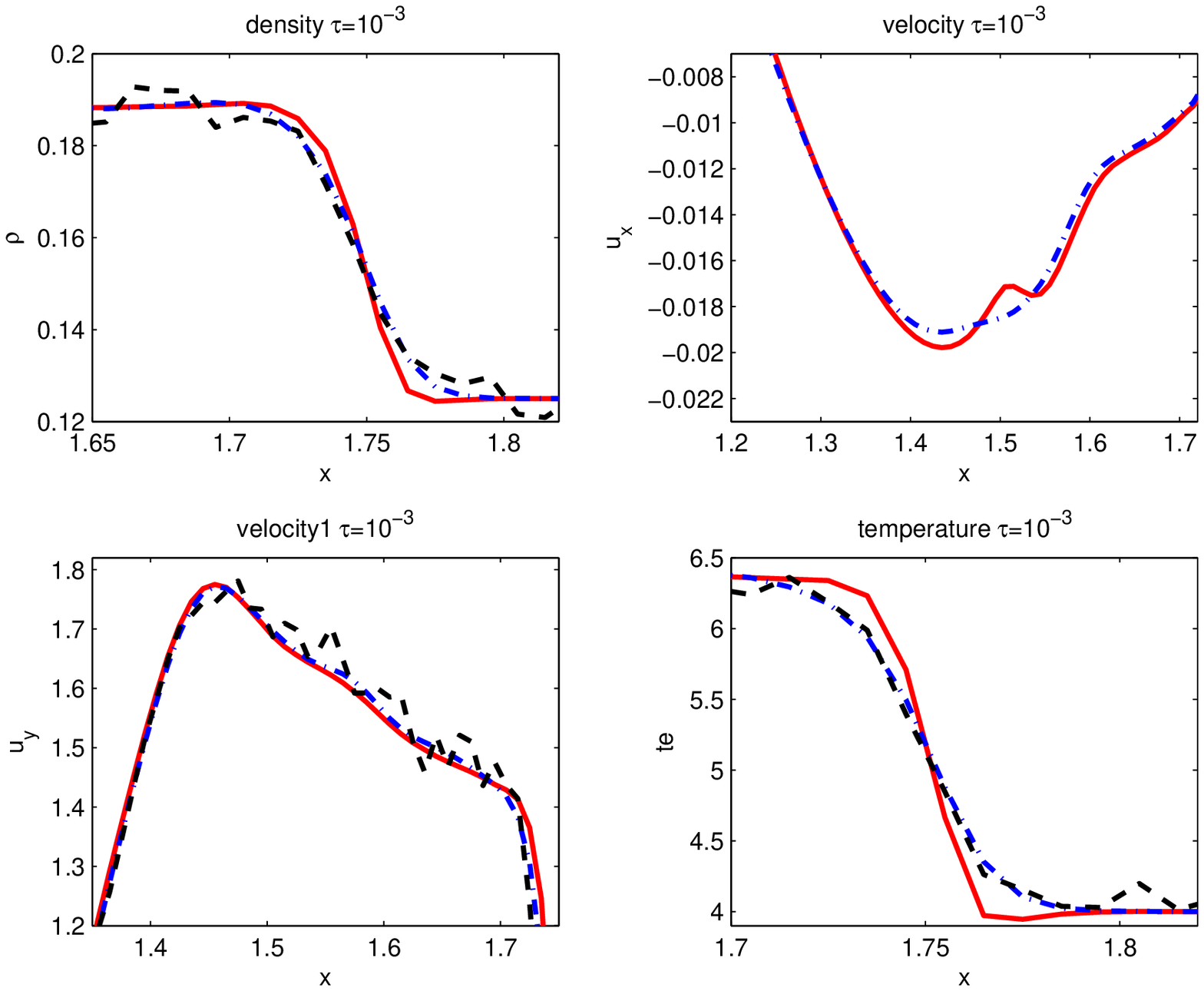}
\caption{2D Sod test: solution at
$t_{\text{final}}=0.07$ and $x=1$ for $\tau=10^{-3}$.
Top-Middle: Density (top left),
velocity in the $x$-direction (top right), velocity in the
$y$-direction (middle left) and temperature (middle right).
Bottom: magnification of the solution.
HOFKS method continuous line, FKS method dash dotted line, DSMC dotted line.
(DSMC results are not plotted for the y-component of velocity because the
noise induced by the method produces oscillations the amplitude of which are
far greater than the scale of the figure.)
} \label{sod2D3}
\end{center}
\end{figure}

In Figure~\ref{sod2D1}, we report the profiles fixing $x=1$ for respectively the
density, the mean velocity in the $x$-direction and in the $y$-direction
and the temperature using a $200\times 200$ mesh for $\tau=10^{-4}$. In
Figure~\ref{sod2D3}, we report the same profiles for the same spatial position,
\textit{i.e.} $x=1$, for the same macroscopic
quantities but for a larger value of the Knudsen number: $\tau=10^{-3}$. In this latter case, for the velocity in the $x$-direction, we did not report the solution for the DSMC method because the number of particles employed does not permit to compute the solution with sufficient precision. Observe, in fact, that in this test case, the final value of the $x$-velocity is of the order of $10^{-3}$, which is a value that due to the statistical fluctuations is very difficult to capture with DSMC methods.

Moreover in Figure~\ref{sod2D1} bottom panels, we report some magnifications of the same profiles which permits to better appreciate the differences between the methods. We observe that, as in the 1D case, the accuracy of the FKS method lies between the first and the second order spatial accuracy for small $\tau$. On the other hand, we cannot see the differences between the HOFKS scheme and the second order MUSCL scheme for $\tau=10^{-4}$ even with the magnifications. This is because the same MUSCL scheme is used for solving the compressible Euler equations and for constructing the HOFKS scheme. In the case $\tau=10^{-3}$, the HOFKS method, as in the 1D case, gives sharper solutions with respect the FKS method. This is also the case for the DSMC method which exhibits a larger numerical diffusion with respect to our kinetic scheme. For larger $\tau$, the HOFKS, the 
%=== E N D   F I G ================================================
FKS and the DSMC almost collapse on the same line, if the number of particles is chosen sufficiently large, and thus we did not report any result.

%=== B E G I N   F I G ================================================
%=== B E G I N   F I G ================================================

To conclude this test case study, we report, in table~\ref{tab:CPU_Sod2D1}, the CPU time $T$, as well as the CPU time per time cycle $T_{\text{cycle}}$, the CPU time per cycle per cell $T_{\text{cell}}$ and the number of cycles needed to perform the computation for different meshes in space for respectively the FKS, the HOFKS and the DSMC method when the Knudsen number is $\tau=10^{-3}$. For the HOFKS and the FKS schemes the meshes are fixed in velocity and
made of $20^{2}$ points. For the DSMC method we choose the same number of particles which has been used to produce the results shown in the previous Figures, \textit{i.e.} $100$ particles in average per cell and $100$ realizations. We computed the computational effort for these three schemes repeating the same test using different meshes in space ranging from $N_x=N_y=25$ to $N_x=N_y=200$. The results of this analysis can be summarized as follows.
The HOFKS method is around $1.5$ times more expensive than the FKS method for all the cases studied. This new scheme is however still very efficient: around $11.5$ minutes to compute the solution on a $200\times 200$ spacial mesh for a kinetic equation on a mono-processor laptop. For comparison the DSMC method, which is also known to be a fast method, requires around $505$ minutes. The ratio between the two methods is about $44$. Moreover DSMC gives less accurate solutions as seen on the $y$-component of velocity which can not be accurately captured.

%=== B E G I N   T A B L E ================================================
\begin{table}[h!]
  \begin{tabular}{|c|c|c|c|c|c|}
    \hline
    \textbf{Cell x \#}  $N_c$ & \textbf{\# Deg. freedom}  $N_{tot}$ & \textbf{Cycle} & \textbf{Time} & \textbf{Time/cycle} & \textbf{Time/cell} \\
    $N_x\times N_y$ & $N_x\times N_y\times N_{v_x}^2$ &  $N_{\text{cycle}}$  &  $\text{T}$ (s)  & $T_{\text{cycle}}$ (s) & $T_{\text{cell}}$ (s) \\
    \hline
    \hline
    \multicolumn{6}{|c|}{\textbf{FKS scheme}} \\
    \hline
    $25 \times 25$ & $25 \times 25\times 20^2$ & $13$   & $\sim 1.5$s   & $0.12$  & $1.9\times 10^{-4}$ \\
    $=625$ & $=250000$ & & & & \\
    $50 \times 50$ & $50 \times 50\times 20^{2}$ & $25$   & $6$s   & $0.24$  & $1.04\times 10^{-4}$ \\
    $=2500$ &$=10^6$ & &  & & \\
    $100 \times 100$ & $100 \times 100 \times 20^{2}$  & $50$  & $50$s   & $1$ & $1.0\times 10^{-4}$ \\
    $=10000$ &$=4 \times 10^6$ & &  &  & \\
    $200 \times 200$ & $200 \times 200\times 20^{2}$ & $100$  & $440$s   & $4.4$ & $1.1\times 10^{-4}$ \\
    $=40000$ &$=16 \times 10^6$ & &  &  & \\
    \hline
    \hline
    \multicolumn{6}{|c|}{\textbf{HOFKS scheme}} \\
    \hline
    $25 \times 25$ & $25 \times 25\times 20^2$ & $13$   & $\sim 2$s   & $0.15$  & $\sim 2.0\times 10^{-4}$ \\
    $=625$ & $=250000$ & & & & \\
    $50 \times 50$ & $50 \times 50\times 20^{2}$ & $25$   & $9$s   & $0.36$  & $1.44\times 10^{-4}$ \\
    $=2500$ &$=10^6$ & &  & & \\
    $100 \times 100$ & $100 \times 100 \times 20^{2}$  & $50$  & $77$s   & $1.54$ & $1.54\times 10^{-4}$ \\
    $=10000$ &$=4 \times 10^6$ & &  &  & \\
    $200 \times 200$ & $200 \times 200\times 20^{2}$ & $100$  & $690$s   & $6.9$ & $1.72\times 10^{-4}$ \\
    $=40000$ &$=16 \times 10^6$ & &  &  & \\
    \hline
    \hline
    \multicolumn{6}{|c|}{\textbf{DSMC scheme}} \\
    \hline
     & $N_c \times N_{average \times cell} $  & &  & & \\
    $25 \times 25$ &  $25 \times 25 \times 100^2$ & $11$   & $73$s   & $6.63$  & $0.0106$ \\
    $=625$ &$=6.25\times 10^{6}$ & & & & \\
    $50 \times 50$ & $50 \times 50 \times 100^2\times 50 \times 50$  & $22$   & $540$s   & $24.54$  & $0.0098$ \\
    $=2500$ & $=2.5\times 10^{7}$& &  & & \\
    $100 \times 100$ & $100 \times 100\times 100^2$  & $45$  & $3700$s   & $82.22$ & $0.0082$ \\
    $=10000$ & $=10^{8}$ & & $\sim 61$mn &  & \\
    $200 \times 200$ & $200 \times 200\times 100^2$ & $90$  & $30300$s   & $336.66$ & $0.0084$ \\
    $=40000$ & $=4\times 10^{8}$ & &$\sim 505$mn  &  & \\
    \hline
\end{tabular}
\caption{ \label{tab:CPU_Sod2D1}
    2D Sod shock tube. Computational effort for the FKS, HOFKS and DSMC schemes for $\tau=10^{-3}$.
    The time per cycle is obtained by $T_{\text{cycle}} = \text{T}/N_{\text{cycle}}$ and
    the time per cycle per cell by $T_{\text{cell}} = \text{T}/N_{\text{cycle}}/N_c$.
}
\end{table}
%=== E N D   T A B L E ================================================

%\clearpage
\subsection{2D Implosion problem}
\label{subsec_impl2D}

Finally we consider a 2D/2D implosion problem on the square $[0,2]\times[0,2]$ discretized with $100^{2}$ points.
As in the previous test, the velocity space is a square but with larger bounds, \textit{i.e.} $-20$ and $20$,
which means $[-20,20]^{2}$, discretized with $N_v=30$ points in
each direction which gives $40^{2}$ points. The domain is divided into two parts,
a disk centered at point $(1,1)$ of radius $R_d=0.2$ is filled with a
gas with density $\rho_{L}=0.125$, mean velocity $u_{L}=0$ and
temperature $T_{L}=4$, whereas the gas in the rest of the domain
is initiated with $\rho_{R}=1$, $T_{R}=4$, velocity in the $x$-direction $u_x=1$ for $x\in [0,1]$ and $u_x=-1$ for $x\in [1,2]$ while the velocity in the $y$-direction is initiated with $u_y=1$ for $y\in [0,1]$ and $u_y=-1$ for $y\in [1,2]$.
The final time is $t_{\text{final}}=0.07$.

%---- BEG FIG ------------
\begin{figure}[h!]
\begin{center}
\includegraphics[scale=0.365]{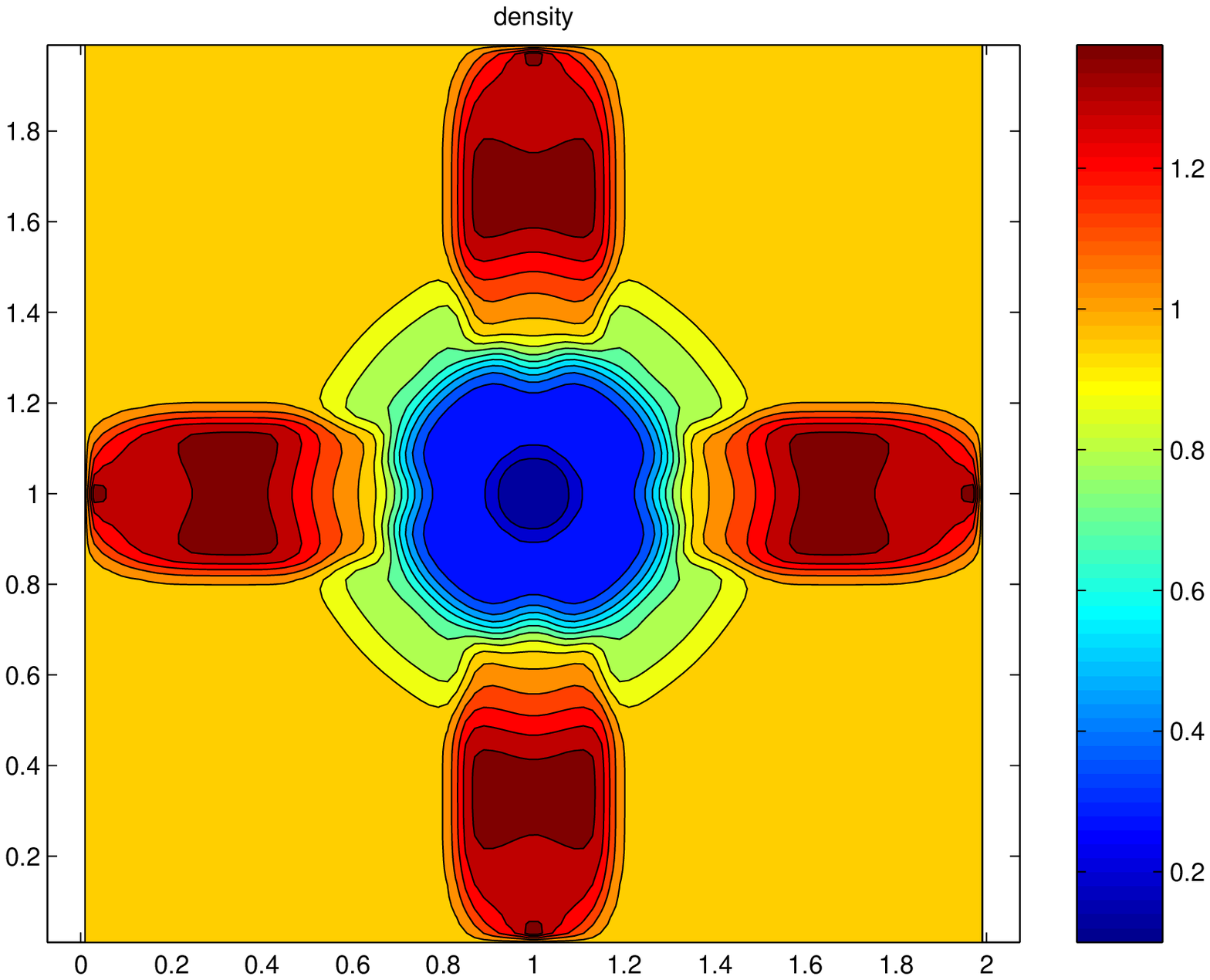}
\includegraphics[scale=0.365]{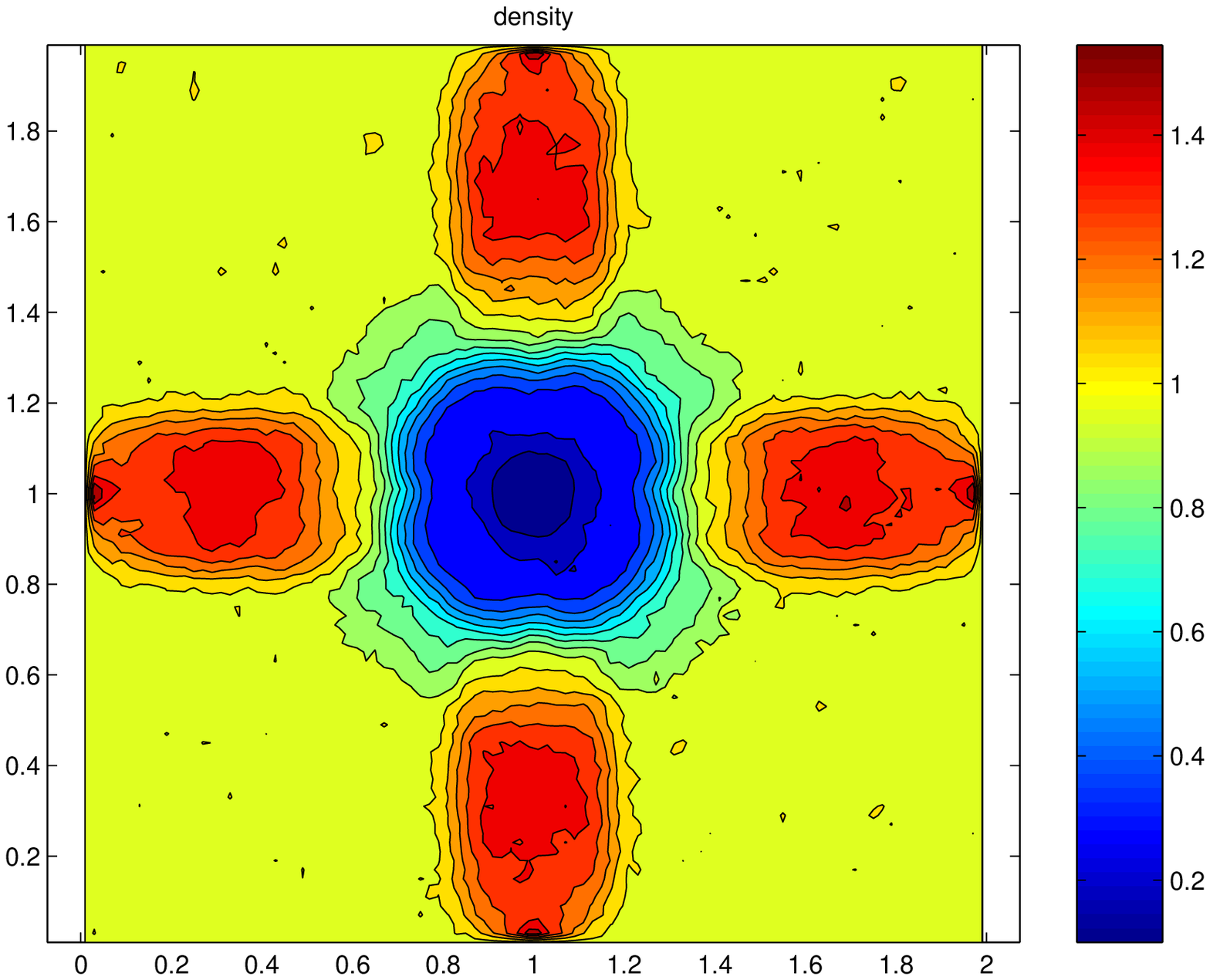}
\includegraphics[scale=0.365]{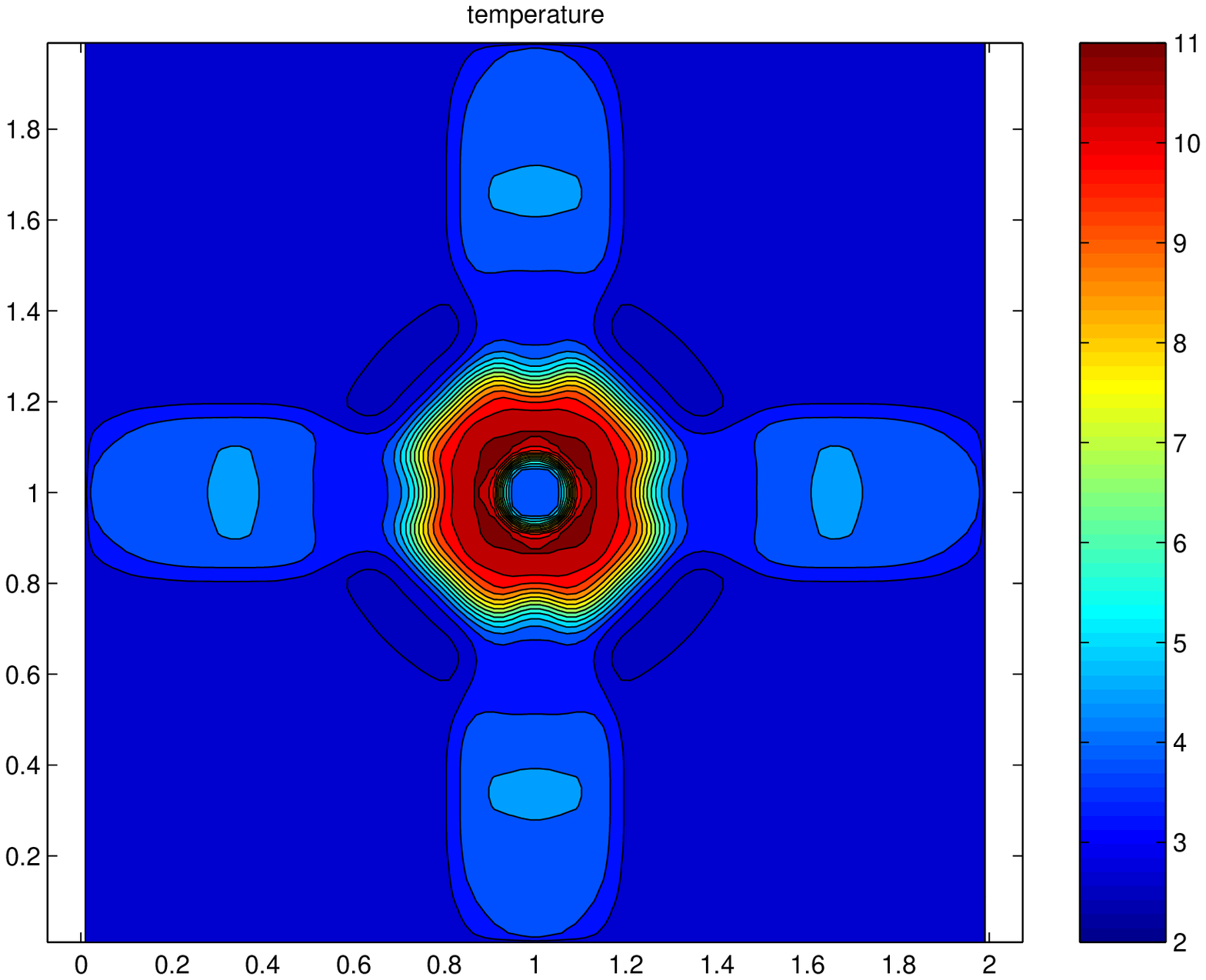}
\includegraphics[scale=0.365]{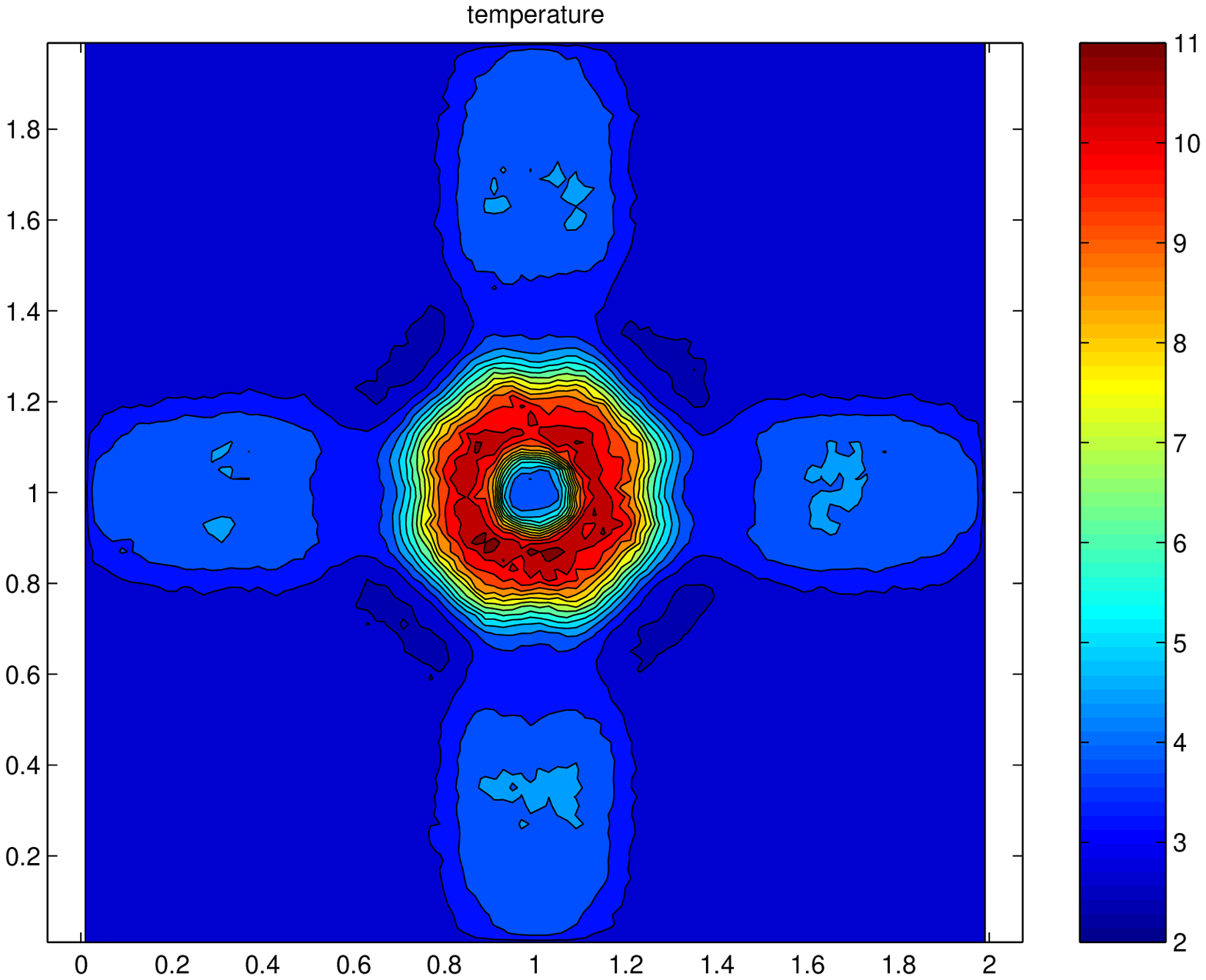}
\includegraphics[scale=0.365]{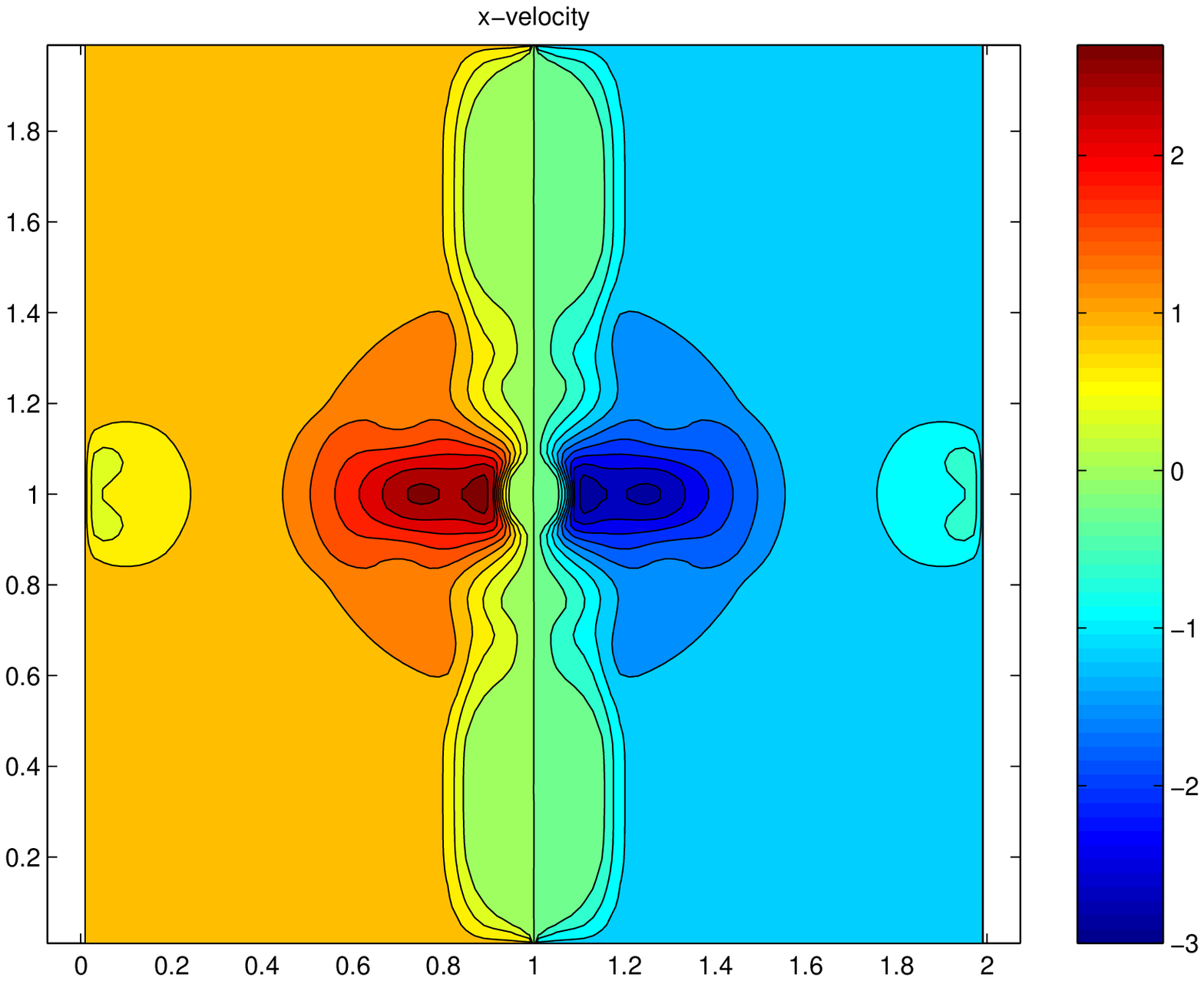}
\includegraphics[scale=0.365]{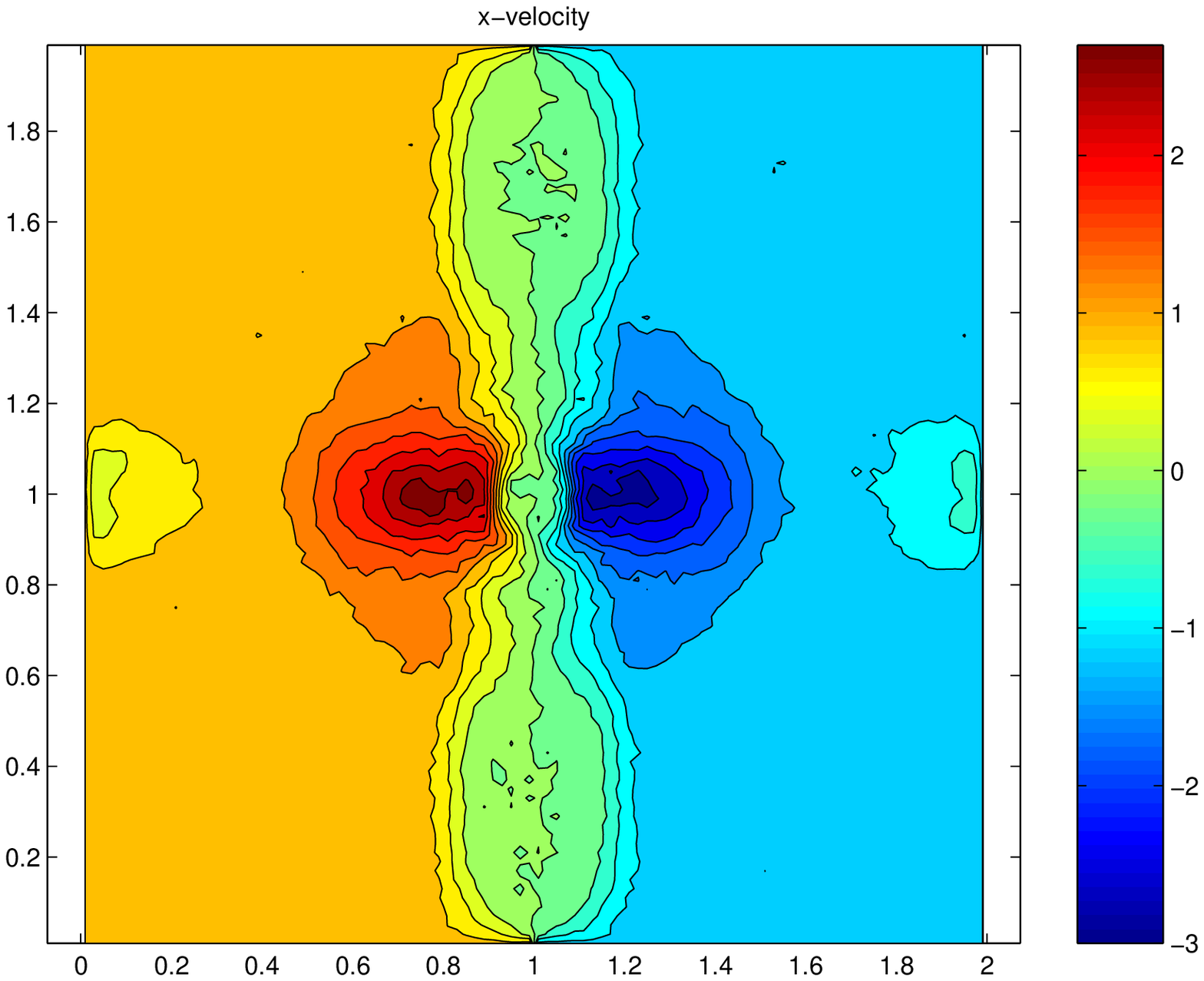}
\includegraphics[scale=0.365]{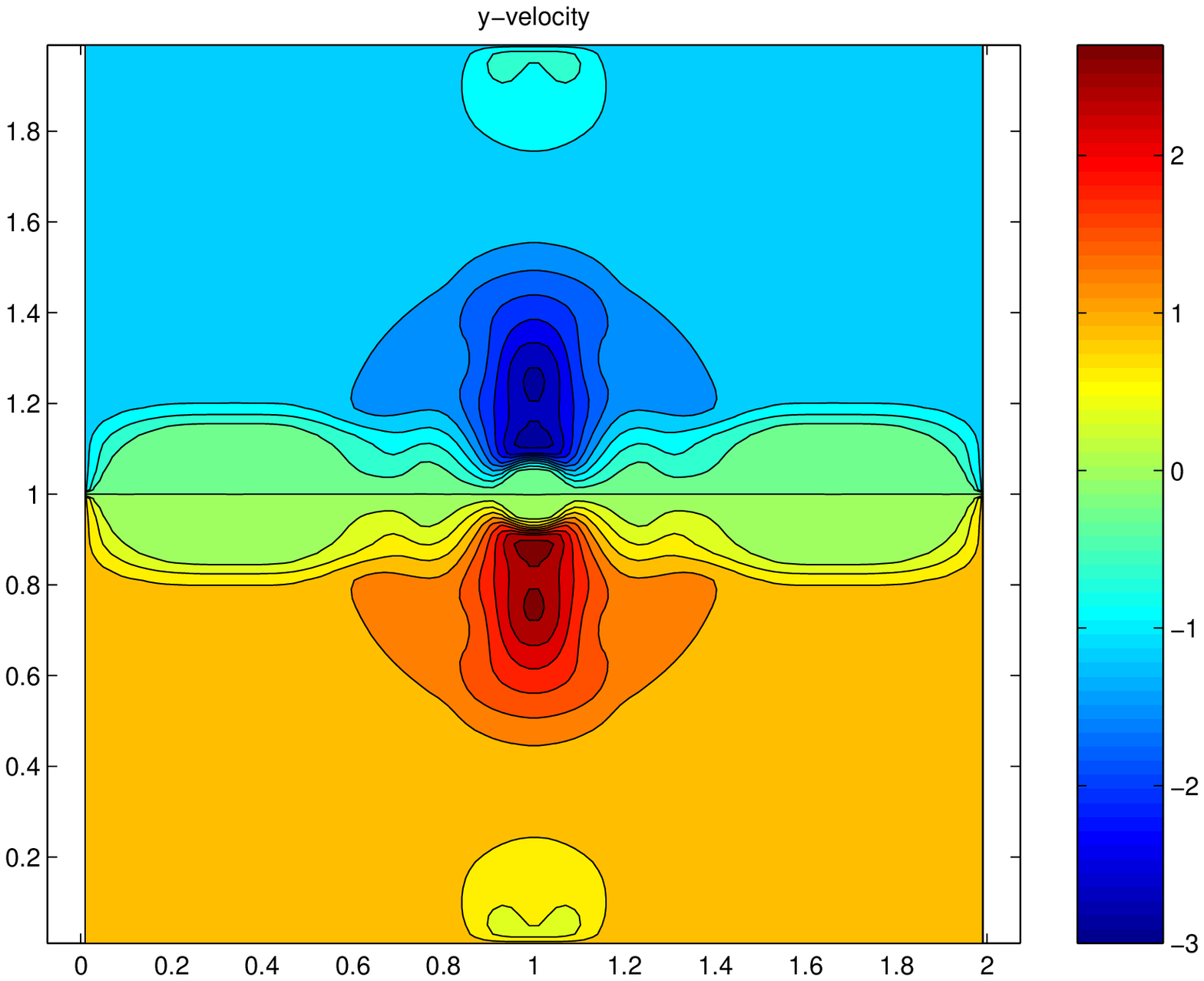}
\includegraphics[scale=0.365]{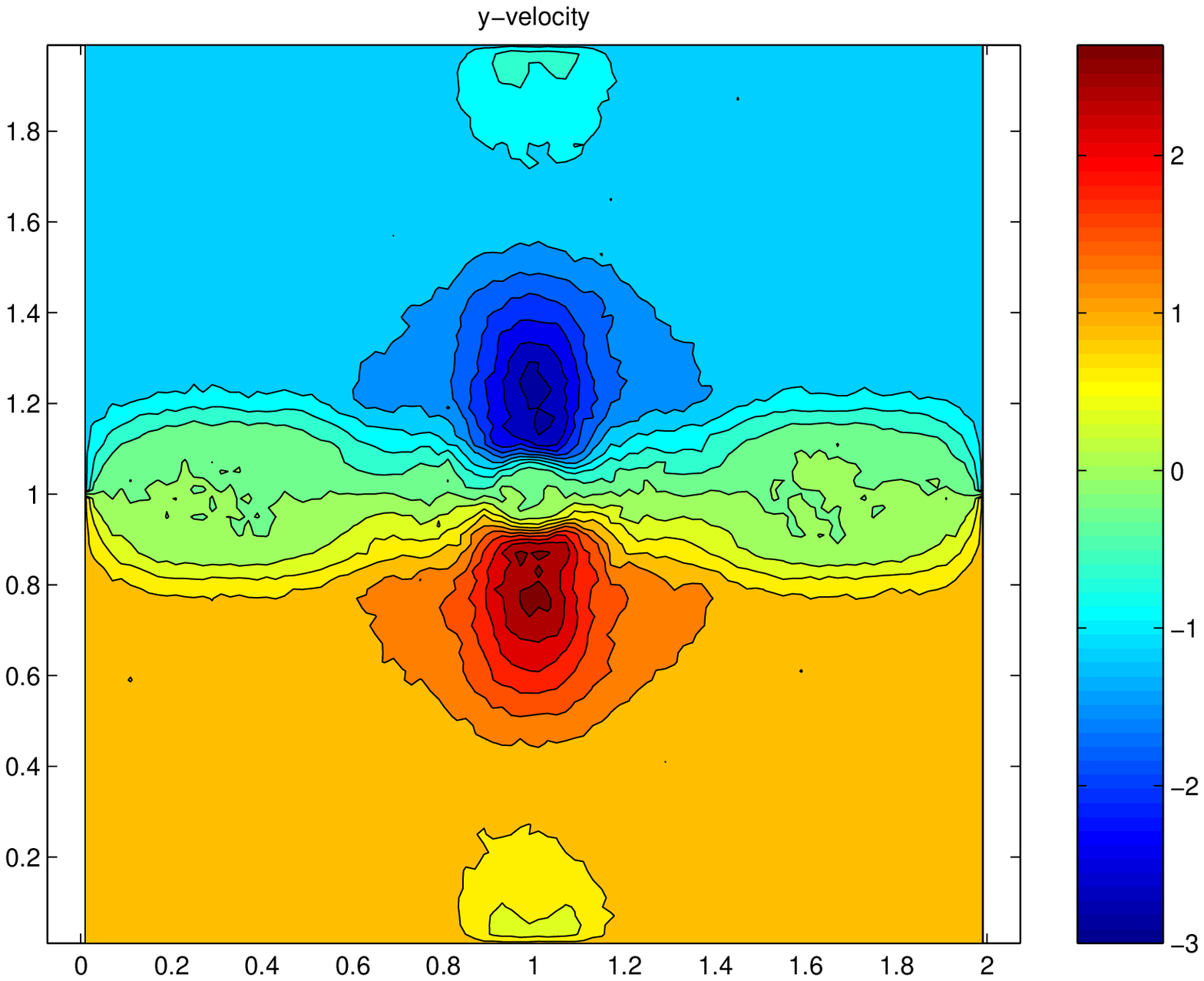}
\caption{2D Implosion test at
$t_{\text{final}}=0.07$ with $\tau=10^{-3}$. HOFKS scheme (left), DSMC scheme (right).
From top to bottom: density, temperature, $x$-velocity, $y$-velocity.
} \label{Implosion}
\end{center}
\end{figure}
%---- END FIG ------------

%=== B E G I N   T A B L E ================================================
\begin{table}[h]
  \begin{tabular}{|c|c|c|c|c|c|}
    \hline
    \textbf{Cell x \#}  $N_c$ & \textbf{\# Deg. freedom} & \textbf{Cycle} & \textbf{Time} & \textbf{Time/cycle} & \textbf{Time/cell} \\
    $N_x\times N_y$ &  &  $N_{\text{cycle}}$  &  $\text{T}$ (s)  & $T_{\text{cycle}}$ (s) & $T_{\text{cell}}$ (s) \\
    \hline
    \hline
    \multicolumn{6}{|c|}{\textbf{HOFKS scheme}} \\
    \hline
    & $N_x \times N_y \times  N_{v_x}^2 $ &  &   & &  \\
    $25 \times 25$ & $25 \times 25\times 30^2$ & $21$   & $\sim 3$s   & $0.14$  & $\sim 2.22\times 10^{-4}$ \\
    $=625$ & $=562500$ & & & & \\
    $50 \times 50$ & $50 \times 50\times 30^{2}$ & $43$   & $22$s   & $0.51$  & $2.04\times 10^{-4}$ \\
    $=2500$ &$=2.25 \times 10^6$ & &  & & \\
    $100 \times 100$ & $100 \times 100 \times 30^{2}$  & $85$  & $180$s   & $2.12$ & $2.08\times 10^{-4}$ \\
    $=10000$ &$=9 \times 10^6$ & &  &  & \\
    $200 \times 200$ & $200 \times 200\times 30^{2}$ & $170$  & $1350$s   & $7.94$ & $1.98\times 10^{-4}$ \\
    $=40000$ &$=36 \times 10^6$ & & $\sim 22.5$mn &  & \\
    \hline
    \hline
    \multicolumn{6}{|c|}{\textbf{DSMC scheme}} \\
    \hline
    & $N_x \times N_y \times N_{average \times cell}$ &  &   & &  \\
    $25 \times 25$ &  $ 25 \times 25 \times 50\times 200$ & $9$   & $74$s   & $8.22$  & $0.0132$ \\
    $=625$ &$=6.25\times 10^{6}$ & & & & \\
    $50 \times 50$ & $50 \times 50 \times 50\times 200$  & $20$   & $620$s   & $31$  & $0.0124$ \\
    $=2500$ & $=2.5\times 10^{7}$& &  & & \\
    $100 \times 100$ & $100 \times 100 \times 50\times 200$  & $40$  & $3842$s   & $96.05$ & $0.0096$ \\
    $=10000$ & $=10^{8}$ & &$\sim 64$mn  &  & \\
    $200 \times 200$ & $200 \times 200 \times 50\times 200$ & $81$  & $30600$s   & $377.77$ & $0.0094$ \\
    $=40000$ & $=4\times 10^{8}$ & & $\sim 8.5$h  &  & \\
    \hline
\end{tabular}
\caption{ \label{tab:CPU_Imp2D1}
    2D Implosion test. Computational effort for the HOFKS and the DSMC scheme.
    The time per cycle is obtained by $T_{\text{cycle}} = \text{T}/N_{\text{cycle}}$ and
    the time per cycle per cell by $T_{\text{cell}} = \text{T}/N_{\text{cycle}}/N_c$.
}
\end{table}
%=== E N D   T A B L E ====================================================

We report the results for the Knudsen number equal to $\tau=10^{-3}$ comparing our scheme to the DSMC method for the Boltzmann-BGK equation.
For this latter, we employed on average $200$ particles per cell and the solution is averaged over $50$ realizations.\\
In Figure~\ref{Implosion} we report the isolines of density, $x$-velocity, $y$-velocity and temperature,
for the HOFKS method on the left and the DSMC method on the right.
We observe that the two methods furnish the same results except for the statistical noise of the DSMC method. On the other hand, the computational costs of the two approaches are still very different. We report, as for the 2D Sod test, in table~\ref{tab:CPU_Imp2D1}, the CPU time $T$, as well as the CPU time per time cycle $T_{\text{cycle}}$, the CPU time per cycle per cell $T_{\text{cell}}$ and the number of cycles needed to perform the computation for different meshes in space. For the DSMC method we choose the same number of particles which has been used to produce the results shown in the figures, \textit{i.e.} $200$ particles in average per cell and $50$ realizations. The results of this last test can be summarized as follows: The HOFKS method takes around $22$ minutes for computing the solution on a $200\times 200$ mesh for a kinetic equation on a mono-processor laptop. The augmentation of the computational cost with respect to the 2D Sod test is essentially due to: first the augmentation of the mesh points in which the velocity space is discretized and second the reduction of the time step. The reduction of the time step is caused by the macroscopic solver of the compressible Euler equations which needs smaller time steps to ensure stability. For comparison the DSMC method, which still furnishes fluctuating and somewhat less accurate solutions, requires around $510$ minutes. The ratio of CPU time is around $22$ in favor of the HOFKS for a better overall accuracy.

%%%%%%%%%%%%%%%%%%%%%%%%%%%%%%%%%%%%%%%%%%
%%%%%%%%%%%%%%%%%%%%%%%%%%%%%%%%%%%%%%%%%%
%%%%%%%%%%%%%%%%%%%%%%%%%%%%%%%%%%%%%%%%%%
\clearpage
\section{Conclusions} \label{sec_conclu}

In this work we have presented an high order in space extension of a new super
efficient numerical method for solving kinetic equations. The method
is based on a splitting between the collision and the transport
terms. The collision part is solved on a grid while the transport
linear part is solved exactly by following the characteristics
backward in time. The key point is that, conversely to
semi-Lagrangian methods, we do not need to reconstruct the
distribution function at each time step. This permits to tremendously reduce the computational cost with respect other existing methods for kinetic equations. In this paper, in order to solve the limitations in term of spatial accuracy close to the thermodynamical equilibrium of the original Fast Kinetic Scheme, we coupled the solution of the FKS method with the solution of the compressible Euler equations. Then, we matched the moments obtained from the solution of the macroscopic equations with the moments obtained from the solution of the equilibrium part of the kinetic equation. Finally, we recovered the solution as a convex combination of the two contributions: the macroscopic and the microscopic parts. This improvement permits to preserve the desired accuracy in space for all the different regimes studied.

The numerical results show that the HOFKS method performs as the FKS method for large values of the Knudsen number and as a high order shock capturing scheme for small Knudsen numbers. Moreover, the method requires a small computational effort. Numerical experiments have shown that the computational cost of the new method is around $1.5$ times larger than the previous FKS method for a clear gain in accuracy when reached fluid regimes. Most importantly, we showed that this new class of fast kinetic schemes is more accurate and around $25-65$ times faster than DSMC methods which are known to be efficient schemes. This important result opens the gate to extensive realistic numerical simulations of far from equilibrium physical models.

In this work, we only focused on the spatial accuracy of the method and not on the time accuracy, we remind to future works for the development of schemes which are both accurate in time and in space. We also would like to extend the method to non uniform meshes and different discretizations of the velocity space. Finally, we want to apply this method to other kinetic equations as the Boltzmann or the Vlasov equation.

%========================================================================
% BIBLIOGRAPHY
%

%
% END DOCUMENT
%
%=============================================================================

\end{document}